\input amstex
\magnification=\magstephalf
\documentstyle{amsppt}
\input pictex
\catcode`\@=11
\def\logo@{}
\catcode`\@=\active

\NoBlackBoxes
\TagsAsMath
\TagsOnRight
\parindent=.5truein
\parskip=3pt
\hsize=6.5truein
\vsize=9.0truein
\raggedbottom
\loadmsam
\loadbold
\loadmsbm
\define     \x{\times}
\let        \< = \langle
\let        \> = \rangle
\let        \| = \Vert
\define     \a          {\alpha}
\redefine   \b          {\beta}
\redefine   \B          {\Cal B}
\redefine   \dl         {\delta}
\redefine   \D          {\Cal D}
\redefine   \Dl         {\Delta}
\redefine   \e          {\epsilon}
\define     \F          {\Cal F}
\define     \G          {\Gamma}
\redefine   \g          {\gamma}
\define     \br         {\bold R}
\define     \R         {\bold R}
\define     \Z          {\bold Z}
\define     \s          {\sigma}
\define     \th         {\theta}
\redefine   \O          {\Omega}
\redefine   \o          {\omega}
\redefine  \w            {\omega}
\redefine   \t          {\tau}
\define     \ub         {\subset}
\define     \cl         {\Cal L}
\define     \E          {\Cal E}
\define     \A          {\Cal A}
\redefine   \i          {\infty}
\define     \pt         {\partial}
\define     \bl         { L}
\define     \ube         {\subseteq}

\topmatter

\title
The Stable Manifold Theorem for Semilinear Stochastic Evolution Equations and 
Stochastic 
Partial 
Differential Equations\\
\medskip
I: The Stochastic Semiflow 
 \endtitle

\author
Salah-Eldin A\. Mohammed{\footnote "$\sp {*}$" {The research of
this author is supported in part by NSF Grants DMS-9703852,
DMS-9975462 and DMS-0203368.\hfil \hfil}}, Tusheng Zhang
{\footnote"$\sp
{**}$"{The research of this author is supported in part by EPSRC Grant
GR/R91144.\hfil\hfil} and Huaizhong Zhao{\footnote"$\sp
{***}$"{The research of this author is supported in part by EPSRC Grants
GR/R69518 and GR/R93582.\hfil\hfil \newline
     AMS 1991 {\it subject classifications.\/} Primary 60H10, 60H20;
secondary 60H25. \newline {\it Key words and phrases.\/} Stochastic semiflow,
$C^k$ cocycle, stochastic evolution equation (see), stochastic partial differential 
equation (spde).
 \hfil\hfil} } }
\endauthor


\abstract
The main objective of this work is to characterize the pathwise local structure
of solutions of semilinear stochastic evolution equations (see's) and
   stochastic partial differential equations (spde's)
    near stationary solutions.  Such characterization is realized through the 
long-term 
behavior of the solution field near
stationary points.  The analysis falls in two parts I, II.  In Part I (this paper),  
we prove a general 
existence and compactness theorem for $C^k$-cocycles  of semilinear see's and 
spde's. Our results cover a large class
of semilinear see's as well as certain semilinear spde's with non-Lipschitz terms 
such
as stochastic reaction diffusion equations and the stochastic Burgers equation with 
additive infinite-dimensional noise. In Part II of this work ([M-Z-Z]), 
we  establish a {\it local stable manifold theorem\/} for non-linear see's and 
spde's.
\endabstract

\endtopmatter

\rightheadtext{Stochastic Flows for SPDE's}

\leftheadtext{S.-E.A. Mohammed, T.S. Zhang  and H.Z. Zhao}

\document
\baselineskip=16truept



\subheading{1. Introduction}

  The construction of local stable and unstable manifolds near hyperbolic equilibria 
is a 
fundamental problem in deterministic
and stochastic dynamical systems. The significance of these invariant manifolds  
consists 
in a characterization of the local
behavior of the dynamical system in terms of long-time asymptotics of its 
trajectories 
near a stationary point.  In recent
years, it has been established 
that local stable/unstable manifolds exist  for finite-dimensional stochastic 
ordinary 
differential equations (sode's) ([M-S.2])
and stochastic systems with finite memory (sfde's)([M-S.1]).  On the other hand, 
existence 
of such manifolds  for stochastic
evolution equations (see's) and stochastic partial differential equations (spde's) 
has 
been an open problem since the early
nineties ([F-S], [B-F], [B-F.1]). 

The main objective of the present work is to establish the existence of local stable 
and 
unstable manifolds near stationary solutions of  semilinear
stochastic
evolution equations (see's) and stochastic partial differential equations (spde's).  
Our 
approach consists in the following
two major undertakings:

\item{$\bullet$} A construction of a sufficiently Fr\'echet differentiable cocycle 
for  
mild/weak trajectories of the see or the spde.
This is achieved in the see case by  a combination of a chaos-type expansion and 
suitable 
lifting techniques, and for
spde's 
by using stochastic variational representations and methods from deterministic 
pde's. 
Part I of this work is devoted to  detailing the construction of the cocycle.

\item{$\bullet$} The application of classical non-linear ergodic theory techniques 
developed by Oseledec [O] and Ruelle [Ru.2]
in order to study the local structure  of the above cocycle in a neighborhood of a 
hyperbolic stationary point.  Stationary
points
correspond to stationary solutions and hyperbolicity is characterized via the 
Lyapunov 
spectrum of the linearized cocycle
along the stationary trajectory. This is the subject  of Part II of this work 
([M-Z-Z]).

In [F-S], Flandoli and Schauml\"offel established the existence of a
random evolution operator and its Lyapunov spectrum for a linear stochastic heat 
equation with finite-dimensional noise, on 
a bounded Euclidean domain. For linear see's with finite-dimensional noise, 
a stochastic semi-flow (i.e. random evolution operator) was obtained in [B-F].
Subsequent work on the dynamics of non-linear 
spde's has
focused mainly on the question of existence of continuous semiflows 
and the existence and uniqueness of invariant measures.

The  problem of existence of semiflows for see's and spde's is a non-trivial one, 
mainly
due to the well-established fact that finite-dimensional methods for constructing 
(even 
continuous)  stochastic
flows  break down in the infinite-dimensional setting  of spde's and see's. In 
particular, 
Kolmogorov's 
continuity theorem fails for random fields parametrized by infinite-dimensional 
Hilbert 
spaces
(cf.  [Mo.1], pp. 144-149, [Sk], [Mo.2], [F.1], [F.2], [D-Z.1], pp. 246-248). 

As indicated above, the existence of a smooth semiflow is a necessary tool for 
constructing the stable and unstable
manifolds
near a hyperbolic stationary random point, ala  work of Ruelle ([Ru.1], [Ru.2]). In 
this 
article, we show the existence of
smooth perfect
cocycles for mild solutions of semilinear see's in Hilbert space (Theorem 2.6). Our 
construction employs a ``chaos-type"
representation coupled with lifting and variational techniques using the linear 
terms of 
the see (Theorems 2.1-2.4). This
technique bypasses the need for Kolmogorov's continuity theorem and appears to be 
new.  
Applications to specific
classes of
spde's are given.
 In particular, we obtain smooth stochastic semiflows for semilinear spde's driven 
by 
cylindrical  Brownian motion with a covariance Hilbert space $K$ (Theorem 3.5). 
In these applications, it turns out that in addition to smoothness of the non-linear 
terms, one requires some level of
dissipativity or Lipschitz  continuity in order to guarantee the existence of smooth 
globally defined 
semiflows.  Specific examples of spde's include semilinear parabolic spde's with 
Lipschitz nonlinearities (Theorem 3.5), stochastic reaction diffusion 
equations 
(Theorems 4.1, 4.2) and stochastic
Burgers equations with additive infinite-dimensional noise
(Theorem 4.3).
 
We begin by formulating the ideas of  a {\it stochastic semiflow} and a {\it 
cocycle} 
which are central to the analysis in this work.

Let $(\O,\F,P)$ be a  probability space. Denote by $\bar \F$ the
$P$-completion of $\F$, and let \newline
$(\O,\bar \F,(\F_t)_{t\ge 0},P)$ be  a complete filtered probability space
satisfying the usual conditions ([Pr]).

Denote $\Delta:=\{(s,t)\in \R^2: 0 \leq s \leq t \}$, and $\R^+ := [0, \infty )$. 
For a topological
space $E$, let $\B(E)$ denote its  Borel $\s$-algebra.

Let $k$ be a positive integer and $0 < \e \le 1$.
If $E,N$ are real Banach spaces,
we will denote by $L^{(k)} (E,N)$ the Banach space of all $k$-multilinear maps
$A:E^k \to N$ with the uniform norm $\|A\|:=\sup \{|A(v_1,v_2,\cdots,v_k)|:
v_i \in E, |v_i| \leq 1, i=1,\cdots,k  \}$. Suppose $U \subseteq E$ is
an open set. A map $f:U \to N$ is said to be
{\it of class $C^{k,\e}$} if it is $C^k$ and if $D^{(k)} f: U \to L^{(k)} (E,N)$ is
$\e$-H\"older continuous on bounded sets in $U$.
A $C^{k,\e}$ map $f:U \to N$ is said to be {\it of class $C^{k,\e}_b$}
if all its derivatives $D^{(j)} f, 1 \leq j \leq k$, are globally bounded on $U$,
and $D^{(k)} f$ is $\e$-H\"older continuous on $U$.
A mapping
$\tilde f:[0,T] \times U \to N$ is {\it of class $C^{k,\e}$ in the second
variable uniformly with respect to the first} if for each $t \in [0,T]$,
$\tilde f(t,\cdot)$ is $C^{k,\e}$ on $U$, for every bounded set $U_0 \subseteq
U$ the
spatial partial derivatives $D^{(j)} \tilde f(t,x), j=1,\cdots, k,$ are uniformly
bounded
in $(t,x) \in [0,T]\times U_0$ and the
corresponding $\e$-H\"older constant of $D^{(k)} \tilde f(t,\cdot)|U_0$
is uniformly bounded in $t \in [0,T]$.

The following definitions are crucial to the developments in this article.

\definition{Definition 1.1} 

Let $E$ be a Banach space, $k$ a
non-negative integer and $\e \in (0,1]$.
A {\it stochastic} $C^{k,\e}$ {\it semiflow} on $E$
is a random field $V:\Dl \times E \times \O \to E$ satisfying the
following properties:

\item{(i)} $V$ is $(\B (\Dl) \otimes \B(E) \otimes \F,
         \B(E))$-measurable.

\item{(ii)} For each $\o \in \O$, the map
$\Dl \times E \ni (s,t,x) \mapsto V(s,t,x,\o) \in E $
      is continuous.

\item{(iii)} For fixed $(s,t,\o)\in \Dl \times \O$, the map
$E \ni x \mapsto X(s,t,x,\o) \in E$
     is $C^{k,\e}$.

\item{(iv)} If $0 \le r \leq  s \le t $, $\o \in \O$ and $x \in E$,
     then
     $$
         V(r,t,x,\o) = V(s,t,V(r,s,x,\o),\o).
$$

\item{(v)} For all $(s,x,\o) \in \R^+ \times E \times \O$, one has
$V(s,s,x,\o) = x.$
\enddefinition


\definition{Definition 1.2} 

Let $\theta :\R^+ \times \O \to \O$ be a
$P$-preserving semigroup on
the probability space $(\O,\F,P)$, $E$ a Banach space, $k$ a non-negative
integer and $\e \in (0,1]$. A $C^{k, \e}$ {\it perfect
cocycle\/} $(U, \theta)$  on $E$
is a $(\B (\R^+) \otimes \B(E) \otimes \F, \B(E))$-measurable random field
$U:\R^+ \times E \times \O \to E$
with the following properties:

\item{(i)} For each $\o\in \O$, the map
$\R^+ \times E \ni (t,x) \mapsto U(t,x,\o) \in E$
     is continuous; and for fixed $(t,\o)\in \R^+ \times \O$, the map
$E \ni x \mapsto U(t,x,\o) \in E$ is $C^{k,\e}$.

\item{(ii)} $U(t+s,\cdot,\o)=
U(t,\cdot,\theta (s,\o))\circ U(s,\cdot,\o)$
for all $s,t \in \R^+$ and all $\o \in \O$.

\item{(iii)} $U(0,x,\o)=x$ for all $x \in E, \o \in \O$.
\enddefinition


Note that a cocycle  $(U,\theta)$ corresponds to a one-parameter semigroup on
$E \times \O$. The following figure illustrates the cocycle
property. The vertical solid lines
represent random  copies of $E$ sampled according to
the probability measure $P$.

\newpage




   \centerline{\it The  Cocycle Property}
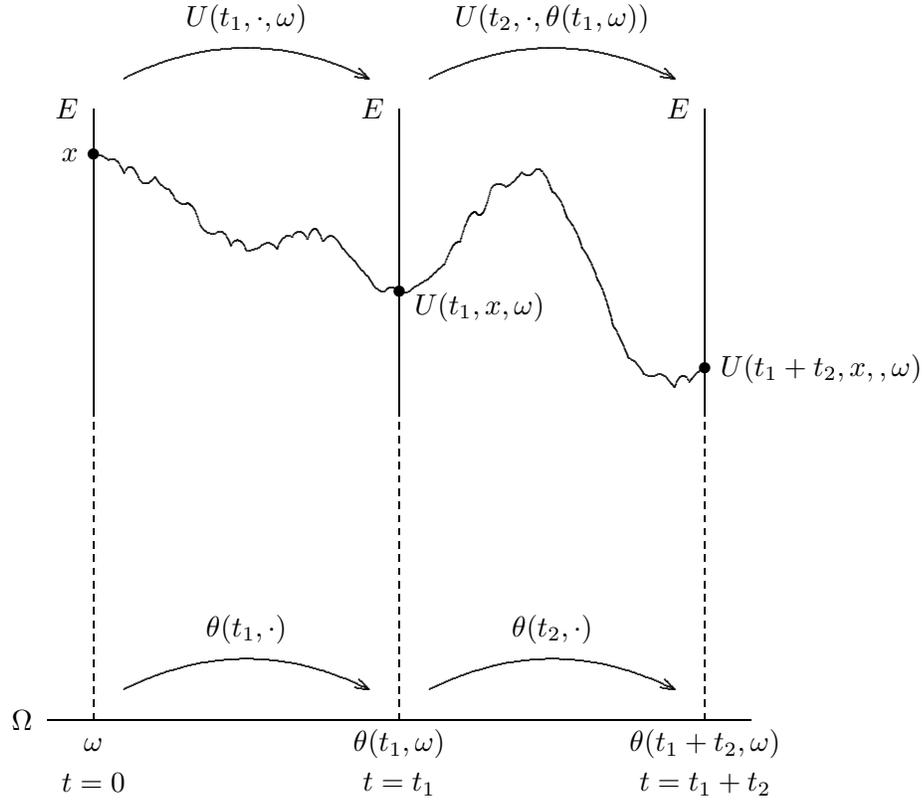
\centerline{
\beginpicture
\setcoordinatesystem units <0.8truein, 0.8truein>
\setplotarea x from -.8 to 5, y from -.8 to 5.0
%
    %
%
\putrule from -.3 0 to 4.3 0
%
%
\putrule from 0 2 to 0 4
\putrule from 2 2 to 2 4
\putrule from 4 2 to 4 4
\put {$E$} [rc] at -0.1 4.0
\put {$E$} [rc] at  1.9 4.0
\put {$E$} [rc] at  3.9 4.0
%
%
\setdashes <2.5truept>
\putrule from 0 0 to 0 2
\putrule from 2 0 to 2 2
\putrule from 4 0 to 4 2
\setsolid
\put {$\Omega$}                 [rc] at -0.4 0
\put {$\omega$}                 [cb] at 0 -.20
\put {$\theta(t_1,\omega)$}     [cb] at 2 -.25
\put {$\theta(t_1+t_2,\omega)$} [cb] at 4 -.25
\put {$t=0$}                    [ct] at 0 -.35
\put {$t=t_1$}                  [ct] at 2 -.35
\put {$t=t_1+t_2$}              [ct] at 4 -.35
%
%
\setquadratic    
\plot 0.2 4.2   1.0 4.4   1.8 4.2 /
\plot 2.2 4.2   3.0 4.4   3.8 4.2 /
\arrow <.07truein> [.4, .7] from 1.79 4.205 to 1.8 4.2
\arrow <.07truein> [.4, .7] from 3.79 4.205 to 3.8 4.2
\put {$U(t_1,\cdot,\omega)$}         [cb] at 1.0 4.5
\put {$U(t_2,\cdot,\theta(t_1,\omega))$} [cb] at 3.0 4.5
%
%
\plot 0.2 0.2   1.0 0.4   1.8 0.2 /
\arrow <.07truein> [.4, .7] from 1.79 0.205 to 1.8 0.2
\put {$\theta(t_1,\cdot)$}         [cb] at 1.0 0.5
\plot 2.2 0.2   3.0 0.4   3.8 0.2 /
\arrow <.07truein> [.4, .7] from 3.79 0.205 to 3.8 0.2
\put {$\theta(t_2,\cdot)$}         [cb] at 3.0 0.5
%
%
\put {$\bullet$}  [cc] at 0.0 3.7
\put {$x$} [rc] at -.1 3.7
\put {$\bullet$}  [cc] at 2.0 2.8
\put {$U(t_1,x,\omega)$}  [lt] at 2.1 2.8
\put {$\bullet$}  [cc] at 4.0 2.3
\put {$U(t_1+t_2,x,,\omega)$} [lc] at 4.1 2.3
\setquadratic
\plot
0.0 3.7036
0.05 3.69397594
0.10 3.66604728
0.15 3.66270941
0.20 3.58075773
0.25 3.61918764
0.30 3.53729456
0.35 3.51397384
0.40 3.55642091
0.45 3.49643116
0.50 3.47170000
0.55 3.39492741
0.60 3.37673181
0.65 3.35683621
0.70 3.23486362
0.75 3.18983706
0.80 3.18377953
0.85 3.19021404
0.90 3.10386363
0.95 3.14065127
1.00 3.07360000
1.05 3.07891442
1.10 3.10582548
1.15 3.12094574
1.20 3.08398777
1.25 3.16166412
1.30 3.16178734
1.35 3.19646998
1.40 3.14732460
1.45 3.21436375
1.50 3.13110001
1.55 3.16731492
1.60 3.13686627
1.65 3.06438080
1.70 3.01408527
1.75 2.95140646
1.80 2.89807113
1.85 2.82760606
1.90 2.80263800
1.95 2.83419373
2.00 2.80
2.05 2.79997586
2.10 2.83430943
2.15 2.86268106
2.20 2.89507114
2.25 2.94386004
2.30 2.98612813
2.35 3.09995580
2.40 3.12882341
2.45 3.27891136
2.50 3.32480000
2.55 3.32623162
2.60 3.47909602
2.65 3.51104498
2.70 3.48843020
2.75 3.52680340
2.80 3.58601635
2.85 3.56812075
2.90 3.60906833
2.95 3.57311084
3.00 3.46040001
3.05 3.40659770
3.10 3.28370648
3.15 3.21143905
3.20 3.08450810
3.25 3.00432635
3.30 2.87750648
3.35 2.76226121
3.40 2.58270325
3.45 2.48994528
3.50 2.36980000
3.55 2.33617760
3.60 2.26687805
3.65 2.25039884
3.70 2.26673741
3.75 2.23859121
3.80 2.17415773
3.85 2.25203441
3.90 2.21231870
3.95 2.27370808
4.00 2.30000000
/
\endpicture
}

    %


    \medskip

The main objective  of this article is to show that under
sufficient regularity conditions on the coefficients, a large class
of semilinear see's and spde's admits
a $C^{k,\e}$ semiflow  $V: \Delta \x H \x \O \to H$ for a suitably chosen state
space $H$ and which satisfies
$V(t_0,t,x,\cdot)=u(t,x)$  for all $x \in H$ and $t \geq t_0$,
a.s., where $u$ is the solution of the see/spde with initial function
$x \in H$ at
$t=t_0$. In the autonomous case, we show further that the semiflow $V$
generates a cocycle $(U,\theta)$  on $H$, in the sense of Definition 1.2 above.  The 
cocycle and its 
Fr\'echet derivative are compact in all cases.

\bigskip
\subheading{2. Flows and cocycles of semilinear stochastic evolution equations}

In this section, we will establish the existence and regularity of semiflows
generated by mild solutions of semilinear stochastic evolution equations. We will
begin with the linear case. In fact, the linear cocycle will be used
to represent the mild solution of the semilinear stochastic evolution equation
via a variational formula which transforms the semilinear stochastic evolution 
equation to a random integral equation
(Theorem 2.5). The latter equation plays a key role in establishing the regularity 
of the stochastic flow of the semilinear see
(Theorem 2.6).

One should note at this point the fact that  Kolmogorov's continuity  theorem  fails 
for 
random fields parametrized by
infinite-dimensional spaces.  As a simple example, consider the random field $I: 
L^2([0,1], \R) \to L^2 (\O,\R)$ defined by 
the Wiener integral
$$ I(x):= \int_0^1 x(t) \, dW(t), \quad x \in L^2([0,1], \R).$$
The above random field has no continuous (or even linear!) measurable selection 
$ L^2([0,1], \R)  \times \O \to \R$ ([Mo.1], pp. 144-148; [Mo.2]).

\medskip
\noindent
{\it (a) Linear stochastic evolution equations.}

We will first prove the existence of semiflows associated with mild
solutions of
linear stochastic evolution equations of the form:
$$
\left. \aligned
du(t,x,\cdot)=&-Au(t,x,\cdot) dt+Bu(t,x,\cdot)\, dW(t),\quad t > 0\\
u(0,x,\omega)=&x.
\endaligned \right \} \tag{2.1}
$$
In the above equation  $A: D(A) \subset H \to H$ is a
closed linear operator on a separable real Hilbert space $H$. Assume
that $A$ has
a complete
orthonormal system of eigenvectors  $\{ e_n: n\geq 1\}$ with
corresponding positive
eigenvalues $\{ \mu_n, n\geq 1\}$;  i.e.,
$Ae_n=\mu_n e_n, \,\, n \geq 1.$
Suppose $-A$ generates  a strongly continuous semigroup
of bounded linear operators $T_t: H \to H, \, t\geq 0$.
Let $E$ be a separable Hilbert space and  $W(t), t\geq 0,$ be a E-valued Brownian 
motion defined on the canonical filtered
Wiener space
$(\O, \F, (\F_t)_{t\geq 0},P)$
and with  a separable covariance Hilbert space $K$. Here $K\subset E$ is a 
Hilbert-Schmidt
embedding. Indeed, $\O$ is the space
of all continuous paths $\o: \R \to E$ such that $\o (0)=0$ with the
compact open
topology, $\F$ is its Borel $\sigma$-field,  $\F_t$ is the sub-$\sigma$-field
generated by all evaluations $\O \ni \o \mapsto \o(u) \in E, u \leq t$, and $P$
is Wiener measure on $\O$. The Brownian motion is given by
$$ W(t,\o):=\o(t), \quad \o \in \O,\, t \in \R,$$
and may be represented by
$$
W(t)=\sum_{k=1}^{\i} W^k(t) f_k, \quad t \in \R,
$$
   where $\{f_k: k \geq 1\}$ is a complete orthonormal basis of $K$, and
$W^k, k \geq 1,$ are standard independent one-dimensional Wiener
processes ([D-Z.1], Chapter 4). Note that, in general, the above series converges 
absolutely in $E$ but not in $K$.


     Denote by $L_2(K,H) \subset L(K,H)$
the Hilbert space of all Hilbert-Schmidt operators $S:K \to H$, given the norm
$$\|S\|_2 := \biggl [\sum_{k=1}^\infty |S(f_k)|^2\biggr ]^{1/2},$$
where   $|\cdot|$ is the norm on $H$.
Suppose $B: H \to L_2(K,H)$ is  a bounded linear operator. The stochastic
integral in (2.1) is defined in the following sense ([D-Z.1], Chapter 4):

Let $F:[0,a] \times \O \to L_2(K,H)$ be $\B([0,a]\otimes \F,
\B(L_2(K,H)))$-measurable,
     $(\F_t)_{t \geq 0}$-adapted and such that
$\displaystyle \int_0^a E\|F(t)\|_{L_2(K,H)}^2 \,dt < \infty$. Define
$$\int_0^a F(t) \, dW(t):= \sum_{k=1}^\infty \int_0^a F(t)(f_k) \, dW^k (t)$$
where the $H$-valued stochastic integrals on the right hand side are
with respect
to the  one-dimensional Wiener processes $W^k,$
$ \, k\geq 1$.
Note that the above series converges in $L^2 (\O, H)$ because
$$
\sum_{k=1}^\infty E|\int_0^a F(t)(f_k) \, dW^k (t)|^2= \int_0^a
E\|F(t)\|^2_{L_2(K,H)} \,dt <
\infty.
$$

Throughout the rest of the article, we will denote by
$\theta: \R \times \O \to \O$ the standard $P$-preserving ergodic
Wiener shift on
$\O$:
$$\theta (t,\o)(s):=\o (t+s)-\o(t), \quad t,s \in \R.$$
Hence $(W,\theta)$ is a {\it helix}:
$$
W(t_1+t_2,\o)-W(t_1,\o)= W(t_2, \theta (t_1,\o)), \quad t_1,t_2 \in \R,\, \o
\in \O.
$$
As usual, we let $L(H)$ be the Banach space of all bounded linear operators $H
\to H$ given the
uniform operator norm $\|\cdot\|_{L(H)}$.  Denote by $L_2(H) \subset L(H)$ the 
Hilbert
space  of all
Hilbert-Schmidt operators $S: H \to H$.
%
%
It is easy to see that if $S \in L_2(H)$ and
     $T \in L(H)$, then $\|S\|_{L(H)} \leq \|S\|_2$, $T \circ S$ (and
$S\circ T$) $\in L_2(H)$ and $\|T\circ S\|_{L_2(H)} \leq \|T\|_{L(H)}
\|S\|_{L_2(H)}$.

     A {\it mild solution} of (2.1) is a family of $(\B(\R^+)\otimes \F,
\B(H))$-measurable,
$(\F_t)_{t\geq 0}$-adapted processes
$u(\cdot,x,\cdot): \R^+ \times \O \to H, \,\, x \in H,$ satisfying
the following
stochastic integral equation:
$$
u(t,x,\cdot)=T_tx+\int_0^t T_{t-s}Bu(s,x,\cdot)\, dW(s), \quad t \geq
0. \tag{2.2}
$$

\medskip
     The next lemma describes a canonical lifting of the strongly continuous
semigroup $T_t: H \to  H, \, t \geq 0,$ to a strongly continuous semigroup of
bounded linear
operators $\tilde T_t: L_2(K,H) \to L_2(K,H), t \geq 0.$

\medskip

\proclaim{Lemma 2.1} 

Define the family of maps
$\, \,\tilde T_t: L_2(K,H) \to L_2(K,H),\, t \geq 0$, by
$$\tilde T_t (C):= T_t \circ C, \quad C \in L_2(K,H), \, t \geq 0.$$
Then the following is true:
\item{(i)} $\tilde T_t, \, t \geq 0$, is a strongly continuous semigroup of
bounded linear operators on $L_2(K,H)$; and
$\|\tilde T_t\|_{L(L_2(K,H))}=\|T_t\|_{L(H)}$ for all $t \geq 0$.
\item{(ii)} If $-\tilde A:\Cal D(\tilde A) \subset L_2(K,H) \to
L_2(K,H)$ is the
infinitesimal
generator of $\tilde T_t, \, t \geq 0$, then
$$
\Cal D(\tilde A)=\{C: C\in L_2(K,H), \,  C(K) \subseteq \Cal D(A),
A\circ C \in L_2(K,H) \}
$$
and
$$\tilde A (C) = A\circ C $$
for all  $\, C \in \Cal D(\tilde A)$.
\item{(iii)} $\tilde T_t, t \geq 0,$ is a contraction semigroup if
$T_t, t \geq 0,$ is.
\endproclaim


\demo{Proof} 

Observe that each $\tilde T_t: L_2(K,H) \to L_2(K,H), t
\geq 0,$ is a bounded
linear map of $L_2(K,H)$ into itself. Indeed, it is easy to see that
$$
\|\tilde T_t(C)\|_{L_2(K,H)} \leq \|T_t\|_{L(H)} \|C\|_{L_2(K,H)}, \quad C \in
L_2(K,H),\, t \geq 0;
\tag{2.3}
$$
and hence $\|\tilde T_t \|_{L(L_2(K,H))} \leq \|T_t\|_{L(H)}$ for all
$t \geq 0$. This implies assertion (iii). The reverse inequality
$$\|T_t\|_{L(H)} \leq \|\tilde T_t \|_{L(L_2(K,H))}, \quad t \geq 0,$$
is not hard to check. Hence the last assertion in (i) holds.

We next verify the semi-group property of $\tilde T_t, t \geq 0$. Let
$t_1, t_2 \geq 0,\, C \in L_2(K,H)$. Then
$$
(\tilde T_{t_2} \circ \tilde T_{t_1})(C)=
T_{t_2} \circ (T_{t_1}\circ C)
= T_{t_1+t_2} \circ C= \tilde T_{t_1+t_2}( C).
$$
Note also that $\tilde T_0= I_{L(L_2(K,H))},$ the identity map
     $L_2(K,H) \to L_2(K,H)$. Therefore, $\tilde T_t,\, t \geq 0,$ is 
a semigroup
on $L_2(K,H)$. To prove the strong continuity of $\tilde T_t,\, t \geq 0,$
we will show that
$$ \lim_{t \to 0+} \tilde T_t(C)=C \tag{2.4}$$
for each $C \in L_2(K,H)$.  To prove the above relation, let $C \in L_2(K,H)$
and recall that $\{f_k: k \geq 1\}$ is a complete orthonormal basis of $K$.
    From the strong continuity of $T_t, t \geq 0$, it follows that
$$\lim_{t \to 0+} |T_t (C(f_k))-C(f_k)|^2_H =0 \tag{2.5}$$
for each integer $k \geq 1$. Furthermore,
$$
|T_t (C(f_k))-C(f_k)|^2_H \leq 2 [\sup_{0\leq t \leq a} \|T_t\|_{L(H)}^2
+1]|C(f_k)|_H^2, \quad k
\geq 1. \tag{2.6}
$$
Since $C$ is Hilbert-Schmidt, (2.6) implies that the series
$\sum_{k=1}^\i |T_t (C(f_k))-C(f_k)|^2_H$ converges uniformly w.r.t.
$t$. Therefore,
from (2.5), (2.6) and dominated convergence, it follows that
$$
\align
\lim_{t \to 0+} \|\tilde T_t(C)-C\|_{L_2(K,H)}^2 =& \lim_{t \to
0+}\sum_{k=1}^\i
|T_t
(C(f_k))-C(f_k)|^2_H\\
=& \sum_{k=1}^\i \lim_{t \to 0+} |T_t (C(f_k))-C(f_k)|^2_H =0. \tag{2.7}
\endalign
$$
Therefore, (2.4) holds and $\tilde T_t, t \geq 0,$ is strongly continuous.

We next prove assertion (ii) of the lemma. Let $-\tilde A: \Cal D(\tilde A)
\subset L_2(K,H) \to
L_2(K,H)$ be the infinitesimal generator of $\tilde T_t, t \geq 0$.
We begin with a proof of
the inclusion
$$
\{C: C\in L_2(K,H), \,  C(K) \subseteq \Cal D(A), A\circ C \in L_2(K,H) \}
\subseteq \Cal D(\tilde A).  \tag{2.8}
$$
Let $C\in L_2(K,H)$ be such that  $ C(K) \subseteq \Cal D(A)$ and
$A\circ C \in L_2(K,H)$. We will show that
$$\lim_{t \to 0+} \frac{\tilde T_t (C)-C}{t}= A \circ C  \tag{2.9} $$
%
in $L_2(K,H)$.  To prove (2.9), note first that
$$
\align
\sup_{0 \leq t \leq a} \frac{1}{t} |T_t (C(f_k))-C(f_k)|_H
=&\sup_{0 \leq t \leq a}\frac{1}{t}\biggl | \int_0^t T_s( A(C(f_k)))\,ds\biggl
|_H \\
&\leq \sup_{0 \leq t \leq a}\|T_s\|_{L(H)} |A(C(f_k))|_H
\tag{2.10}
\endalign
$$
because $C(f_k) \in \Cal D(A)$ for every $k \geq 1$. Since
$$
\|A\circ C\|_{L_2(K,H)}= \sum_{k=1}^{\i} |A(C(f_k))|_H^2 < \infty, \tag{2.11}
$$
it follows from (2.10), (2.11) and dominated convergence that
$$
\align
\limsup_{t \to 0+}\biggl \|\frac{\tilde T_t (C)-C}{t}-A\circ C \biggr
\|^2_{L_2(K,H)}=& \limsup_{t \to 0+} \sum_{k=1}^{\i} \biggl |\frac{T_t
(C(f_k))-C(f_k)}{t}-A(C(f_k))\biggr |^2_{H}\\
&\leq \sum_{k=1}^{\i} \limsup_{t \to 0+}
     \biggl |\frac{T_t (C(f_k))-C(f_k)}{t}-A(C(f_k))\biggr |^2_{H}\\
=&0.
\tag{2.12}
\endalign
$$
     This proves (2.9). In particular, $C \in \Cal D (\tilde A)$ and $\tilde
A(C)=A\circ C$.

It remains to prove the inclusion
$$
\Cal D(\tilde A) \subseteq \{C: C\in L_2(K,H), \,  C(K) \subseteq \Cal D(A),
     A\circ C \in L_2(K,H) \}. \tag{2.13}
$$
Suppose $C \in \Cal D(\tilde A)$. We will show that $C(K) \subseteq \Cal D(A)$,
$A \circ C \in L_2(K,H)$ and $\tilde A(C)=A\circ C$. Since
$$
\align
\lim_{t \to 0+}\biggl \|\frac{\tilde T_t (C)-C}{t}-\tilde A(C)\biggr
\|^2_{L_2(K,H)}
=& \lim_{t \to 0+}\sum_{k=1}^{\i}
        \biggl |\frac{T_t (C(f_k))-C(f_k)}{t}-\tilde A(C)(f_k)\biggr
|^2_{H}=0. \tag{2.14}
\endalign
$$
we have  that
$$\lim_{t \to 0+}
        \biggl |\frac{T_t (C(f_k))-C(f_k)}{t}-\tilde A(C)(f_k)\biggr
|^2_{H}=0\tag{2.15}
$$
for every $k \geq 1$. Therefore, $C(f_k) \in \Cal D(A)$ and $\tilde
A(C)(f_k)=A(C(f_k))$ for each $k \geq 1$. Now pick any $f \in K$ and write
$$f^n:= \sum_{k=1}^n <f,f_k> f_k, \quad n \geq 1.$$
Then $C(f^n)=\sum_{k=1}^n <f,f_k> C(f_k) \in \Cal D(A), \, n \geq 1,$ and
$C(f)= \lim_{n \to \infty} C(f^n)$ in $H$. Now since $\tilde A(C) \in L_2(K,H)
\subseteq L(K,H)$, it follows that
$$
\align
\tilde A(C)(f) =&\lim_{n \to \i} \tilde A(C)(f^n)\\
=&\lim_{n \to \i} \sum_{k=1}^n <f,f_k> \tilde A(C)(f_k)\\
=&\lim_{n \to \i} \sum_{k=1}^n <f,f_k>  A(C(f_k))\\
=& \lim_{n \to \i} A(C(f^n)).
\endalign
$$
Since $A$ is a closed operator, the above relation implies that $C(f) \in \Cal
D(A)$ and $A(C(f))=\tilde A(C)(f)$. As  $\tilde A(C) \in L_2(K,H)$,
and $f \in K$ is arbitrary, it follows that  $C(K) \subseteq \Cal D(A)$,
$\,A\circ C \in L_2(K,H)$ and $\tilde A(C)=A \circ C$. This proves (2.13) and
completes the proof of the lemma.
\hfill \qed
\enddemo

\bigskip

Our main results in this section give  regular versions $u: \R^+ \times H \times \O 
\to H$ 
of mild solutions of
(2.1) such that $u(t,\cdot,\o) \in L(H)$ for all $(t,\o) \in \R^+ \times \O$ 
(Theorems 
2.1-2.3).  These regular versions are shown to be $L(H)$-valued cocycles with 
respect to 
the
Brownian
shift $\theta$ (Theorem 2.4). In order to formulate these regularity results, we 
will 
require the following lemma:

\proclaim{Lemma 2.2} 

Let $ B:H \to L_2(K,H)$ be continuous linear,
and $v:\R^+ \times \O \to L_2(H)$
be a  $(\B(\R^+)\otimes \F, \B(H))$-measurable,  $(\F_t)_{t\geq 0}$-adapted
process such that
$\displaystyle \int_0^a E\|v(t)\|_{L_2(H)}^2 \,dt < \infty$ for each $a > 0$.
Then  the random field
$\displaystyle \int_0^t \tilde T_{t-s}(\{[B\circ v(s)](x)\})\, dW(s), \,\, x
\in H, \,t \geq  0,$
admits a jointly measurable version which will be denoted by
$\displaystyle \int_0^t
T_{t-s}Bv(s)\, dW(s)$ (by abuse of notation) and has the following properties:

\item{(i)} $\displaystyle \biggl [\int_0^t T_{t-s}Bv(s)\, dW(s)\biggr
](x)=\int_0^t \tilde
T_{t-s}(\{[B\circ v(s)](x)\})\, dW(s) $
for all $x \in H, t \geq 0,$ a.s..
\item{(ii)} For a.a. $\o \in \O$ and each $t \geq 0$, the map
$$
H \ni x \mapsto \biggl [ \biggl (\int_0^t T_{t-s}Bv(s)\, dW(s)\biggr )(\o)\biggr 
](x) \in H
$$
is Hilbert-Schmidt.
%
%
\endproclaim

\newpage

\demo{Proof} 

To prove the lemma, we will define
$\displaystyle \int_0^t T_{t-s}Bv(s)\, dW(s)$ as an It\^o  stochastic integral
with values in the Hilbert space $L_2(H)$ in the  sense of [D-Z.1],
Chapter 4). To do this, we
will introduce the following notation.

For any $V \in L_2(H)$ and $ B \in L(H, L_2(K,H)$, define the linear map
$B\star V: K \to L_2(H)$ by
$$(B \star V)(f)(x):= B(V(x))(f), \quad f \in K,\, x \in H. \tag{2.16}$$
Then $B \star V \in L_2(K,L_2(H))$ because of the following computation
$$
\allowdisplaybreaks
\align
\|B \star V\|^2_{L_2(K,L_2(H))} =& \sum_{k=1}^{\infty} \|(B \star
V)(f_k)\|^2_{L_2(H)}\\
=& \sum_{k=1}^{\infty} \sum_{n=1}^{\infty} |(B \star V)(f_k)(e_n)|^2_{H}\\
=&\sum_{k=1}^{\infty} \sum_{n=1}^{\infty} |B(V(e_n))(f_k)|^2 \\
=&\sum_{n=1}^{\infty} \sum_{k=1}^{\infty} |B(V(e_n))(f_k)|^2 \\
=&\sum_{n=1}^{\infty}  \|B(V(e_n))\|^2_{L_2(K,H)} \\
&\leq \|B\|^2_{L(H, L_2(K,H))} \|V\|^2_{L_2(H)} < \infty.
\endalign
$$
Now let $v:\R^+ \times \O \to L_2(H)$ be as in the lemma. Denote by
$\tilde{\tilde T_t} :L_2(K, L_2(H)) \to L_2(K, L_2(H))$ the induced
lifting of $T_t: H \to H,
\,\, t \geq 0,$ via Lemma 2.1; i.e.
$$
\tilde{\tilde T_t}(C)(f):=T_t\circ C(f), \quad C \in L_2(K,L_2(H)),\, f \in K.
$$
Fix $t \in [0,a]$. Then the process
$[0,t] \ni s \mapsto \tilde{\tilde T}_{t-s}(B\star v(s)) \in L_2(K,L_2(H))$
is $(\Cal F_s)_{0\leq s \leq t}$-adapted  and square-integrable, viz.
$$
\align
E \int_0^t \|\tilde {\tilde T}_{t-s}(B\star &v(s))\|^2_{L_2(K,L_2(H))}\, ds \\
&\leq \|B\|^2_{L(H, L_2(K,H))} \sup_{0 \leq u \leq t} \|T_u\|^2_{L(H)} \int_0^t
E\|v(s)\|_{L_2(H)}^2 \,ds < \infty.
\endalign
$$
In view of this, the $L_2(H)$-valued It\^o  stochastic integral
$\displaystyle\int_0^t \tilde{\tilde T}_{t-s}(B\star v(s))\, dW(s)$
is well-defined ([D-Z.1],
Chapter 4). For simplicity of notation, we will denote this
stochastic integral by
$$
\int_0^t T_{t-s}Bv(s)\, dW(s):=\int_0^t \tilde{\tilde T}_{t-s}(B\star
v(s))\, dW(s). \tag{2.17}
$$
This gives the required version of the  the random field
$\displaystyle \int_0^t \tilde T_{t-s}(\{[B\circ v(s)](x)\})\, dW(s), \,\, x
\in H, \,t \geq  0,$ because
$$
\allowdisplaybreaks
\align
\biggl [\int_0^t T_{t-s}Bv(s)\, dW(s)\biggr ](x)&:=\biggl [\int_0^t
\tilde{\tilde
T}_{t-s}(B\star v(s))\, dW(s)\biggr ](x) \\
=&\sum_{k=1}^{\i} \int_0^t \biggl [\tilde{\tilde T}_{t-s}(B\star
v(s))(f_k)\biggr ](x)\,
dW^k(s)\\
=&\sum_{k=1}^{\i} \int_0^t T_{t-s}\{B\star v(s))(f_k)(x)\}\, dW^k(s)\\
=&\sum_{k=1}^{\i} \int_0^t T_{t-s}\{B(v(s)(x))(f_k)\}\, dW^k(s)\\
=&\sum_{k=1}^{\i} \int_0^t \tilde T_{t-s}\{[B\circ v(s)](x)\}(f_k)\, dW^k(s)\\
=&\int_0^t \tilde T_{t-s}\{[B\circ v(s)](x)\}\, dW(s)
\endalign
$$
for all $x \in H$ and $t \geq 0$ a.s.. In the above computation, we have used
     the fact that for fixed $x \in H$, the It\^o stochastic integral commutes
with
the continuous linear evaluation map $L_2(H) \ni T \mapsto T(x) \in H$.
%
%
%
%
\hfill \qed
\enddemo

\bigskip

\proclaim{Theorem 2.1} 

Assume that  for some $\a \in (0,1)$,
$A^{-\a}$ is trace-class, i.e.,
$\displaystyle \sum_{n=1}^{\infty} \mu_n^{-\a} <\infty. $
Then the mild solution of the linear stochastic evolution equation (2.1) has a
$(\B(\R^+)\otimes \B(H) \otimes \F, \B(H))$-measurable
     version $u: \R^+ \times H \times \O \to H$ with the following properties:

\item{(i)} For each $x \in H$, the process $u(\cdot,x,\cdot): \R^+
\times \O \to
H$
is $(\B(\R^+)\otimes \F, \B(H))$-measurable,  $(\F_t)_{t\geq 0}$-adapted and
satisfies the stochastic integral equation (2.2).
\item{(ii)} For almost all $\o \in \O$, the map $[0,\infty) \times H  \ni (t,x)
\to u(t,x,\omega)\in H$ is jointly continuous.
Furthermore, for any fixed $ a \in \R^+$,
$$
\displaystyle E \sup_{0\leq t \leq a}\|u(t,\cdot,\cdot)\|^{2p}_{L(H)} < \infty,
$$
whenever $p \in (1, \a^{-1}]$.
%
%
\item{(iii)}For each $t>0$ and almost all $\o \in \O$, $u(t,\cdot,\omega):
H\rightarrow H$ is  a
Hilbert-Schmidt operator with the following representation:
$$
\align
u(t,\cdot,\cdot)
=T_t+\sum_{n=1}^{\infty}\int_0^t T_{t-s_1}B
&\int_0^{s_1}T_{s_1-s_2}B\,\cdots\\
&\cdots \int_0^{s_{n-1}}T_{s_{n-1}-s_n}BT_{s_n} \,dW(s_n)\cdots \, dW(s_2)\,
dW(s_1).\tag{2.18}
\endalign
$$
In the above equation, the iterated It\^o stochastic integrals are
interpreted in the sense of
Lemma 2.2,  and the convergence of the series holds in the
Hilbert space
${L_2(H)}$ of Hilbert-Schmidt operators on $H$.
\item{(iv)}
For almost all $\o \in \O$, the path $[0, \infty) \ni t \mapsto u(t,\cdot,
\o)-T_t \in L_2(H)$ is
continuous. In particular, the path $(0, \infty) \ni t \mapsto
u(t,\cdot, \o) \in
L_2(H)$
is continuous for a.a. $\o \in \O$. Furthermore,  the process $u: (0,\infty)
\times \O \to L_2(H)$
is $(\F_t)_{t\geq 0}$-adapted and $(\B((0,\infty))\otimes \F,
\B(L_2(H)))$-measurable.
\endproclaim


\demo{Proof} 

Under the hypotheses on $A$, it is well known that the see (2.1)
has a unique
$(\F_t)_{t\geq
0}$-adapted mild solution $u$ satisfying the integral equation (2.2) in $H$.
Moreover, (2.2) and a simple application of the It\^o isometry together with
Gronwall's lemma implies that
$$\sup_{0\leq t\leq a \atop |x| \leq 1}E[|u(t,x,\cdot )|^2]<\infty$$
for each $a\in (0,\infty)$.

Applying  (2.2) recursively, we obtain by induction
$$
\align
 u(&t,x,\cdot)\\
=&T_tx+\sum_{k=1}^{n}\biggl [\int_0^t
T_{t-s_1}B\int_0^{s_1}T_{s_1-s_2}B \cdots
\int_0^{s_{k-1}}
T_{s_{k-1}-s_k}BT_{s_k}\,dW(s_k)\cdots dW(s_2)\,
dW(s_1)\biggr ]x
\\
& +\int_0^t T_{t-s_1}B\int_0^{s_1}T_{s_1-s_2}B \cdots
\int_0^{s_n}T_{s_n-s_{n+1}}Bu(s_{n+1},x,\cdot)\, dW(s_{n+1})\cdots dW(s_2)\,
dW(s_1).
\endalign
$$
Set $C_t:=\displaystyle \sup_{0\leq s\leq t}||T_sB||^2_{L_2(K,H)}$ for
each  $t > 0$.
Therefore,
$$
\align
     E&[|\int_0^t T_{t-s_1}B\int_0^{s_1}T_{s_1-s_2}B \cdots
\int_0^{s_n}T_{s_n-s_{n+1}}Bu(s_{n+1},x,\cdot) dW(s_{n+1})\cdots dW(s_2)\,
dW(s_1)|^2]\\
=&\int_0^t d{s_1}E[||T_{t-s_1}B\int_0^{s_1}T_{s_1-s_2}B\cdots
\int_0^{s_n}T_{s_n-s_{n+1}}Bu(s_{n+1},x,\cdot) 
dW(s_{n+1})\cdots dW(s_2)||_{L_2(K,H)}^2]\\
\leq & C_t \int_0^t d{s_1}E[|\int_0^{s_1}T_{s_1-s_2}B\cdots
\int_0^{s_n}T_{s_n-s_{n+1}}Bu(s_{n+1},x,\cdot) dW(s_{n+1})\cdots dW(s_2)|^2]\\
\leq&  \cdot\cdot\cdot \cdot \cdot\\
\leq&  C_t^n \int_0^t d{s_1}\int_0^{s_1}d{s_2}\cdot\cdot\cdot
\int_0^{s_n}E[|u(s_{n+1},x,\cdot)|^2]  ds_{n+1}\leq   C_t^n M
\frac{t^n}{n!}\rightarrow 0.
\endalign
$$
This gives  the following series representation of $u(t,x, \cdot)$:
$$
\align
&u(t,x,\cdot)=T_tx
+\sum_{n=1}^{\infty}\biggl [\int_0^t
T_{t-s_1}B\int_0^{s_1}T_{s_1-s_2}B\cdots
\int_0^{s_{n-1}}\\
&
\hskip3.4cm T_{s_{n-1}-s_n}BT_{s_n}\, dW(s_n)\cdots
dW(s_2)\,dW(s_1)\biggr ]x
\tag{2.19}
\endalign
$$
for each $x \in H$. The above series of iterated It\^o stochastic integrals
converges in $L^2(\O,H)$
uniformly in compacta in $t$ and for $x$ in bounded
sets in $H$.

Using the fact that $A^{-1}$ is trace class, we will show further that the
series expansion (2.18) actually holds  in the Hilbert space
$L^2(\O,L_2(H))$. To see this,
first observe that $T_t$ and all the terms in the series on the right hand side
of (2.18) are
Hilbert-Schmidt for any fixed $t > 0$.
We use the comparison test to conclude that the series on the right hand side
of (2.18) converges (absolutely) in $L^2(\O,L_2(H))$.  Fix $a > 0$. Then by
successive applications of the It\^o isometry (in $L_2(H)$), one gets
$$
\allowdisplaybreaks
\align
E& \|\int_0^t T_{t-s_1}B\int_0^{s_1}T_{s_1-s_2}B\cdot\cdot\cdot
\int_0^{s_n}T_{s_n-s_{n+1}}BT_{s_{n+1}}dW(s_{n+1})\cdots
dW(s_2)\,dW(s_1)\|_{L_2(H)}^2\\
&\leq 4\int_0^a d{s_1}E[\|T_{t-s_1}B\int_0^{s_1}T_{s_1-s_2}B\cdot\cdot\cdot
\int_0^{s_n}\\
&
\hskip3cm T_{s_n-s_{n+1}}BT_{s_{n+1}}dW(s_{n+1})\cdots
dW(s_2)\|_{L_2(K,L_2(H))}^2]\\
&\leq C'_a \int_0^a d{s_1}E[\|\int_0^{s_1}T_{s_1-s_2}B\cdot\cdot\cdot
\int_0^{s_n}T_{s_n-s_{n+1}}BT_{s_{n+1}}dW (s_{n+1})\cdots dW(s_2)\|_{L_2(H)}^2]
\\
&\leq  \cdot\cdot\cdot \cdot \\
&\leq  C_a^{'n} \int_0^a d{s_1}\int_0^{s_1}d{s_2}\cdot\cdot\cdot
\int_0^{s_n}E[\|BT_{s_{n+1}}\|_{L_2(K,L_2(H))}^2]  ds_{n+1}\\
&\leq  C_a^{'n} \frac{a^n}{n!}\int_0^a \|T_s\|_{L_2(H)}^2 \,ds =C_a^{'n}
\frac{a^n}{n!} \int_0^a \sum_{k=1}^{\infty} e^{-2\mu_k s} \,ds\\
&\leq  C_a^{'n} \frac{a^n}{n!}\sum_{k=1}^{\infty}
\frac{1}{2\mu_k}, \tag{2.20}
\endalign
$$
for each integer $n \geq 1$, where $ C_a^{'}:=4C_a$.
This implies that the expansion (2.18) converges in $L^2(\Omega,
{L_2(H)})$ for each
$ t > 0$.  Hence assertion (iii) of the theorem holds.

We next prove assertion (iv). Consider the series in (2.18) and let
$\Phi^n(t)\in
L_2(H)$ be its
general term, viz.
$$
\Phi^n(t):= \int_0^t T_{t-s_1}B \int_0^{s_1}T_{s_1-s_2}B
\cdot\cdot\cdot \int_0^{s_{n-1}}T_{s_{n-1}-s_n}BT_{s_n} \,dW(s_n)\cdots
dW(s_2)\,dW(s_1),$$
for $t \geq 0, n \geq 1$. Note the relations
$$
\left. \aligned
\Phi^{n}(t)=& \int_0^t T_{t-s_1}B \Phi^{n-1}(s_1)\, dW(s_1),\quad n \geq 2, \\
\Phi^{1}(t)=& \int_0^t T_{t-s_1}B T_{s_1}\, dW(s_1),
     \endaligned \right \} \tag{2.21}
     $$
for $t \geq 0$.

First, we show by induction that for each $n \geq 1$, the process
$\Phi^n:[0,\infty) \times \O \to L_2(H)$ has a version with a.a. sample paths
continuous on $[0,\infty)$. In view of (2.21), this will follow from
Proposition (7.3) ([D-Z.1],
p. 184) provided we
show that
$$ \int_0^a E\| \Phi^{n-1}(t)\|_{L_2(H)}^{2p}\,dt < \infty \tag{2.22}$$
for all  integers $n > 1$ and $p \in (1, \a^{-1}]$. For later use, we will
actually prove the
stronger estimate
$$
E\sup_{0 \leq s \leq t}\| \Phi^{n}(s)\|_{L_2(H)}^{2p} \leq K_1 \frac{(K_2
t)^{n-1}}{(n-1)!}, \quad
t \in [0,a], \tag{2.22'}
$$
for all  integers $n\geq 1$, and $p \in (1, \a^{-1}]$, where $K_1, K_2$ are
positive constants
depending
only on $p$ and $a$. We use induction on $n$ to establish (2.22$'$). To
check (2.22$'$)
for $n=1$, choose $p \in (1, \a^{-1}]$, and consider the following easy
estimates:
$$
\align
\biggl \{\int_0^a \|T_s\|_{L_2(H)}^{2p}\, ds \biggr \}^{1/p} =& \biggl
\{\int_0^a \biggl [
\sum_{k=1}^{\infty} e^{-2\mu_k s} \biggr ]^p \,ds \biggr \}^{1/p}\\
&\leq \sum_{k=1}^{\infty} \biggl \{\int_0^a e^{-2\mu_k p s}\, ds \biggr
\}^{1/p}\\
&\leq \frac{1}{(2p)^{1/p}} \sum_{k=1}^{\infty} \mu_k^{-1/p}\\
&\leq \frac{1}{(2p)^{1/p}} \sum_{k=1}^{\infty} \mu_k^{-\a} < \infty.
\endalign
$$
Now use the second equality in (2.21) and Proposition (7.3) ([D-Z.1], p. 184) to
get the following estimate:
$$
E\sup_{0 \leq s \leq t}\| \Phi^{1}(s)\|_{L_2(H)}^{2p} \leq C_1 \int_0^t
\|T_{s_{1}}\|_{L_2(H)}^{2p}\,  ds_{1}
\leq \frac{C_1}{2p} \biggl [\sum_{k=1}^{\infty} \mu_k^{-\a}\biggr
]^p \tag{2.23}
$$
for all $t \in [0,a]$ and for $p \in (1, \a^{-1}]$. The constant $C_1$ does not
depend on $t \in
[0,a]$. Since $A^{-\a}$ is trace-class, the above inequality implies
that (2.22$'$)
holds with
$K_1:= \displaystyle\frac{C_1}{2p} \biggl [\sum_{k=1}^{\infty}
\mu_k^{-\a}\biggr ]^p$, for all
$n \geq 1, $ and any $p \in (1, \a^{-1}]$.  Now suppose that (2.22$'$)
holds for some integer $n \geq 1$ and all $p \in (1, \a^{-1}] $. Then the first
equality in (2.21) and
Proposition (7.3) ([D-Z.1], p. 184) imply  that there is a positive constant
$K_2:=K_2(p,a)$ such that
$$
\align
E \sup_{0 \leq s \leq t}\| \Phi^{n+1}(s)\|_{L_2(H)}^{2p}&\leq K_2 \int_0^{t}
E\|\Phi^{n}(s_1)\|_{L_2(H)}^{2p} \, ds_1  \\
&\leq K_2 \int_0^t   K_1 \frac{(K_2 s_1)^{n-1}}{(n-1)!}\,ds_{1} =K_1 \frac{(K_2
t)^{n}}{n!},
\tag{2.24}
\endalign
$$
for all $t \in [0,a]$ and $p \in (1, \a^{-1}]$. Therefore by
induction, (2.22$'$)
(and hence (2.22)) holds
for all integers $n \geq 1$ and  any $ p \in (1, \a^{-1}]$.

    From the first equality in (2.21), (2.22) and Proposition 7.3
([D-Z.1], p. 184),
it follows that each $\Phi^n:[0,\infty) \times \O \to L_2(H)$ has a
version with
a.a. sample paths
continuous on $[0,\infty)$.
    From the  estimate $(2.22')$, it is easy to see that
the series $\displaystyle \sum_{n=1}^{\infty} \Phi^n$ converges absolutely in
$L^{2p} (\O, C([0,a],L_2(H)))$ for each $a > 0$ and $p \in (1, \a^{-1}]$.
     This gives  a continuous modification for the sum $\displaystyle
\sum_{n=1}^{\infty} \Phi^n$
of the series in (2.18).  Hence the $L_2(H)$-valued process
$$u(t,\cdot, \cdot)-T_t= \displaystyle \sum_{n=1}^{\infty} \Phi^n(t), \quad t
\geq 0,$$
has a version with almost all sample-paths continuous on $[0, \infty)$. This
proves
the first assertion in (iv). To prove the second assertion in (iv),
it suffices to
show that the mapping $ (0,\infty) \ni t \mapsto T_t \in L_2(H)$ is
locally Lipschitz.
To see this,
let $0 < t_1 < t_2 \leq a < \infty$. Then
$$
\align
\|T_{t_2}-T_{t_1}\|_{L_2(H)}^2 &\leq \sum_{k=1}^{\i} [e^{-\mu_k t_2} -
e^{-\mu_k t_1}]^2\\
&\leq (t_2-t_1)^2 \sum_{k=1}^{\i} \mu_k^2 e^{-\mu_k t_1}\\
&\leq \frac{3}{4t_1}(t_2-t_1)^2\sum_{k=1}^{\i} \mu_k^{-1}.
\endalign
$$
Since $A^{-1}$ is trace-class, the above inequality implies that the mapping
$ (0,\infty) \ni t \mapsto T_t \in L_2(H)$ is locally Lipschitz.  The second
assertion in (iv) now
follows immediately from this and the first assertion.

     The measurability assertions in (iv) follow directly from the relation
$$
u(t,\cdot, \cdot)=T_t +\displaystyle \sum_{n=1}^{\infty} \Phi^n(t), \quad t >
0,
$$
and the fact that, as $L_2(H)$-valued It\^o stochastic integrals, the processes
$\Phi^n:
(0,\infty) \times \O \to L_2(H), \, n \geq 1,$
are  $(\F_t)_{t> 0}$-adapted and $(\B((0,\infty))\otimes \F,
\B(L_2(H)))$-measurable.

The evaluation map
$$
\gather
     L_2(H) \times H \to H\\ (S,x) \mapsto S(x)
\endgather
$$
is continuous bilinear. Therefore the first assertion in (iv) implies that
     the map $[0,T]\times H  \ni (t,x) \to u(t,x,\omega)-T_t(x) \in H$
is jointly continuous for almost all $\o \in \O$. Since $[0,T]\times
H  \ni (t,x)
\to T_t(x) \in H$
is jointly continuous (by strong continuity of the semigroup $T_t,
t\geq 0$), the first
assertion in (ii) follows.

Finally, it remains to prove the estimate in (ii). In view of
$(2.22')$, the series
in (2.18) converges
absolutely in $L^{2p} (\O, C([0,a],L_2(H))), \, p \in (1,\a^{-1}]$. Therefore,
$$
\align
\displaystyle \bigl \{E \sup_{0\leq t \leq a}\|u(t,\cdot,\cdot
)\|^{2p}_{L(H)}\bigr \}^{1/(2p)}
&\leq \sup_{0\leq t \leq a}\|T_t\|_{L(H)} + \sum_{n=1}^{\infty}\, \bigl
\{E\sup_{0\leq t \leq
a}\|\Phi^n(t)\|_{L_2(H)}^{2p} \bigr \}^{1/(2p)}\\
&\leq \sup_{0\leq t \leq a}\|T_t\|_{L(H)}+  K_1^{1/(2p)}\sum_{k=1}^{\infty}
\biggl \{ \frac{(K_2 a)^{n-1}}{(n-1)!}\biggr \}^{1/(2p)} < \infty.
\endalign
$$
This proves the estimate in (ii), and the proof of Theorem 2.1 is complete.
%
%
\hfill \qed
\enddemo

\bigskip

\proclaim{Theorem 2.2} 

Assume the following:
\item{(i)} $A^{-1}$ is a trace class operator, i.e.,
$\displaystyle \sum_{n=1}^{\infty} \mu_n^{-1} <\infty.$
\item{(ii)} $T_t \in L(H), \, t \geq 0,$ is a strongly continuous contraction
semigroup.

Then the mild solution of the linear stochastic evolution equation (2.1) has a
version $u: \R^+ \times
H \times \O \to H$ which satisfies the assertions (i),
(iii) and (iv) of Theorem 2.1. Furthermore, for almost all $\o \in \O$, the map
$[0,\infty) \times H
\ni (t,x) \to u(t,x,\omega)\in H$ is jointly continuous, and for any fixed $ a
\in \R^+$,
$$
\displaystyle E \sup_{0\leq t \leq a}\|u(t,\cdot,\cdot )\|^{2}_{L(H)} < \infty.
$$
\endproclaim

\medskip

\demo{Proof} 

The proof follows that of Theorem 2.1. We will only
highlight the differences.

We assume Hypotheses (i) and (ii). By the proof of Theorem 2.1, Hypothesis (i)
implies that the solution of (2.1) admits a version
$u: \R^+ \times H \times \O \to H$ which satisfies assertions (i) and (iii)
of Theorem 2.1.

Use the notation in the proof of Theorem 2.1. In particular, one has
$$u(t,\cdot, \cdot)-T_t= \displaystyle \sum_{n=1}^{\infty} \Phi^n(t), \quad t
\geq 0,$$
where the series converges in $L^2(\O,L_2(H))$ for each $t \geq  0$. Since
$A^{-1}$ is trace-class, then \newline
$\displaystyle \int_0^a \|T_s\|^2_{L_2(H)} \, ds < \infty$. Using
this, the fact that $T_t, \, t \geq 0,$ is a contraction semigroup,
and Theorem 6.10 ([D-Z.1], p. 160), it follows that
$\Phi^1: [0,\infty) \times \O \to L_2(H)$ has a sample-continuous version.
Furthermore, there is a
positive constant $K_3$ such that
$$
E\sup_{0 \leq s \leq t}\| \Phi^{1}(s)\|_{L_2(H)}^{2} \leq K_3 \int_0^t
\|T_{s}\|_{L_2(H)}^{2}\,  ds
< \infty,
\tag{2.25}
$$
for all $t \in [0,a]$ ([D-Z.1], Theorem 6.10, p. 160). We will show that the
series $\displaystyle \sum_{n=2}^{\i} \Phi^n$ converges in
$L^{2p}(\O, C([0,a],L_2(H))$ for all $p \geq 1$. Therefore, the series
$\displaystyle \sum_{n=1}^{\i} \Phi^n$ converges in
$L^{2}(\O, C([0,a],L_2(H)))$.

By Lemma (7.2), ([D-Z.1], p. 182), we have
$$
\align
E\| \Phi^{1}(t)\|_{L_2(H)}^{2p}& \leq K_4 \biggl [\int_0^t
\|T_{s}\|_{L_2(H)}^2\,
     ds \biggr ]^p
     \tag{2.26}
\endalign
$$
for all $t \in [0,a]$ and all $p \geq 1$. The constant $K_4$ depends on $p$ but
is independent of
$t \in [0,a]$. Since $A^{-1}$ is trace-class, the above inequality, Proposition
7.3 ([D-Z.1], p. 184) and an induction argument imply the following inequality:
$$
E\sup_{0 \leq s \leq t}\| \Phi^{n}(s)\|_{L_2(H)}^{2p} \leq K_5 \frac{(K_6
t)^{n-1}}{(n-1)!}, \quad
t \in [0,a],
$$
for all  integers $n\geq 2$, and $p \geq 1$, where $K_5, K_6$ are positive
constants depending only
on $p$ and $a$ (cf. (2.22$'$) in the proof of Theorem 2.1).
The rest of the proof of the theorem follows from the above inequality by a
similar argument to the
one in the  proof of Theorem 2.1. \hfill \qed
\enddemo
\bigskip

\proclaim{Theorem 2.3} 

Assume that $\,\,\,\displaystyle \sum_{n=1}^{\infty} 
\mu_n^{-1}\|B(e_n)\|^2_{L_2(K,H)} 
<\infty.$

Then the mild solution of the linear stochastic evolution equation (2.1) has a
$(\B(\R^+)\otimes \B(H) \otimes \F, \B(H))$-measurable
     version $u: \R^+ \times H \times \O \to H$ with the following properties:

\item{(i)} For each $x \in H$, the process $u(\cdot,x,\cdot): \R^+
\times \O \to
H$
is $(\B(\R^+)\otimes \F, \B(H))$-measurable,  $(\F_t)_{t\geq 0}$-adapted and
satisfies the stochastic integral equation (2.2).
\item{(ii)} For almost all $\o \in \O$, the map $[0,\infty) \times H  \ni (t,x)
\to u(t,x,\omega)\in H$ is jointly continuous.
Furthermore, for any fixed $ a \in \R^+$,
$$
\displaystyle E \sup_{0\leq t \leq a}\|u(t,\cdot,\cdot)\|^{2}_{L(H)} < \infty,
$$
\item{(iii)}For each $t>0$ and almost all $\o \in \O$, $u(t,\cdot,\omega):
H\rightarrow H$ is  a
bounded linear operator with the following representation:
$$
\align
u(t,\cdot,\cdot)
=T_t+\sum_{n=1}^{\infty}\int_0^t T_{t-s_1}B
&\int_0^{s_1}T_{s_1-s_2}B\,\cdots\\
&\cdots \int_0^{s_{n-1}}T_{s_{n-1}-s_n}BT_{s_n} \,dW(s_n)\cdots \, dW(s_2)\,
dW(s_1).
\endalign
$$
In the above equation, the iterated It\^o stochastic integrals are
interpreted in the sense of
Lemma 2.2,  and the convergence of the series holds in the
Hilbert space
${L_2(H)}$ of Hilbert-Schmidt operators on $H$. If in addition, $T_t: H \to H$ is 
compact for each  $t > 0$, then so is $u(t,\cdot,\omega):H\rightarrow H$  
for almost all $\o \in \O$. 
\item{(iv)} For almost all $\o \in \O$, the path 
$[0, \infty) \ni t \mapsto u(t,\cdot,\o)-T_t \in L_2(H)$ is
continuous.  Furthermore, the process $u: (0,\infty)
\times \O \to L(H)$
is $(\F_t)_{t\geq 0}$-adapted and $(\B((0,\infty))\otimes \F,
\B(L(H)))$-measurable.
\endproclaim

\medskip


\demo{Proof} 

The proof follows along the same lines as that of Theorem 2.1. Just observe that 
the hypothesis of Theorem 2.3 implies the following
integrability property
$$ \int_0^a \|BT_t\|_{L_2(H, L_2(K,H))}^2 \,dt < \infty$$
for any $a > 0$.
\qed
\enddemo

\remark{Remark}

It is easy to see that the hypothesis  of Theorem 2.3 is satisfied if one assumes 
that  
the mapping $B: H \to L_2(K,H)$ is 
Hilbert-Schmidt.   By contrast to the hypotheses of Theorems 2.1, 2.2,  the 
assumption 
in Theorem 2.3   does not entail any dimension 
restriction if the operator $A$ is a differential operator on a Euclidean domain. 
Furthermore, one does not require even discreteness of the spectrum of $A$ if we 
assume 
that $B: H \to L_2(K,H)$ is Hilbert-Schmidt. However, in this case, one gets a flow 
of 
{\it bounded linear} (but not necessarily compact) maps $u(t,\cdot, \o) \in L(H), \, 
t > 
0, \o \in \O$.
\endremark
\bigskip 

\bigskip

We will continue to assume  the hypotheses of Theorem 2.1, 2.2 or 2.3.
 
     Let $u: \R^+ \times \O \to L(H)$ be the regular version of
the mild solution of (2.1) given by Theorem 2.1, 2.2 or 2.3.
Our next result in this section identifies $u$ as a {\it fundamental solution}
(or {\it parametrix})
for (2.1).


      Consider the following stochastic integral equation:
$$
\left. \aligned
v(t)=&T_t+\int_0^t T_{t-s}Bv(s) dW(s), \quad t > 0\\
v(0)=&I,
\endaligned \right \} \tag{2.27}
$$
where  $I$ denotes the identity operator
on $H$ and the stochastic integral is interpreted as an It\^o integral in the
Hilbert space ${L_2(H)}$.


\remark{Remark} 

The initial-value problem (2.27) cannot be
interpreted strictly in the Hilbert space ${L_2(H)}$ since $v(0)=I  \not \in
{L_2(H)}$.  On the other hand, one cannot view the equation  (2.27) in the Banach
space ${ L}(H)$, because the latter Banach space is not
sufficiently ``smooth'' to allow for a satisfactory theory of stochastic
integration (cf. [B-E]).
\endremark

 \bigskip

     We say that a stochastic process  $v:[0,\infty) \times \O \to L(H)$
is a {\it solution} to equation (2.27) if

\item{(i)} $v: (0,\infty) \times \O \to L_2(H)$ is $({\Cal F}_t)_{t >0}$-adapted,
and $(\B((0,\infty))\otimes \F, \B(L_2(H)))$-measurable.
\item{(ii)}  $v \in L^2((0,a)\times \O, L_2(H))$ for all $a \in (0,\infty)$.
\item{(iii)} $v$ satisfies (2.27) almost surely.

\bigskip
\medskip

\proclaim{Theorem 2.3$'$} 

Assume the hypotheses of Theorems 2.1,
2.2 or 2.3.  Let $u$ be the regular
version of the mild
solution of (2.1) given therein. Then $u$ is  the unique solution of (2.27) in
$L^2 ((0,a)\times \O, L_2(H))$ for $a > 0$.
\endproclaim

\demo{Proof} 

Assume the hypotheses of Theorem 2.1, 2.2 or 2.3.
Let $u$ be the regular version of the mild solution of (2.1) given by these
theorems.

Note first that $u: (0,\infty) \times \O \to L_2(H)$ is $({\Cal F}_t)_{t \geq
0}$-adapted, and
$(\B((0,\infty))\otimes \F, \B(L_2(H)))$-measurable. This follows
from assertion (iv) in Theorem 2.1.

In the proofs of Theorems 2.1, 2.2, we have shown
that the series $\displaystyle \sum_{n=1}^{\infty} \Phi^n$ converges absolutely
in $L^2 (\O,
C([0,a],L_2(H)))$, and hence also in $L^2 (\O, L^2((0,a),L_2(H)))$, because of
the continuous linear imbedding
$$
L^2 (\O, C([0,a],L_2(H))) \subset L^2 (\O, L^2((0,a),L_2(H)))\equiv L^2 ((0,a)
\times \O,L_2(H)).
$$
Thus
$$
\align
\displaystyle \int_0^a E \|u(t,\cdot,\cdot )\|^2_{L_2(H)}\,dt &\leq
2 \int_0^a\|T_t\|_{L_2(H)}^2\,dt + 2\biggl [\sum_{n=1}^{\infty}
\biggl \{\int_0^a E \|\Phi^n (t)\|_{L_2(H)}^2 \, dt\biggr \}^{1/2}
\biggr ]^2 < \infty.
\endalign
$$
In particular, the It\^o  stochastic integral $ \displaystyle \int_0^t
T_{t-s}Bu(s)\, dW(s)$ is
well-defined in $L_2(H)$ for each $t \in [0,a]$ (Lemma 2.2).

We next show that $u$ solves the operator-valued stochastic integral
equation (2.27).
To see this, use the fact that
$$
\biggl [\int_0^t T_{t-s}Bu(s,\cdot,\cdot)\, dW(s)\biggr ](e_n)= \int_0^t
T_{t-s}Bu(s,e_n,\cdot)\,
dW(s), \quad n \geq 1,
$$
and the integral equation (2.2) to conclude that
$$
u(t,\o)(e_n)=T_te_n +\biggl [(\o) \int_0^t T_{t-s}Bu(s,\cdot,\cdot)\,
dW(s)\biggr ](e_n),
\quad t \geq  0, \, n \geq 1,  \tag{2.28}
$$
holds for {\it all} $\o$ in a sure event  $\O^* \in \F$   which is
independent of $n$ and $t \geq 0$. Since $\{e_n: n \geq 1\}$
is a complete orthonormal system in $H$, it follows from (2.28) that
for all $\o
\in \O^*$, one has
$$
u(t,\o)(x)=T_tx +(\o)\int_0^t \bigl [T_{t-s}Bu(s,\cdot,\cdot)\, dW(s)\bigr
](x), \quad t \geq
0, \, n \geq 1  \tag{2.29}
$$
for all $x \in H$. Thus $u$ is a solution of (2.27).

Finally we show that (2.27) has a unique $(\F_t)_{t > 0}$-adapted solution in
$L^2((0,a)\times \O,
L_2(H))$. Suppose $v_1, v_2$ are two such solutions of (2.27).
Then
$$
E\|v_1(t)-v_2(t)\|_{L_2(H)}^2 \leq \|B\|_{L_2(K,H)}\sup_{0 \leq u \leq a}
\|T_u\|_{L(H)} \int_0^t
E\|v_1(s)-v_2(s)\|_{L_2(H)}^2\, ds \tag{2.30}
$$
for all $t \in (0,a]$. The above inequality implies that
$E\|v_1(t)-v_2(t)\|_{L_2(H)}^2=0$
for all $t > 0$ and uniqueness holds. \qed

\enddemo


\medskip

From now on and throughout this section, we will impose the following 

\flushpar
{\it Condition} (B): 

\item{(i)} The operator $B: H\to L_2(K,H)$  can be extended to a bounded linear 
operator $H\to L(E,H)$, which will also be 
denoted by $B$. 
\item{(ii)} The series $\sum \limits _{k=1}^{\infty}
||B_k^2||_{L(H)}$  converges, where the bounded 
linear operators $B_k: H\to H$ are defined by $B_k(x) :=B(x)(f_k), x\in H,\, k \geq 
1$. 

\bigskip

\proclaim{Theorem 2.4} 

Suppose the hypotheses of Theorems 2.1, 2.2 or 2.3, and Condition (B)
are satisfied. Then the mild solution of (2.1) admits a version
$u: \R^+ \times \O \to L(H)$ satisfying Theorems 2.1 or 2.2 and is such that
\item{(i)}$(u,\theta)$ is a perfect $L(H)$-valued cocycle:
$$u(t+s,\o)= u(t,\theta (s,\o))\circ u(s,\o) \tag{2.31}$$
for all $s,t \geq 0$ and all $\o \in \O$;
\item{(ii)} $\displaystyle \sup_{0 \leq s \leq t \leq a}
\|u(t-s,\theta(s,\o))\|_{L(H)}
< \i,$ for all $\o \in \O$ and all $a > 0$.
\endproclaim


\demo{Proof} 

In view of Theorem 2.3$'$, $u$ satisfies the stochastic
integral equation
$$
\left. \aligned
u(t)=&T_t+\int_0^t T_{t-s}Bu(s) \,dW(s), \quad t > 0 \\
u(0)=&I
\endaligned \right \}
\tag{2.32}
$$
with $u(t) \in L_2(H)$ a.s. for all $t > 0$.

Our strategy for proving the cocycle property (2.31) is to approximate
the cylindrical
Wiener process $W$ in (2.32) by a suitably defined family of smooth processes
$W_n:\R^+ \times \O \to E, \, n \geq 1,$  prove the cocycle property for the
corresponding
approximating solutions and then pass to the limit in $L_2(H)$ as $n$ tends to
$\infty$.

%
%


Define   $W_n$ on  $\R^+ \times \O$,  $ n \geq 1,$ by
$$
W_n(t, \o):= n \int_{t - 1/n}^t W(u,\o) \,du - n \int_{- 1/n}^0 W(u,\o) \,du,
\quad t \geq 0,\, \o
\in \O. \tag{2.33}
$$
It is easy to see that each $W_n$ is a helix:
$$W_n (t,\theta (t_1,\o))= W_n(t+t_1,\o)-W_n (t_1,\o),  \tag{3.34}$$
and
$$W_n' (t,\theta (t_1,\o))= W_n'(t+t_1,\o) \tag{2.35}$$
     for all $t,t_1 \geq 0,\, \o \in \O, \, n \geq 1$. In (2.35), the prime $'$
denotes differentiation with respect to $t$.

For each $k \geq 1$, recall the definition of  $B_k:H \to H$ in Condition (B)(ii).
For each integer $n \geq 1$, define the process
$u_n: \R^+ \times \O \to L_2(H)$
to be the unique $(\B((0,\infty))\otimes \F, \B(L_2(H)))$-measurable,
$({\Cal F}_t)_{t > 0}$-adapted solution of the random integral equation:
$$
\left. \aligned
u_n(t,\o)=&T_t+\int_0^t T_{t-s}\circ \{ [B\star u_n (s,\o)](
W_n'(s,\o))\}\, ds\\
&\qquad \quad -\frac{1}{2} \int_0^t \sum_{k=1}^{\infty}T_{t-s}\circ
B_k^2\circ u_n(s,\o) \,ds, \quad t > 0\\
u_n(0,\o)=&I,
\endaligned \right \} \tag{2.27)(n}
$$
for $\o \in \O$. Recall that the operation $\star$ is defined by (2.16) in the
proof of Lemma 2.2.

Then
$$ \lim_{n \to \infty} \sup_{0 < t \leq a} \|u_n(t)-u(t)\|^2_{L_2(H)} =0,
\tag{2.36}$$
in probability, for each $a > 0$. The convergence (2.36) follows by modifying 
the proof (in $L_2 (H)$)
of the Wong-Zakai approximation theorem for
stochastic evolution equations in ([Tw], Theorem 3.4.1). 
(Cf. [I-W], Theorem 7.2, p. 497).

Next, we show that for each $n \geq 1$, $(u_n,\theta)$ is a perfect cocycle.
Fix $n \geq 1, \, t_1 \geq 0$ and $\o \in \O$. Using  $(2.27)(n)$, it follows
that
%
$$
\allowdisplaybreaks
\align
u_n(t,&\theta (t_1,\o))\circ u_n(t_1,\o)\\
=&T_t\circ u_n(t_1,\o) \\
&
+\int_{t_1}^{t+t_1}
T_{t+t_1-s}\circ \{[B\star (u_n (s-t_1,\theta (t_1,\o))\circ u_n (t_1,\o))](
W_n'(s-t_1,\theta (t_1,\o)))\}\, ds\\
& -\frac{1}{2} \int_{t_1}^{t+t_1} \sum_{k=1}^{\infty}T_{t+t_1-s}\circ
B_k^2 \circ (u_n(s-t_1,\theta (t_1,\o))\circ u_n (t_1,\o))\,ds \\
=&T_{t+t_1} + \int_0^{t_1} T_{t+t_1-s}\circ \{[B\star u_n (s,\o)](
W_n'(s,\o))\}\, ds
\\
&
-\frac{1}{2} \int_0^{t_1} \sum_{k=1}^{\infty}T_{t+t_1-s}\circ  B_k^2
\circ u_n(s,\o)
\,ds \\
&+ \int_{t_1}^{t+t_1} T_{t+t_1-s}\circ \{[B\star (u_n (s-t_1,\theta
(t_1,\o))\circ u_n
(t_1,\o))](W_n'(s-t_1,\theta (t_1,\o)))\}\, ds\\
& -\frac{1}{2} \int_{t_1}^{t+t_1} \sum_{k=1}^{\infty}T_{t+t_1-s}\circ
B_k^2 \circ (u_n(s-t_1,\theta (t_1,\o))\circ u_n (t_1,\o)) \,ds, \quad t > 0.
\endalign
$$
Hence, using $(2.27)(n)$ and (2.35), we obtain
$$
\allowdisplaybreaks
     \align
u_n(t,&\theta (t_1,\o))\circ u_n(t_1,\o)-u_n(t_1+t,\o)\\
=&\int_{t_1}^{t+t_1} T_{t+t_1-s}\circ \{[B\star (u_n (s-t_1,\theta
(t_1,\o))\circ
u_n (t_1,\o)-u_n(s,\o))]\\
&
\hskip4cm ( W_n'(s-t_1,\theta (t_1,\o)))\}\, ds\\
& -\frac{1}{2} \int_{t_1}^{t+t_1} \sum_{k=1}^{\infty}T_{t+t_1-s}\circ
B_k^2 \circ [u_n(s-t_1,\theta (t_1,\o))\circ u_n (t_1,\o)-u_n(s,\o)] \,ds \\
=&\int_0^{t} T_{t-s}\circ\{[ B\star (u_n (s,\theta (t_1,\o))\circ
u_n (t_1,\o)-u_n(s+t_1,\o))]( W_n'(s,\o))\}\, ds \\
&-\frac{1}{2} \int_0^{t} \sum_{k=1}^{\infty}T_{t-s}\circ B_k^2 \circ
[u_n (s,\theta
(t_1,\o))\circ u_n (t_1,\o)-u_n(s+t_1,\o)] \,ds,
\endalign
$$
for all $t > 0$. The above identity and a simple application of Gronwall's
lemma yields
$$u_n(t,\theta (t_1,\o))\circ u_n(t_1,\o)-u_n(t_1+t,\o)=0 \tag{2.37}$$
for all $t, t_1 \geq 0$ and all $\o \in \O$. Hence $(u_n,\theta)$ is a perfect
cocycle in $L(H)$.
Using (2.36) and passing to the limit in $L(H)$ as $n \to \infty$ in the
above identity implies that $(u,\theta)$ is a crude
$L(H)$-valued cocycle. In order to obtain a perfect version  of this cocycle,
     it is sufficient to
prove that there is a sure event $\O^* \in \F$ (independent of $t_1 \in \R^+$)
such that $\theta (t,\cdot)(\O^*) \subseteq \O^*$ for all $t \geq 0$, and
there is a subsequence $u_{n'}$ of $u_n$ such that
$$ \lim_{n', m' \to \infty}\sup_{0 < t \leq a} \|u_{n'}(t,
\o)-u_{m'}(t,\o)\|^2_{L_2(H)} =0,
\tag{2.38}$$
for each $a > 0$ and all $\o \in \O^*$. Set
$v_n(t_1,t,\o):= u_n(t-t_1,\theta (t_1,\o)), \, t \geq t_1 \geq 0$.
Then $v_n$ solves the integral equation
$$
     \align
v_n(t_1,t,\o)=&T_{t-t_1}+\int_{t_1}^t T_{t-s}\circ \{[B\star v_n
(t_1,s,\o)]( W_n'(s,\o))\}\,
ds\\
&\qquad \quad -\frac{1}{2} \int_{t_1}^t \sum_{k=1}^{\infty}T_{t-s}\circ
B_k^2 \circ v_n(t_1,s,\o) \,ds, \quad\\
v_n(t_1,t_1,\o)=&I,
\endalign
$$
for $ t \geq t_1 \geq 0$.
The above equation implies that $v_n(t_1,t,\o)$ is continuous in $(t_1,t)$
for each $\o \in \O$. Furthermore, if we apply the approximation scheme
 (in $L_2 (H)$)  to the above integral equation, we get a subsequence
$\{v_{n'}\}_{n'=1}^{\i}$ of
     $\{v_{n}\}_{n=1}^{\i}$ such that  for a.a. $\o \in \O$
$$ \lim_{n', m' \to \infty}\sup_{0 < t_1 \leq t \leq a} \|v_{n'}(t_1,t,
\o)-v_{m'}(t_1,t,\o)\|^2_{L_2(H)} =0, \tag{2.39}$$
for each $a > 0$.  Now define $\O^*$ to be the set of all
$\o \in \O$ such that
the subsequence $\{v_{n'} (t_1,t,\o):  n' \geq 1\}$ converges in
$L(H)$ uniformly
in $(t_1,t)$ for $0 < t_1 \leq t \leq a$ and all $a > 0$. Therefore $\O^*$ is a
$\theta (t,\cdot)$-invariant sure event. Define
$$u(t,\o):= \lim_{n'\to \infty} v_{n'} (0,t,\o)$$
for all $t \geq 0$ and all $\o \in \O^*$. Hence $(u,\theta)$ is a
perfect cocycle
in $L(H)$. This proves assertion (i) of the theorem.

To prove the second assertion of the theorem, fix $s \geq 0$ and define
$\hat v_n(s,t,\o):= \hat u_n(t-s,\theta (s,\o))=u_n(t-s,\theta
(s,\o))-T_{t-s}, \, t \geq s
\geq 0$. It is easy to see that  $\hat v_n$ solves the integral equation
$$
     \align
\hat v_n(s,t,\o)=&\int_{s}^t T_{t-\lambda }\circ \{[B\star \hat v_n
(s,\lambda,\o)](
W_n'(\lambda,\o))\}\, d\lambda + \int_{s}^t T_{t-\lambda }B
T_{\lambda -s}( W_n'(\lambda,\o))\, d\lambda \\
&\qquad \quad -\frac{1}{2} \int_{s}^t \sum_{k=1}^{\infty}T_{t-\lambda}\circ
B_k^2 \circ \hat v_n(s,\lambda,\o) \,d\lambda, \quad\\
\hat v_n(s,s,\o)=&0 \in L_2(H),
\endalign
$$
for $ t \geq s \geq 0$.
The above equation implies that the map $\Delta \ni (s,t) \mapsto \hat v_n(s,t,\o) 
\in 
L_2(H)$ is continuous for each $\o \in \O$. Applying  the approximation scheme 
again,
   there is a subsequence $\{\hat v_{n'}\}_{n'=1}^{\i}$ of
     $\{\hat v_{n}\}_{n=1}^{\i}$ such that for a.a. $\o \in \O$, one has
$$ \lim_{n', m' \to \infty}\sup_{0 \leq  s \leq t \leq a} \|\hat v_{n'}(s,t,
\o)-\hat v_{m'}(s,t,\o)\|^2_{L_2(H)} =0,
$$
for each $a > 0$.  Define $\hat \O^*$ to be the set of all $\o \in
\O$ such that
the subsequence $\{\hat v_{n'} (s,t,\o):  n' \geq 1\}$ converges in
$L_2(H)$ uniformly
in $(s,t)$ for $0 \leq s \leq t \leq a$ and all $a > 0$. Therefore
$\hat \O^*$ is a
$\theta (t,\cdot)$-invariant sure event. Define
$$\hat u(t,\o):= \lim_{n'\to \infty}\hat v_{n'} (0,t,\o)$$
for all $t \geq 0$ and all $\o \in \hat \O^*$. Therefore, the map
$\Delta \ni (s,t) \mapsto \hat u(t-s,\theta(s,\o)) \in L_2(H)$ is jointly 
continuous. 
In
particular,
$\displaystyle \sup_{0 \leq s \leq t \leq a} \|\hat
u(t-s,\cdot,\theta(s,\o))\|_{L(H)} < \i,$
for all $\o \in \hat\O^*$ and all $a > 0$.
Using the fact that
$\displaystyle \sup_{0 \leq s \leq t \leq a} \|T_{t-s}\|_{L(H)} < \i$, it
follows that $u(t,\o):= \hat u(t,\o)+T_t, \, \o \in \O^* \cap \hat \O^*,$
     gives a
version of the cocycle that also satisfies assertion (ii) of the theorem. This
completes the
proof of the theorem. \hfill \qed
\enddemo

\remark{Remarks}

     \item{(i)} Results analogous to Theorem  2.4 hold if $B$ is
replaced by the an affine linear map $B(x):=B_0+ B_1(x), \,\, x \in H$, where $B_0 
\in 
L(E,H)$
and $B_1: H \to L(E,H)$ satisfies Condition (B). In this case, one gets a
cocycle $(u, \theta)$ where each map $u(t,\cdot, \o): H \to H$ is of the form
$u(t,\cdot, \o)=u_0 (t,\cdot, \o)+ u_1 (t,\o)$ with
$u_0(t,\cdot, \o) \in L_2(H)$ and $u_1 (t,\o) \in H$ for $t > 0, \o \in \O$.
This follows using minor modifications of the above arguments.

\item{(ii)} It is possible to replace $B$ in the see (2.1) by an adapted random 
field   
$B: \R^+ \times H \times \O \to L(E,H)$ satisfying appropriate integrability and 
regularity conditions, which is such that 
$B(t,\cdot, \o): H \to L(E,H)$ satisfies Condition (B) for each $t \geq 0,\, \o \in 
\O$. 
The conclusions of Theorems 
2.1-2.3, 2.3$'$ 
will still hold in this case. However,
the stochastic semiflow will only satisfy Definition 1.1 (rather than the cocycle 
property 
in Definition 1.2). On the other hand if
$B$ is stationary, then the cocycle property should hold (on a suitably enlarged 
probability space) (Theorem 2.4).  

\item{(iii)} Theorems 2.1-2.4, 2.3$'$ also hold if the operator $A$ is allowed to 
have a 
non-zero discrete spectrum
$\{ \mu_n: n \geq 1\}$ which is bounded below. This yields a splitting
$A= A_0 + A_1$ where $\sigma(A_0)$ consists of positive eigenvalues
and $\sigma (A_1)$
of finitely many negative eigenvalues.
\endremark

\bigskip
\noindent
{\it (b) Semilinear stochastic evolution equations:}

In this section, we continue to assume that the operators $A,B$, the cylindrical
Brownian motion $W$, the canonical filtered Wiener space
$(\O, \F, (\F_t)_{t\geq 0},P)$ and the Brownian shift $\theta : \R
\times \O \to \O$
are as defined in part (a) of this section and satisfy the conditions therein. The 
semigroup generated
by $-A$ is denoted
as before by $T_t, t \geq 0$. Furthermore, we let
$F:H \to H$ be a (Fr\'echet) $C^1$ non-linear map satisfying the
following locally
Lipschitz and linear growth hypotheses:
$$
\left. \aligned
|F(v)| &\leq C(1+|v|), \quad v \in H\\
|F(v_1)-F(v_2)| &\leq L_n |v_1-v_2|, \quad v_i \in H,
|v_i| \leq n, i=1,2,
\endaligned \right \} \tag{2.40}
$$
for some positive constants $C,L_n, n \geq 1 $.

Consider the semilinear stochastic evolution equation:
$$
\left. \aligned
du(t)&=-Au(t) dt+F(u(t))dt+ Bu(t)\,dW(t), \quad t > 0,\\
u(0)&= x \in H,
\endaligned \right \} \tag{2.41}
$$
where the operators $A,B$ satisfy the hypotheses of Theorem 2.4.

Our main objective in this section is to establish the existence of a
$C^k$ perfect
cocycle $(U, \theta)$ for the above stochastic evolution equation.
First we define
a {\it mild solution} of (2.41) as a family of $(\B(\R^+) \otimes
\F, \B(H))$-measurable, $(\F_t)_{t\geq 0}$-adapted processes
$u(\cdot,x,\cdot): \R^+ \times \O \to H, \,\, x \in H,$ satisfying
the following 
stochastic integral equations:
$$
u(t,x,\cdot)=T_t(x)+ \int_0^t T_{t-s}(F(u(s,x,\cdot)))\, ds +
\int_0^t T_{t-s}Bu(s,x,\cdot)\,
dW(s), \quad t \geq  0, \tag{2.42}
$$
a.s. ([D-Z.1], Chapter 7, p. 182).

To fix notation, denote by $\phi: \R^+ \times \O \to L(H)$ the
perfect cocycle generated by
the linear stochastic evolution equation
$$
\left. \aligned
d\phi(t)=&-A\phi(t) dt + B\phi(t)\,dW(t), \quad t > 0,\\
\phi(0)=&I \in L(H),
\endaligned \right \} \tag{2.43}
$$
and obtained via Theorem 2.4. That is,
$\phi (t,\o):= u(t,\cdot,\o) \in L_2(H), t > 0, \o \in \O,$
in the notation of part (a) of this section.

Our first step in the construction of a non-linear cocycle of (2.41) is
to observe
that mild solutions of (2.41) correspond to solutions of a {\it
random} integral
equation on $H$. This is shown in the following theorem:

\proclaim{Theorem 2.5} 

Suppose the hypotheses of Theorem 2.4
are satisfied. Then every
$(\B(\R^+)\otimes \B(H) \otimes \F, \B(H))$-measurable,
$(\F_t)_{t\geq 0}$-adapted
solution field $U(t,x,\o)$ of the $H$-valued random integral equation
$$
U(t,x,\o)=\phi (t,\o)(x)+ \int_0^t \phi (t-s,\theta (s,\o))(F(U(s,x,\o))) \,ds,
\quad t \geq 0,\, x \in H, \tag{2.44}
$$
is a mild solution of the semilinear stochastic evolution equation (2.41).
\endproclaim


\demo{Proof} 

Let $U$ be a solution of (2.44) with the given
measurability properties. It is
sufficient to prove that $U(\cdot,x,\cdot)$ satisfies the stochastic integral
equation (2.42). Substituting from  the identity:
$$
\phi(t,\o)(x)=T_t(x)+(\o) \int_0^t T_{t-s} B \phi (s,\cdot)(x)\,
dW(s), \quad t \geq 0, x\in
H,
$$
into (2.44), gives the following a.s. relations 
$$
\allowdisplaybreaks
\align
U(t,x,\cdot)=&T_t(x)+\int_0^t T_{t-s} B \phi (s,\cdot)(x)\, dW(s)+
\int_0^t T_{t-s}(F(U(s,x,\cdot))) \,ds \\
&\qquad +\int_0^t \int_0^{t-s} T_{t-s-s'}B \phi (s',\theta
(s,\cdot))(F(U(s,x,\cdot)))\,dW(s',\theta (s,\cdot)) \,ds\\
=&T_t(x)+\int_0^t T_{t-s} B \phi (s,\cdot)(x)\, dW(s)+
\int_0^t T_{t-s}(F(U(s,x,\cdot))) \,ds \\
&\qquad +\int_0^t \int_0^{t-s} T_{t-s-s'}B \phi (s',\theta (s,\cdot))
(F(U(s,x,\cdot)))\, dW(s'+s) \,ds\\
=&T_t(x)+\int_0^t T_{t-s} B \phi (s,\cdot)(x)\, dW(s)+
\int_0^t T_{t-s}(F(U(s,x,\cdot))) \,ds \\
&\qquad +\int_0^t \int_s^{t} T_{t-\lambda}B \phi (\lambda -s,\theta
(s,\cdot))(F(U(s,x,\cdot)))\,dW(\lambda) \,ds\\
=&T_t(x)+\int_0^t T_{t-s} B \phi (s,\cdot)(x)\, dW(s)+
\int_0^t T_{t-s}(F(U(s,x,\cdot))) \,ds \\
&\qquad +\int_0^t \int_0^{\lambda} T_{t-\lambda}B \phi (\lambda -s,\theta
(s,\cdot))(F(U(s,x,\cdot)))\,ds\, dW(\lambda)\\
=&T_t(x)+\int_0^t T_{t-s}(F(U(s,x,\cdot))) \,ds \\
&\qquad+ \int_0^t T_{t-\lambda} B\{ \phi (\lambda)(x) + \int_0^{\lambda}
\phi (\lambda -s,\theta (s,\cdot))(F(U(s,x,\cdot)))\, ds\}\,dW(\lambda)\\
=&T_t(x)+ \int_0^t T_{t-s}(F(U(s,x,\cdot)))\, ds + \int_0^t
T_{t-\lambda}BU(\lambda,x,\cdot)\,
dW(\lambda)  
\endalign
$$
for $t\geq 0$.  Hence $U$ satisfies (2.42) and is therefore a mild solution of
(2.41). \hfill \qed
\enddemo

\medskip

Our next theorem shows that the random integral equation (2.44) admits a unique
$(\B(\R^+)\otimes \B(H) \otimes \F, \B(H))$-measurable,
$(\F_t)_{t\geq 0}$-adapted solution $U: \R^+ \times H \times \O \to H$. The fact 
that 
$(U,
\theta)$ is a smooth perfect cocycle can be read off from (2.44), as in the proof of 
Theorem 2.6 below.  

For any positive integer $j$, denote by 
$L^{(j)}_2 (H,H) \subset  L^{(j)}(H,H)$ the space of all Hilbert-Schmit 
$j$-multilinear 
maps 
$A \in L^{(j)}(H,H)$ given the Hilbert-Schmidt norm
$$\|A\|_{L^{(j)}_2 (H,H)} := \sum_{n_i  \geq 1 \atop {1 \leq i \leq j}} |A(e_{n_1}, 
e_{n_2}, \cdots,e_{n_j}) |_H^2  < \infty$$
where $\{e_{n_i} : n_i \geq 1 \}$ is a complete orthonormal system in $H$ for each 
$1 
\leq 
i \leq j$.

\proclaim{Theorem 2.6} 

Assume that the operators $A,B$ satisfy the
hypotheses of Theorem 2.4.
Suppose that  $F$ satisfies the linear growth and Lipschitz conditions (2.40).
Then the mild solution of (2.41) has a
$(\B(\R^+)\otimes \B(H) \otimes \F, \B(H))$-measurable
version  $U:\R^+ \times H \times \O \to H$ with the following properties:

\item{(i)}For each $x \in H$, $U(\cdot,x,\cdot): \R^+ \times \O \to H$
     is  $(\F_t)_{t\geq 0}$-adapted and satisfies (2.42) a.s..
\item{(ii)}$(U, \theta)$ is a  perfect $C^{0,1}$ cocycle (in the sense of Definition
1.2).
\item{(iii)} For each $(t,\o) \in (0, \infty) \times \O$, the map
$H \ni x \mapsto U(t,x,\o) \in H$  takes bounded sets
into relatively compact sets.
\vskip 0.2cm
Moreover, if we assume that  $F$ is $C^{k,\e}\,$ on $H$ for a positive integer
$k$ and $\e \in (0,1]$, then the mild solution $(U, \theta)$ also
enjoys  the following
properties:

\item{(iv)}$(U, \theta)$ is a $C^{k,\e}$ perfect cocycle.
\item{(v)} For each $(t,x,\o) \in \R^+ \times H \times \O$, the Fr\'echet
derivatives $D^{(j)} U(t,x,\o) \in L_2^{(j)}(H,H), \mathbreak
1 \leq j \leq k$,
and each map
$$
[0,\i) \times H \times \O \ni (t,x,\o) \mapsto  D^{(j)} U(t,x,\o) \in
L^{(j)}(H,H), \quad 1 \leq j \leq k,
$$
is strongly measurable.
\item{(vi)} For any positive $a, \rho$,
$$
E \log^+ \biggl
\{\sup_{0 \leq t_1, t_2 \leq a\atop x \in H}\frac{| U(t_2,x,\theta
(t_1,\cdot))|}{(1 + |x|)}\biggr \}
< \i
$$
and
$$
E \log^+ \sup_{0 \leq t_1, t_2 \leq a\atop |x| \leq \rho, \,  1 \leq j \leq k}
\bigl
\{\|D^{(j)}U(t_2,x,\theta (t_1,\cdot))\|_{L^{(j)}(H,H)}\bigr \} < \i .
$$
    \endproclaim

\newpage

\demo{Proof} 

In view of Theorem 2.5, we construct a version of the
mild solution of
(2.41) by applying
the classical technique of successive approximations to the integral
equation
     (2.44). Define the sequence $U_n: \R^+ \times H \times \O \to H, n \geq 1,$
by
$$
\left. \aligned
U_{n+1}(t,x,\o)=& \phi(t,\o)(x)+ \int_0^t \phi (t-s,\theta (s,\o))
(F(U_n(s,x,\o))) \,ds, \\
U_1(t,x,\o)&:= \phi(t,\o)(x)
\endaligned \right \} \tag{2.45}
$$
for all $(t,x,\o) \in \R^+ \times H \times \O$. Fix an arbitrary bounded
open set
$S$ in $H$.  Let $C_b^0 (S,H)$ denote the space of all continuous  maps
     $f:S \to H$ such that $f(S)$ is relatively compact in $H$.
%
Give  $C_b^0 (S,H)$ the  supremum  norm
$$
\,\, \|f\|_{C_b^0}:= \displaystyle \sup_{x \in S} |f(x)|_H, \quad f \in C_b^0 (S,H).
$$
\noindent
It is not hard to see that $C_b^0 (S,H)$ is a Banach space.
For fixed $\o \in \O$ and any $a > 0$, we will view the
sequence (2.45)  as a uniformly convergent sequence of bounded measurable
paths $[0,a] \ni t \mapsto U_n (t,\cdot,\o) \in  C_b^0 (S,H)$ in the
Banach
space $C_b^0 (S,H)$. To see this, we use induction on $n$.  In view of
Theorem 2.4 (ii),
define the finite random constant $\|\phi\|_{\i}:= \displaystyle \sup_{0
\leq s \leq t \leq a}
\|\phi (t-s,\theta (s,\o))\|_{L(H)}, \o \in \O$. Let $C$ be the positive constant
appearing in
(2.40). Define
$$ M_1=\sup_{x\in S}[ |x|+ C a] \|\phi\|_{\i} e^{C\|\phi\|_{\i}a}, \quad \o \in 
\O.$$
For some integer $n \geq 1$, consider the following induction hypothesis:

\noindent
{\it Hypotheses H(n)}:

\item{(i)} For each $(t,\o)  \in (0,a] \times \O$, $U_n (t,\cdot,\o) \in  C_b^0 
(S,H)$;
\item{(ii)}$ |U_n(t,x,\o)| \leq [|x|+ Ca] \|\phi\|_{\i} e^{C\|\phi\|_{\i}t}$ for
all $(t,x,\o) \in [0,a] \times H \times \O$;
\item{(iii)} $|U_{n+1} (t,x,\o)- U_{n} (t,x,\o)| \leq C [1+
\|\phi\|_{\i}|x|]L^{n-1} \|\phi\|_{\i}^n \displaystyle \frac{t^n}{n!},
\quad (t,x,\o) \in [0,a] \times H \times \O$, where $L$ is the Lipschitz constant of 
$F$ on 
the ball
$B(0,M_1)\subset H$.
\item{(iv)} $U_n:\R^+ \times H \times \O \to H$ is
$(\B(\R^+)\otimes \B(H) \otimes \F, \B(H))$-measurable,
      and for each $x \in H$, $U_n(\cdot,x,\cdot): \R^+ \times \O \to H$
     is  $(\F_t)_{t\geq 0}$-adapted.

\noindent

\medskip

We will first check that $H(1)$ is satisfied. Since $\phi (t, \cdot,\o): H
\to H$
is continuous linear for each $(t,\o) \in [0,a]\times \O$, it is clear
that $H(1)(i)$ and $H(1)(ii)$ are satisfied. Using (2.45) and the linear
growth property of $F$,
it follows that
$$
\align
|U_2 (t,x,\o)-U_1(t,x,\o)| \leq C\|\phi\|_{\i}\int_0^t [1+|\phi
(s,\o)(x)|_H] \,ds
\leq  C [1+ \|\phi\|_{\i}|x|] \|\phi\|_{\i}t,
\endalign
$$
for all $(t,x,\o) \in [0,a]\times H \times \O$. Therefore, $H(1)(iii)$ holds. To see 
the
measurability (inductive 
hypothesis $H(1)(iv)$), use the definition
of $U_1$ in (2.45) and Theorem 2.1.

Now assume that $H(n)$ holds for some integer $n \geq 1$. In particular,
for each $(t,\o) \in (0,a]\times \O$, $U_n(t,\cdot,\o)$ maps $S$ into a relatively 
compact 
set in
$H$. Therefore, the map
$$
     H \ni x \mapsto \int_0^t \phi(t-s,\theta (s,\o))(F(U_n (s,x,\o))\, ds\in
H
$$
takes $S$ into a relatively compact set in $H$, because, for fixed 
$(t,\o) \in (0,a]\times \O$,
the integrand
$$
     H \ni x \mapsto  \phi(t-s,\theta (s,\o))(F(U_n (s,x,\o))) \in H
$$
     has the same property, and is uniformly bounded in $(s,x) \in [0,t] \times S$ 
(H(n)(ii)). 
Hence, $U_{n+1}(t,\cdot,\o)(S)$ is relatively
compact
in $H$ for each $(t,\o) \in (0,a]\times \O$. Since $U_n(t,\cdot,\o): H \to H$ is
continuous, it is
easy to see from (2.45) that $U_{n+1}(t,\cdot,\o):H \to H$ is also
continuous for each $(t,\o) \in (0,a]\times \O$. Hence,
$H(n+1)(i)$ is satisfied. Using $H(n)(ii)$, the Lipschitz property of $F$
and (2.45), a
straightforward
computation shows that $H(n+1)(iii)$ is satisfied. A similar argument,
using
$H(n)(ii)$, the linear growth property of $F$ and (2.45), shows that
$H(n+1)(ii)$
also holds. To check $H(n+1)(iv)$, note first that for fixed $s \in [0,
t]$, the map
$\O \ni \o \mapsto \phi (t-s,\theta (s,\o)) \in L(H)$ is
$\F_t$-measurable. This follows
from the approximation
argument at the end of the proof of Theorem 2.4. Hence by $H(n)(iv)$, it
follows that for
fixed $s \in [0,t]$, the map $\O \ni \o \mapsto \phi (t-s,\theta
(s,\o))(F(U_n(s,x,\o))) \in
L(H)$ is $\F_t$-measurable. Hence by (2.45), it is easy to see that
$U_{n+1}(t,x,\cdot)$
is $\F_t$-measurable for fixed $(t,x) \in \R^+ \times H$. Furthermore, the
integrand on the
right-hand-side of (2.45) is
jointly-measurable in $(s,x,\o)$, and therefore $U_{n+1}(t,\cdot,\cdot)$
is jointly
measurable for any fixed $t > 0$. By continuity of the path $\R^+ \ni t
\mapsto U_{n+1}(t,
x,\o) $ for fixed $(x,\o) \in H \times \O$, the joint measurability of $U_{n+1}:\R^+ 
\times 
H \times
\O \to H$ follows.
Hence $H(n+1)(iv)$ is satisfied.
Therefore, $H(n)$ holds by induction for all integers $n \geq 1$.

The inequality $H(n)(iii)$ implies that the series
$\displaystyle \sum_{n=1}^\i [ U_{n+1} (t,\cdot,\o)- U_n(t,\cdot,\o)]$
converges in $C_b^0
(S,H)$ uniformly in $t \in [0,a]$ for each $\o \in \O$. Therefore, the sequence
$\{U_n(t,\cdot,\o)\}_{n=1}^\i$
converges
in $C_b^0 (S,H)$ uniformly in $t \in [0,a]$ for each $\o \in \O$. Its limit
$$
\lim_{n \to \i} U_n(t,\cdot,\o)= U_1(t,\cdot,\o)+ \sum_{n=1}^\i [ U_{n+1}
(t,\cdot,\o)-
U_n(t,\cdot,\o)], \quad (t,\o) \in [0,a]\times \O,
$$
is a solution of the random integral equation (2.44). Call this limit $ 
U(t,\cdot,\o) \in 
C^0_b (S,H)$ for $(t,\o) \in \R^+\times \O$. It is immediately clear from
$H(n)(iv)$ and Theorem (2.5) that $U$ satisfies the measurability requirements and 
assertion (i) of
the theorem.

We next show that $U(t,\cdot,\o): H \to H$ is $C^1$ for fixed $(t,\o) \in
\R^+ \times \O$. For
each $(x,y,\o)\in H \times H \times \O$, denote by $z(\cdot,x,y, \o)$ the
unique solution of the random linear integral equation:
$$
\align
z(t,x,y, \o)=&  \int_0^t \phi (t-s, \theta (s,\o))DF(U(s,x,\o))z(s,x,y,\o)
\,ds
\\
&\quad +\int_0^t \phi (t-s, \theta (s,\o))DF(U(s,x,\o))\phi(s,\o)(y) \,ds, \quad t > 
0.\tag{2.46}
\endalign
$$

If we suppress $y \in H$, we can view (2.46) as a linear integral equation in
$L_2(H)$ with
a unique solution $[0,\infty) \ni t \mapsto z(t,x,\cdot,\o) \in L_2 (H)$
for fixed $(x,\o) \in H \times \O$. This holds easily (by successive approximations)
because
$DF$ is  bounded on  bounded subsets of $H$ and $\{ U(t,x,\omega ); 0 \leq  t \leq 
a, 
|x|\leq M\}$ 
is bounded for any $M>0$,  and $\|\phi\|_{\infty}$ is
finite. We claim that
$U(t, \cdot, \o)$ is Fr\'echet differentiable with Fr\'echet derivative
$D U(t,x,\o) \in L_2(H)$ given by
$$D U(t,x,\o)(y)=z(t,x,y, \o) + \phi (t,\o)(y), \quad y \in H
\tag{2.47}$$
for each $(t,x,\o) \in \R^+ \times H \times \O$. To prove our claim,
define
$$
\mu (t,x,y,h,\o):= U(t,x+hy,\o)-U(t,x,\o)-h[z(t,x,y, \o) + \phi (t,\o)(y)]
\tag{2.48}
$$
for each $(t,x,y,h,\o) \in \R^+ \times H \times H \times \R \times \O$.
Using (2.48),
(2.44) and (2.46), we obtain:
%
$$
\align
\mu(t,x,y,h,\o)=&\int_0^t \phi (t-s, \theta (s,\o))DF(U(s,x,\o))
\mu(s,x,y,h,\o) \,ds\\
&\quad + \int_0^t \phi (t-s, \theta (s,\o))\biggl \{ \int_0^1 DF[\lambda
U(s,x+hy,\o)+ (1-\lambda)U(s,x,\o)]\\
&\qquad \qquad  -DF(U(s,x,\o))\,d\lambda \biggr  \}
(U(s,x+hy,\o)-U(s,x,\o)) \, ds \tag{2.49}
\endalign
$$
for all $(t,x,y,h,\o) \in \R^+ \times H \times H \times \O$.  Set
$$M_2=\sup_{s\leq a,|h|\leq 1, |y|\leq 1}\{ |U(s,x+hy,\o )|\}, \quad \o \in \O.$$
Then $M_2$ is finite for each $\o \in \O$, because of H(n)(ii).
Let $L_1 > 0$ be the  Lipschitz constant of $DF$ on the ball $B(0,M_2)$,
and $\| DF\|$ be the bound of $DF$ on $B(0,M_2)$ . Then (2.49)  
implies the following inequality:
$$
\align
|\mu(t,x,y,h,\o)| &\leq  \|\phi\|_{\i} \|DF\|\int_0^t |\mu(s,x,y,\o)|
\,ds\\
&\quad  + L_1 \|\phi\|_{\i} \int_0^t |U(s,x+hy,\o)-U(s,x,\o)|^2 \, ds
\tag{2.50}
\endalign
$$
for all $t  \in [0,a], x,y \in H, h \in \R,  |y|, |h| \leq 1, \o \in \O$. Using 
(2.44) 
and 
Gronwall's 
lemma, it is easy to see that
$$
|U(t,x+hy,\o)-U(t,x,\o)| \leq |h| \|\phi\|_{\i} |y| e^{\|\phi\|_{\i}
\|DF\|t}\tag{2.51}
$$
for all $t  \in [0,a], x,y \in H, h \in \R,  |y|, |h| \leq 1, \o \in \O$. By (2.50), 
(2.51) 
and 
another simple
application of Gronwall's lemma, we obtain
$$
|\mu(t,x,y,h,\o)| \leq  \frac{|h|^2 |y|^2 \|\phi\|_{\i}^2
L_1}{2\|DF\|}\biggl
[e^{2\|\phi\|_{\i}\|DF\|a}-1 \biggr ] e^{\|\phi\|_{\i}\|DF\|t} \tag{2.52}
$$
for all  $t  \in [0,a], x,y \in H, h \in \R,  |y|, |h| \leq 1, \o \in \O$.  Thus,
$$
\lim_{h \to 0} \frac{1}{h} \sup_{|y|\leq 1\atop{0\leq t \leq
a}}|\mu(t,x,y,h,\o)|=0
\tag{2.53}
$$
for all $x \in H, \o \in \O$. The above relation shows that
$U(t,\cdot,\o):H \to H$
is Fr\'echet differentiable at any $x \in H$ and our claim (2.47) holds.
Now combining (2.46) and (2.47), it follows that $D U(t,x,\o)$ satisfies the
$L(H)$-valued integral equation:
$$
\align
D U(t,x,\o)=& \phi (t,\o)+ \int_0^t \phi (t-s, \theta
(s,\o))DF(U(s,x,\o))
DU(s,x,\o) \,ds \tag{2.54}
\endalign
$$
for each $(t,x,\o) \in \R^+ \times H \times \O$.
In the above integral equation, the ``coefficients''
$$
\gather
     [0,\i) \times \O \ni (t,\o) \mapsto \phi (t,\o)  \in L(H)\\
\Delta \times H \times \O \ni (s,t,x,\o) \mapsto  \phi (t-s,\theta
(s,\o))DF(U(s,x,\o)) \in
L(H)
\endgather
$$
are jointly measurable, where $\Delta=\{(s,t)\in \R^2: 0 \leq s \leq t
\}$.
%
Therefore, the solution map
$$ [0,\i) \times H \times \O \ni (t,x,\o) \mapsto D U(t,x,\o) \in L(H)$$
is jointly measurable.  Furthermore, by continuity of the map
$ H \ni x\mapsto  DF(U(s,x,\o)) \in L(H,\R)$
it follows from (2.54)
that the map $H \ni x \mapsto D U(t,x,\o) \in L(H)$ is continuous for
fixed $t > 0$
and $\o \in  \O$.  Thus $U(t,\cdot,\o): H \to H$ is $C^1$. (In fact,  the
map $H \ni x
\mapsto D U(t,x,\o) \in L_2(H), t > 0,$ is continuous  because  of  the
continuity of the map $H \ni x \mapsto z(t,x,\cdot,\o) \in L_2(H)$ in the
$L_2(H)$-valued integral equation underlying (2.46).)

Suppose further that $F$ is $C^{k,\e}, k \geq 1, \e \in (0,1]$. For $k=1$, assertion 
(vi) 
of the theorem 
follows from (2.44), the linear growth property of $F$, (2.54), Gronwall's lemma and
the fact that $E \|\phi\|_{\i} < \i$. By suppressing $y$ in (2.46) and taking 
higher-order 
Fr\'echet derivatives 
with respect to $x$ of the underlying $L_2 (H)$-valued integral equation, assertions 
(v) 
and (vi) can be 
established  by induction on $k > 1$.

It remains to prove that $(U,\theta)$ is a perfect cocycle on $H$. We
use uniqueness
of solutions of (2.44). Fix $t_1 \geq 0,\, \o \in \O$ and $x \in H$. It
is sufficient to prove that
$$ U(t+t_1,x,\o)= U(t,U(t_1,x,\o),\theta (t_1,\o))  \tag{2.55}$$
for all $t \geq 0$. Define the two mappings $y,z: [0, \i) \to H$ by
$$
y(t):=U(t,U(t_1,x,\o),\theta (t_1,\o)), \quad
z(t):=U(t+t_1,x,\o) \tag{2.56}
$$
for all $t \geq 0$. Since $U$ satisfies (2.44), it follows that
$$
\align
y(t)=&\phi (t,\theta (t_1,\o))(U(t_1,x,\o))+ \int_0^t \phi
(t-s,\theta (s,\theta
(t_1,\o)))(F(U(s,U(t_1,x,\o), \theta (t_1,\o))) \,ds  \\
=&\phi (t,\theta (t_1,\o))(\phi (t_1,\o)(x))+
\int_0^{t_1} \phi (t,\theta (t_1,\o))\{\phi (t_1-s, \theta
(s,\o))(F(U(s,x,\o)))\} \,ds  \\
&\qquad + \int_{t_1}^{t+t_1} \phi (t+t_1-s, \theta (s,\o))(F(y(s-t_1))) \,ds\\
=& \phi (t+t_1,\o)(x) + \int_0^{t_1} \phi (t+t_1-s,\theta
(s,\o))(F(U(s,x,\o))) \,ds \\
&\qquad +  \int_{t_1}^{t+t_1} \phi (t+t_1-s,\theta (s,\o))(F(y(s-t_1))) \,ds
\endalign
$$
for all $ t \geq 0$. Making the substitution $t':=t+t_1$, the above
relation yields
$$
\align
y(t'-t_1)=& \phi (t',\o)(x) + \int_0^{t_1} \phi (t'-s,\theta
(s,\o))(F(U(s,x,\o))) \,ds \\
&\qquad +  \int_{t_1}^{t'} \phi (t'-s,\theta (s,\o))(F(y(s-t_1))) \,ds, \quad
t' > t_1.   \tag{2.57}
\endalign
$$
Using (2.44) and the definition of $z$, it follows that
$$
\align
z(t)=& \phi (t+t_1,\o)(x) + \int_0^{t_1} \phi (t+t_1-s,\theta
(s,\o))(F(U(s,x,\o))) \,ds \\
&\qquad +  \int_{t_1}^{t+t_1} \phi (t+t_1-s,\theta
(s,\o))(F(U(s,x,\o))) \,ds \quad t \geq 0.
\endalign
$$
Therefore,
$$
\align
z(t'-t_1)=& \phi (t',\o)(x) + \int_0^{t_1} \phi (t'-s,\theta
(s,\o))(F(U(s,x,\o))) \,ds \\
&\qquad +  \int_{t_1}^{t'} \phi (t'-s,\theta (s,\o))(F(z(s-t_1))) \,ds, \quad
t' \geq t_1.  \tag{2.58}
\endalign
$$
It is easy to see that (2.57) and (2.58) imply
$$
\align
|y(t'-t_1)-z(t'-t_1)| & \leq \int_{t_1}^{t'} \|\phi (t'-s,\theta (s,\o))\|\cdot
|F(y(s-t_1))-F(z(s-t_1))|\,ds \\
&\leq L\|\phi\|_{\i}\int_{t_1}^{t'}|y(s-t_1))-z(s-t_1)|\,ds, \quad      t_1 \leq t' 
\leq 
t_1 + a,
\tag{2.59}
\endalign
$$
   where $L$ is the Lipschitz constant of $F$ on the bounded set 
$\{ y(s), z(s), 0\leq s\leq a\}$.
   From the above inequality, we get $y(t'-t_1)-z(t'-t_1)=0$ for all
$t' \geq t_1$.
Hence, $y(t)=z(t)$ for all $t \geq 0$. This  implies the perfect cocycle
property (2.55) and  completes the proof of the theorem.
\hfill \qed
\enddemo

\remark{Remarks}

 \item{(i)} From the proof of Theorem 2.6, it is easy to see that the assertions of 
the 
theorem 
still hold if one replaces the linear growth condition $F$ by the
condition that  $F$ carries bounded sets in $H$ into bounded sets, and $U(\cdot, 
\cdot, \o )$ is bounded on bounded subsets of  $[0, \infty) \times H$.
\item{(ii)} In (2.41), it is possible to replace $F$ by a time-dependent $\tilde F: 
\R^+ 
\times H \to H$ of class 
$C^{k,\e}$ in the second variable uniformly with respect to $t$ in compacta. This 
gives 
a $C^{k,\e}$
semiflow $V: \Delta \times H \times \O \to H$ in the sense of Definition 1.1.
\endremark
\medskip

\subheading{3. Semilinear stochastic partial differential equations: Lipschitz 
nonlinearity}

Let $\Cal D$ be a smooth bounded domain in $\R^d$.
Consider the Laplacian  operator:
$$
\Delta :=\sum \limits _{i=1}^d {\partial ^2\over
\partial \xi _i^2} \tag{3.1}
$$
defined on $\Cal D$.  Let $H:=H_0^{k} (\Cal D)$ be the
Sobolev space of order $k > d/2$, i.e., the completion of $C_0^{\infty}(\Cal D)$ 
under the Sobolev norm
$$
\vert \vert u\vert \vert_{H_0^{k}}^2:=\sum_{|\alpha |\leq k}\int_\Cal D 
|D^{\alpha}u(\xi)|^2d\xi,
$$
where $d\xi$ denotes $d$-dimensional Lebesgue measure on $\R^d$. 

  Consider the spde
$$
\left. \aligned
du(t)=&{1\over 2}\Delta u(t) dt+f(u(t))dt+\sum_{i=1}^{\infty} \sigma_i 
u(t)\,dW_i(t), \quad t > 0\\
u(0)=&\psi \in  H_0^{k} (\Cal D) \\
u(t)|_{\partial \Cal D} =&0, \quad t \geq 0,
\endaligned \right \} \tag{3.2}
$$
where $f: \R \to \R$ is a $C^{\i}_b$ function, the $\sigma_i : \Cal D \to \R, i \geq 
1,$ 
are functions in the Sobolev space 
$H_0^{s} (\Cal D)$ with $s > k +\frac{d}{2}$, and the $W_i, i \geq 1,$ are standard 
independent one-dimensional Brownian motions.  Assume that the coefficients 
$\sigma_i$ in (3.2) satisfy the 
following condition
$$ \sum_{i=1}^{\infty} \|\sigma_i\|^2_{H_0^s} < \infty. \tag3.3$$
 Denote by $C^{\infty}_0 (\Cal D)$ the set of all smooth 
test functions $\phi : \Cal D \to \R$ which vanish on $\partial \Cal D$. Let 
$L^{\infty} (\Cal D)$ stand for all essentially bounded measurable functions $\psi : 
\Cal D 
\to 
\R$ with the usual norm 
$$\displaystyle \|\psi \|_{\infty}:= \displaystyle\hbox{essup}_{\xi \in \Cal D} 
|\psi 
(\xi)|. $$ 
An $(\F_t)_{t \geq 0}$-adapted random field 
$u: \R^+ \times \Cal D \times \O \to \R$ is a {\it weak solution} of (3.2) if 
$u(t,\cdot, \o) \in H_0^k (\Cal D)$ for a.a. $\o \in \O$ and the 
following identity holds:
$$
\left. \aligned
d<u(t),\phi>_{L^2} =& \nu <u(t),\Delta \phi>_{L^2} \, dt + <f(u(t)),\phi >_{L^2} \, 
dt 
+ 
\sum 
\limits _{i=1}^{\infty} <\sigma _i
u(t),\phi>_{L^2}\, 
dW_i(t),\\
    u(0) =& \psi  \in  H_0^k (\Cal D ), \\
u(t)|_{\partial \Cal D}=&0, \quad t > 0,
\endaligned \right \} 
$$
for all $\phi \in C^{\infty}_0 (\Cal D)$ a.s..  
In the above equality, $<\cdot,\cdot >_{L^2}$ denotes the inner product on  the 
Hilbert  space $L^2(\Cal D)$ of all
square-integrable functions 
$\psi: \Cal D \to \R$, viz. 
$$ <\psi_1,\psi_2>_{L^2} := \int_{\Cal D} \psi_1(\xi)\psi_2 (\xi) \, d\xi,     \quad 
\psi_1,\psi_2 \in L^2(\Cal D) ,$$
where $d\xi$ stands for $d$-dimensional Lebesgue measure.

We will show that (3.2) admits a unique weak solution $u(t) \in H$ a.s., for each 
$\psi \in 
H$. Furthermore,
the ensemble of all weak solutions of (3.2) generates a 
 a $C^{\i}$ perfect cocycle ({\it also denoted by the same symbol}) $u: \R^+ 
\times H \times \O \to H$ satisfying the assertions of Theorem 3.5 below. In 
particular, 
the 
stochastic semiflow $u(t,\cdot, \omega): H\to H$ takes bounded sets into relatively 
 compact sets in $H$.  

In this section and for the rest of the article, we should emphasize that although 
the 
weak 
{\it solution} 
 $u :\R^+ \times \Cal D \times \O \to \R$ of (3.2) and the associated
stochastic {\it semiflow} $u: \R^+ \times H \times \O \to H$ are denoted by the same 
symbol $u$, the distinction between the two notions should be clear from the  
context.

Set $A:= -{1\over 2}\Delta$ with Dirichlet boundary conditions on $\partial \Cal D$. 
 We will view 
the spde (3.2)  as a semilinear
stochastic evolution equation in $H$ of the form (2.41) (Section 2). First, define 
the 
Nemytskii operator
 $$ 
 F(u)(\xi ):= f(u(\xi )), \quad \xi \in \Cal D, u \in H.   \tag 3.3'
$$
 In Lemma 3.3 below, we will show that $F$ is a $C^{\infty}$ map $H \to H$.
Secondly,  consider a Hilbert space $K$ and assume $\{f_k\}_{k=1}^{\infty}$ is a 
complete orthonormal system in $K$. Then  
$$W(t)= \sum_{k=1}^{\infty} W^k (t) f_k, \quad t \in \R,$$
is a cylindrical Brownian motion with covariance 
space $K$.   As in section 2, $W(t)$ is a $E$-valued Brownian motion on 
the canonical filtered Wiener space
$(\O,\F,(\F_t)_{t \in \R}, P)$ for a 
separable Hilbert space $E$ such that $K\subset E$ is a Hilbert-Schmidt embedding. 
Denote by $\theta : \R \times \O \to \O$ be the standard $P$-preserving (ergodic) 
Brownian 
shift.
It is easy to see that $(W, \theta)$ is a perfect helix on $E$:
$$W(t_1+ t_2, \o)= W(t_2,\theta (t_1,\o))-W(t_1,\o), \quad t_1,t_2 \in \R, \o \in 
\O.$$
Define 
the  linear operator $B: H \to L_2 (K,H)$ by  setting
$$B(u)( f_i) :=\sigma_i u, \quad u \in H=H_0^k(\Cal D), \, i \geq 1.$$
In view of the continuous linear (Sobolev) imbedding
$$
H_0^s(\Cal D)\hookrightarrow C^k(\Cal D),
$$
it is easy to see that $B \in L(H, L_2(K,H))$ and satisfies Condition (B) of section 
2(a).   
Thirdly, observe that weak 
solutions of  the spde (3.2) correspond to mild solutions of the semilinear see: 
$$
\left. \aligned
du(t)=&-Au(t) dt+F(u(t))dt+Bu(t)\,dW(t), \quad t > 0\\
u(0)=&\psi \in  H:= H_0^{k} (\Cal D)
\endaligned \right \} \tag{3.2'}
$$
([D-Z.1], p. 156).

Finally, we will  establish a perfect $C^{\i}$-cocycle on the Sobolev space 
$H=H_0^{k} (\Cal D)$ for mild solutions of the 
 semilinear see $(3.2')$, and hence for weak solutions of the spde (3.2).

We begin with some preparations. Following 
standard notation, let $\alpha$ be a d-tuple of non-negative integers, viz.
$\alpha :=(\alpha_1,\alpha _2,\cdots,\alpha_d)$ and denote $|\alpha|:= 
\alpha_1+\alpha _2+ 
\cdots +\alpha_d$. 
For any $\phi \in C^{|\alpha|}(\Cal D)$, denote
$$(D^{(\alpha)}\phi)(\xi ) \equiv \phi ^{(\alpha)}(\xi ):=\partial _1^{\alpha_1}
\partial _2^{\alpha_2} \cdots \partial _d^{\alpha_d}\phi (\xi ), \quad \xi \in \Cal 
D,
$$
and  for any integer $l>0$, define
$$
\|D^l\phi \|_{L^2}:=\sum _{|\alpha|=l}\|D^{(\alpha) }\phi \|_{L^2}.
$$
\vskip0.2cm
\flushpar
\proclaim{Lemma 3.1} 

Let $\beta_1,\cdots ,\beta _{\mu}$ be d-tuples and
$|\alpha|=|\beta _1|+|\beta _2|+\cdots +|\beta _\mu|$,
then there exists a constant $c>0$ such that
$$
\aligned
||f_1^{(\beta_1)} &f_2^{(\beta_2)}\cdots f_{\mu}^{(\beta_{\mu})}||_{L^2}\\
\leq
&c^{\mu}
||f_1||_{L^{\infty}}^{1-{|\beta_1|\over |\alpha|}}
||f_2||_{L^{\infty}}^{1-{|\beta_2|\over |\alpha|}}
\cdots ||f_\mu||_{L^{\infty}}^{1-{|\beta_\mu|\over |\alpha|}}
||D^{|\alpha|}f_1||_{L^2}^{|\beta_1|\over |\alpha|}
||D^{|\alpha|}f_2||_{L^2}^{|\beta_2|\over |\alpha|}
||D^{|\alpha|}f_\mu||_{L^2}^{|\beta_\mu|\over |\alpha|}.
\endaligned
$$
\endproclaim


\vskip0.2cm
A proof of this lemma is given in [Ta], using Gagliardo-Nirenberg-Moser estimates.

\vskip0.2cm
\flushpar
\proclaim{Lemma 3.2} 

Let $F$ be smooth and assume $F(0)=0$.
Then for $u\in H_0^k(\Cal D)\cap L^{\infty}$,
$$
||F(u)||_{H_0^k(\Cal D)}\leq c
C_k(||u||_{L^{\infty}})(1+||u||_{L^{\infty}})^{k-1}||u||_{H_0^k(\Cal D)},
$$
where
$$
C_k(\lambda)=\sup\limits_{|u|\leq \lambda, \ \  \\ 1\leq \mu \leq k}
|F^{(\mu)}(u)|,
$$
and $c$ is a constant.
\endproclaim

\flushpar
\demo{Proof}

 We  need only prove the assertion of the lemma  for $u\in C_0^{\infty}(\Cal D)$. 
The 
chain rule gives for
  any $d$-tuple $\alpha$ with
$1\leq |\alpha| \leq k$,
$$
D^{\alpha}F(u)=\sum\limits_{\beta_1+\beta_2+\cdots+\beta_{\mu}=\alpha,
\mu \geq 1}c_{\beta}u^{(\beta_1)}u^{(\beta_2)}
\cdots
u^{(\beta_{\mu})}F^{(\mu)}(u).
$$
Hence
$$
||D^{\alpha}F(u)||_{L^2}\leq
C_k(||u||_{L^{\infty}})\sum\limits_{\beta_1+\beta_2+\cdots+\beta_{\mu}
=\alpha, \mu \geq 1}c_{\beta}||
u^{(\beta_1)}u^{(\beta_2)}
\cdots
u^{(\beta_{\mu}}||_{L^2}.
$$
Applying Lemma 3.1 to $f_i=u, i=1,2,\cdots, \mu$, we have
$$
||u^{(\beta_1)}u^{(\beta_2)}\cdots u^{(\beta_{\mu})}||_{L^2}\leq
c^{\mu}
||u||_{L^{\infty}}^{\mu -1}||D^{|\alpha|}u||_{L^2}.
$$
Therefore,
$$
\aligned
\sum\limits_{\beta_1+\beta_2+\cdots+\beta_{\mu}=\alpha, \mu \geq 1}c_{\beta}||
u^{(\beta_1)}u^{(\beta_2)} \cdots u^{(\beta_{\mu})}||_{L^2}
\leq &\sum\limits_{\beta_1+\beta_2+\cdots+\beta_{\mu}=\alpha,
\mu \geq 1}c_{\beta}
c^{\mu}
||u||_{L^{\infty}}^{\mu -1}||D^{|\alpha|}u||_{L^2}\\
\leq &
c ||D^{|\alpha|}u||_{L^2} \sum\limits_{1\leq \mu \leq |\alpha|
}C^{|\alpha|} _\mu
||u||_{L^{\infty}}^{\mu -1}\\
= &
c ||D^{|\alpha|}u||_{L^2} (1+
||u||_{L^{\infty}})^{|\alpha|  -1}\\
\leq &
c||u||_{H^k} (1+
||u||_{L^{\infty}})^{k  -1},
\endaligned
$$
for a constant $c>0$.
Note also that
$$
||F(u)||_{L^2}\leq C_k ||u||_{L^2}\leq C_k ||u||_{H^k}, \quad u \in 
C_0^{\infty}(\Cal 
D).
$$
The assertion of the lemma follows easily from the above inequality. \qed
\enddemo


\proclaim{Lemma 3.3} 

Suppose $k>\frac{d}{2}$, and $f: \R \to \R$ is a
$C^{\i}$ function.
Then the function $F: H_0^{k}(\Cal D)\rightarrow H_0^{k}(\Cal D) $ defined by 
$(3.3')$
   is a $C^{\i}$ map from $H_0^{k}(\Cal D)$ into $H_0^{k}(\Cal D)$.
\endproclaim

\demo{Proof} 

Recall the following Sobolev imbeddings 
$$
H_0^{r}(\Cal D)\hookrightarrow L^{\frac{2d}{d-2r}}(\Cal D),\quad
\quad r <\frac{d}{2},
$$
$$
H_0^{r}(\Cal D)\hookrightarrow L^{\i}(\Cal D),\quad
\quad r>\frac{d}{2}.
$$
Let us first prove that $F\in C^1(H, H )$, where $H:=H_0^{k}(\Cal D)$. Fix $u\in
H$.  We will show that $F$ is Fr\'echet differentiable and
$DF(u)(h)(\xi )\equiv S_u(h)(\xi ) =f^{\prime}(u(\xi ))h(\xi ), \,\,
h \in H, \xi \in \Cal D$. To prove this, note that  only
functions in some ball $B(0,\delta )\subset H$ centered at $0$  are
involved. By the Sobolev imbedding theorem, the range of functions in
$B(0,\delta )$ is contained in a compact interval in $\R$. Thus, we can assume $f 
\in 
C^{\i}_b$ in 
the sequel. We start by proving  that $S_u(h)\in H$ for $h\in
H$. Let $r\leq k$. By the chain and product rules, it follows that
$(S_u(h))^{(r)}$  can be
written as a finite sum whose general term is of the form:
$C(\xi )u^{(l_1)}(\xi )\cdot \cdot \cdot u^{(l_m)}(\xi )h^{(j_1)}(\xi)\cdot \cdot 
\cdot
h^{(j_n)}(\xi )$, where
$C(\cdot )\in L^{\i}(\Cal D)$, and $l_1+\cdot \cdot \cdot +l_m +j_1+\cdot \cdot
\cdot j_n=r$. Since $u^{(l)}\in H_0^{k-l}(\Cal D)$ and $h^{(j)} \in
H_0^{k-j}(\Cal D)$, the Sobolev imbedding theorem implies that  $u^{(l)}\in
L^{\frac{2d}{d-2k+2l}}(\Cal D)$ and $h^{(j)}\in
L^{\frac{2d}{d-2k+2j}}(\Cal D)$.
As
$$
\sum_{i=1}^m (d-2k+2l_i)+\sum_{i=1}^n (d-2k+2j_i)-d \leq
(m+n-1)(d-2k)<0,
$$
we have
$$
\frac{\sum_{i=1}^m (d-2k+2l_i)}{2d}+\frac{\sum_{i=1}^n
(d-2k+2j_i)}{2d}\leq \frac{1}{2}.
$$
By H\"older's inequality and the Sobolev imbedding theorem, this implies  that
$$
|C(\cdot )u^{(l_1)}(\cdot)\cdot \cdot \cdot
u^{(l_m)}(\cdot)h^{(j_1)}(\cdot)\cdot \cdot \cdot h^{(j_n)}(\cdot)|_{L^2(\Cal D)}
\leq c |u|_H^{m} |h|_H^n.
$$
where $c$ is a positive constant. Thus $S_u(h)$ is not only in $H$, but the map $H 
\ni 
h 
\mapsto S_u (h) \in H$ is a
continuous linear operator.  Now
$$
F(u+th)(\xi )-F(u)(\xi )-tS_u(h)(\xi ) =\int_0^t
[f^{\prime}(u(\xi )+sh(\xi ))-f^{\prime}(u(\xi ))]h(\xi )ds
$$
for each $\xi \in \Cal D, u,h \in H, t \geq 0$.
To show that $DF(u)=S_u$, we need to prove that
$$
\lim_{t\rightarrow 0}\sup_{|h|_H\leq 1}\biggl |\frac{1}{t}\int_0^t
[f^{\prime}\circ (u+sh)-f^{\prime}\circ (u)]\cdot h\,ds\biggr |_H=0.
$$
It is sufficient to establish
$$
\lim_{s\rightarrow 0}\sup_{|h|_H\leq 1}| [f^{\prime}\circ
(u+sh)-f^{\prime}\circ (u)] \cdot h
|_H=0.
$$
The above relation will hold  if we show that
$$
\lim_{s\rightarrow 0}\sup_{|h|_H\leq 1}|
[(f^{\prime}\circ (u+sh)-f^{\prime}\circ (u))\cdot h]^{(r)} |_{L^2(\Cal D)}=0
$$
for $r\leq k$.

Elementary computations show that
$[(f^{\prime}\circ (u+sh)-f^{\prime}\circ (u))\cdot h]^{(r)}(\xi )$
is a finite sum consisting of terms which are either of the form
$$
G_1(\xi ):=(f^{(l)}\circ (u+sh)-f^{(l)}\circ (u))(\xi )u^{(l_1)}(\xi
)...u^{(l_m)}(\xi )h^{(j_1)}(\xi )...h^{(j
_n)}(
\xi ),l\leq
r+1,
$$
or of the form
$$
G_2(\xi ) :=s^a C(\xi ) u^{(l_1)}(\xi )...u^{(l_m)}(\xi
)h^{(j_1)}(\xi )...h^{(j_n)}(\xi ), a\geq 1,C(\cdot )\in L^{\i}(\Cal D),
$$
where $l_1+...+l_m+j_1+...+j_n=r$.
For terms like $G_1$, using the Lipschitz continuity of $f^{(l)}$ it
follows that
$$
|G_1(\xi )|\leq cs |h|_{L^{\i}(\Cal D)}
|u^{(l_1)}(\xi )...u^{(l_m)}(\xi )h^{(j_1)}(\xi )...h^{(j_n)}(\xi )|.
$$
Using H\"older's inequality and the Sobolev imbedding theorem, and arguing as in the
proof of $S_u(h)\in H$, we obtain the following estimate
$$
|G_1|_{L^2(\Cal D)}\leq c s |u|_H^m |h|_H^{n+1}
$$
where $c$ is a positive constant.
Hence,
$$
    \lim_{s\rightarrow 0}\sup_{|h|_H\leq 1}|G_1|_{L^2(\Cal D)}=0.
$$
Similar arguments lead also to
$$
    \lim_{s\rightarrow 0}\sup_{|h|_H\leq 1}|G_2|_{L^2(\Cal D)}=0.
$$
Therefore,
$$
\lim_{s\rightarrow 0}\sup_{|h|_H\leq 1}|
[(f^{\prime}\circ (u+sh)-f^{\prime}\circ (u))\cdot h]^{(k)}
|_{L^2(\Cal D)}=0, \quad r\leq k,
$$
which   completes the proof that $F: H \to H$ is Fr\'echet differentiable.  The
fact that $F$ is $r$-times differentiable for $r \geq 2$ can be proved inductively  
using
similar but  lengthier computations. Details are left to the reader.
\hfill \qed
\enddemo

%
%
%

Using It\^o's formula, it is easy to see that the solution of the following 
$H$-valued linear stochastic differential equation 
$$
du^*=Bu^*dW(t),\ \ u^*(0)=\psi \in H :=H_0^k(\Cal D)
$$
is given by
$$
u^*(t,\psi, \o) (\xi) :=Q(t, \xi, \o)\psi (\xi), \quad \xi \in \Cal D, \psi \in H, t 
\geq 
0,
$$
where the process $Q: \R^+ \times \Cal D\times \O \to \R$ is defined by 
$$
Q(t,\xi, \o):= \exp \biggl \{\sum \limits _{i=1}^{\infty}\sigma_i(\xi) W_i(t,\o)- 
{1\over 
2}\sum \limits _{i=1}^{\infty}\displaystyle
\sigma_i^2(\xi)t \biggr \}, \, t \geq 0, \xi \in \Cal D, \o \in \O.
$$
 Using the perfect helix property of  $(W, \theta)$, the reader may easily check the 
following cocycle identity for $Q$:
$$
Q(t_1+t_2,\xi,\o) = Q(t_2,\xi,\theta (t_1,\o))Q(t_1,\xi,\o) \quad t_1,t_2 \geq 0, \o 
\in \O. 
$$
The above identity immediately implies that  $u^*: \R^+ \times H \times \O \to H$ is 
a 
perfect linear cocycle with
 respect to the Brownian shift  $\theta $.

We now prove the following proposition:

\vskip 0.2cm
\flushpar
\proclaim{Proposition 3.4} 

Assume $f\in C_b^k(R)$, $k>{d\over 2}$, and the forgoing conditions on the 
coefficients of 
the spde (3.2).
 Let $S$ be a bounded subset of $H_0^k(\Cal D)$. Then for any $T>0$ and almost all
$\omega \in \Omega$, the weak solution $u(t,\psi )$ of the spde (3.2) satisfies
$$
\sup_{\psi \in S}\sup_{0\leq t\leq a}\|u(t,\psi )\|_{H_0^{k}(\Cal D)} 
\leq C(\omega,a),
$$
for any $ a \in \R^+$, where  $C(\omega, a )$ is a random positive constant.
\endproclaim

\vskip 0.2cm
\flushpar
\demo{Proof}

 Let $u(t,\psi )$ be the weak solution of the spde (3.2) with initial function $\psi 
\in  
H_0^k(\Cal D)$.
Pick a sequence $\{\psi_n: n \geq 1 \}$ of smooth functions in $C^{\infty}_b (\Cal 
D)$ 
such 
that $\psi_n \to \psi$ as $n \to \infty$ in  $H_0^k(\Cal D)$. Let 
$u_n(t,\xi):= u(t,\psi_n)(\xi), \,\, t \geq 0, \xi \in \Cal D, n \geq 1$. 
Then each $u_n, \, n \geq 1,$ is a strong solution of the spde (3.2). Define 
$v_n(t,\xi ):=Q(t,\xi)^{-1}u_n(t,\xi ), \, t \geq 0, \xi \in \Cal D$.  
Using the relations
$$
dQ(t,\xi)= \sum\limits_{i=1}^{\infty}\sigma _i(\xi)Q(t,\xi)\, dW_i(t), \quad t > 0, 
\xi \in 
\Cal D, $$
$$
dQ(t,\xi)^{-1}= \sum\limits_{i=1}^{\infty}\sigma_i^2(\xi) Q(t,\xi)^{-1}\,dt 
-\sum\limits_{i=1}^{\infty}\sigma _i(\xi)Q(t,\xi)^{-1} \, dW_i(t),  \quad  t > 0, 
\xi 
\in 
\Cal D,
$$
and It\^o's formula, it follows that 
$$
\aligned
dv_n(t,\xi )
=&Q(t,\xi)^{-1}\frac{1}{2}\Delta u_n(t,\xi )dt+Q(t,\xi)^{-1}f(u_n(t,\xi ))dt 
+Q(t,\xi)^{-1}
u_n(t,\xi )\sum \limits_{i=1}^{\infty} {\sigma _i}(\xi )dW_i(t)\\
&+  u_n(t,\xi )Q(t,\xi)^{-1}\sum\limits _{i=1}^{\infty} \sigma_i^2(\xi )dt-u_n(t,\xi 
)Q(t,\xi)^{-1} \sum\limits _{i=1}^{\infty}
\sigma _i (\xi )dW_i(t)\\
&- u_n(t,\xi )Q(t,\xi)^{-1}\sum\limits _{i=1}^{\infty} \sigma _i^2(\xi )dt
\endaligned
$$
a.s. for all $t > 0, \xi \in \Cal D,$.

 Therefore, for each $n \geq 1$, $v_n (t,\xi, \o)$ satisfies the following parabolic 
equation with random coefficients:
$$
\aligned
{\partial v_n \over \partial t} =&
\frac{1}{2}\Delta v_n+ <\nabla lnQ(t,\xi), \nabla v_n>_{\R^d} -[\frac{1}{2} 
Q(t,\xi)\Delta 
Q(t,\xi)^{-1} \\
& \quad + <\nabla Q(t,\xi), \nabla Q(t,\xi)^{-1}>_{\R^d}]v_n 
+Q(t,\xi)^{-1}f(Q(t,\xi)v_n),  \quad t > 0,\\
v_n(0,\xi)=&\psi_n(\xi).
\endaligned\tag{3.4}^n
$$
Let $v$ denote the unique weak solution of the parabolic random pde 
$$
\aligned
{\partial v \over \partial t} =&
\frac{1}{2}\Delta v+ <\nabla lnQ(t,\xi), \nabla v>_{\R^d} -[\frac{1}{2} 
Q(t,\xi)\Delta 
Q(t,\xi)^{-1} \\
& \quad + <\nabla Q(t,\xi), \nabla Q(t,\xi)^{-1}>_{\R^d}]v 
+Q(t,\xi)^{-1}f(Q(t,\xi)v),  \quad t > 0,\\
v(0,\xi)=&\psi (\xi) 
\endaligned
\tag{3.4}
$$
for $t > 0, \xi \in \Cal D,$ with $\psi \in H_0^k(\Cal D)$.
Since the coefficients of (3.4) are smooth, it is well known  that
$$
\displaystyle \lim_{n\rightarrow \infty}\sup_{0\leq t\leq a}\|v(t, \cdot, \o)-v_n 
(t, 
\cdot,\o)\|_{H_0^{k}(\Cal D)}=0
$$  
for each $\o \in \O$ and any $a \in \R^+$. By rewriting (3.4), it is easy to see 
that 
$v$ satisfies the random pde   
$$
{\partial v \over \partial t}=\frac{1}{2}Q(t,\xi)^{-1}\Delta (Q(t,\xi) v)
+Q(t,\xi)^{-1}f(Q(t,\xi)v), \quad t > 0. \tag{3.5}
$$

Since $\psi\in H_0^k(\Cal D), k>{d\over 2}$, then by virtue of the Sobolev imbedding 
of 
$H_0^k(\Cal D)$ into 
$L^{\infty}(\Cal D)$, we can view (3.4) as a  random reaction diffusion equation in 
$L^{\infty}(\Cal D)$ 
whose non-linear term has linear growth and is globally Lipschitz. Hence, using 
standard 
heat-kernel estimates,
it follows that   
$$\displaystyle \sup_{t \in [0,a]} \|v(t,\cdot, \o)\|_{\i} < \i   \tag{3.6}$$
 for all $\o \in \O$. 

Now put $f \equiv 0$ in the random pde (3.5) and use uniqueness of solutions 
together 
with 
the identity 
$$
Q(t_1+t_2,\xi,\o) = Q(t_2,\xi,\theta (t_1,\o))Q(t_1,\xi,\o) \quad t_1,t_2 \geq 0, \o 
\in \O 
$$
in order to conclude that  weak solutions of the linear spde  
$$
du(t,\xi )=\frac{1}{2}\Delta u(t,\xi ) dt +B u(t,\xi )dW(t), \quad
u(0,\xi )=\psi(\xi ) \tag{3.7}
$$
yield a stochastic linear semiflow $\phi : \R^+ \times H_0^k(\Cal D) \times \O \to 
H_0^k(\Cal D)$ such that $(\phi, \theta)$ is a perfect
$L(H_0^k(\Cal D))$-valued cocycle. Full details of the argument are given in the 
proof 
of 
Theorem 4.1 in the next section.

It is easy to see that the weak solution $u$ of the spde (3.2)  satisfies the 
following 
random integral equation:
$$
u(t,\xi, \o )= \phi (t, \psi,\o)(\xi )+\int_0^t \phi(t-s,\theta 
(s,\omega))F(u(s,\xi,\omega))\, ds,  \tag{3.8}
$$
for $t \geq 0, \psi \in H_0^k(\Cal D), \xi \in \Cal D$.

Now, using Lemma 3.2 together with (3.6), one gets a positive random constant 
$C^1_k$ 
such that
$$
||F(u(t,\cdot, \o))||_{H_0^k} \leq C_{k}^1(\o,a) ||u(t, \cdot, \o)||_{H^k_0}
$$
for all $t \in [0,a], \o \in \O$ and any $a \in \R^+$.

Finally, the assertion of the proposition follows from (3.8) and a simple 
application 
of 
Gronwall's lemma.
\hfill \qed
\enddemo 

\medskip

\proclaim{Theorem 3.5} 

Suppose $k> \frac{d}{2}$. Assume $f:R\to R$ is 
a $C_b^{\infty}$ function. Assume all the forgoing conditions on the coefficients 
and the noise term in the 
spde (3.2).   Then for each 
$\psi \in H_0^{k} (\Cal D))$ the spde (3.2) has a  unique weak $(\F_t)_{t \geq 
0}$-adapted 
solution 
 $u(\cdot,\psi,\cdot): \R^+ \times \O \to H_0^{k} (\Cal D))$. Furthermore, the 
family of 
weak solutions $u(\cdot,\psi,\cdot), \, \psi \in H_0^{k} (\Cal D),$ 
admits a $(\B(\R^+)\otimes \B(H_0^{k} (\Cal D)) \otimes \F, \B(H_0^{k} (\Cal 
D)))$-measurable
version  $u:\R^+ \times H_0^{k} (\Cal D) \times \O \to H_0^{k}(\Cal D)$ having
the following properties:

\item{(i)}For each $\psi \in H_0^{k} (\Cal D)$, $u(\cdot,\psi,\cdot): \R^+
\times \O \to H_0^{k} (\Cal D)$ is  $(\F_t)_{t\geq 0}$-adapted.
\item{(ii)}$(u, \theta)$ is a $C^\i$ perfect cocycle on $H_0^{k} (\Cal D)$
(in the sense of
Definition 1.2).
\item{(iii)} For each $(t,\o) \in (0,\infty) \times \O$, the map
$H_0^{k} (\Cal D) \ni \psi \mapsto u(t,\psi,\o) \in H_0^{k} (\Cal D)$
takes bounded
sets into relatively compact sets.
\item{(iv)} For each $(t,\psi,\o) \in (0,\i) \times H_0^{k} (\Cal D) \times \O$, and
any integer $r \geq 1$,  the Fr\'echet
derivative $D^{(r)} u(t,\psi,\o) \in L^{(r)}_2(H_0^{k} (\Cal D),H_0^{k} (\Cal D))$, 
and the 
map
$$
[0,\i) \times H_0^{k} (\Cal D) \times \O \ni (t,\psi,\o) \mapsto  D^{(r)}
u(t,\psi,\o) \in L^{(r)}(H_0^{k} (\Cal D), H_0^{k} (\Cal D))
$$
is strongly measurable.
\item{(v)} For any positive $a, \rho$ and any positive integer $r$,
$$
E \log^+ \biggl
\{\sup_{0 \leq t_1, t_2 \leq a\atop \psi \in H_0^{k} (\Cal D)}
\frac{\| u(t_2,\psi,\theta (t_1,\cdot))\|_{H_0^{k} (\Cal D)}}{(1 +
\|\psi\|_{H_0^{k} (\Cal D)})}\biggr \}
< \i
$$
and
$$
E \log^+ \sup_{0 \leq t_1, t_2 \leq a\atop \|\psi\|_{H_0^{k} (\Cal D)} \leq \rho}
\bigl \{\|D^{(r)}u(t_2,\psi,\theta
(t_1,\cdot))\|_{L^{(r)}(H_0^{k} (\Cal D),H_0^{k} (\Cal D))}\bigr \} < \i.
$$
\endproclaim

\demo{Proof}

It is easy to see
that the linear cocycle $(\phi, \theta)$ of the spde (3.7) in the proof of 
Proposition 3.4 satisfies all the assertions
in Theorem 2.1 and Theorem 2.2. The theorem now holds  because of  Proposition 3.4 
and 
the 
remark following 
Theorem 2.6.
\hfill \qed
\enddemo



\medskip
\noindent
\subheading {4. Semilinear stochastic partial differential equations: Non-Lipschitz 
nonlinearity}
\medskip

 In this section,  we will study two types of  semilinear stochastic partial
differential equations with non-Lipschitz nonlinearities  and infinite dimensional 
noise. 
  
The two classes of spde's considered are {\it  stochastic reaction diffusion 
equations}
and {\it stochastic  Burgers equation with additive noise}. 
We prove the  existence of a compacting  $C^1$-cocycle in each case.

\medskip
\flushpar
{\it (a) Stochastic reaction diffusion equations}

   This class of spde's has dissipative nonlinear terms and infinite dimensional 
spatially 
smooth white noise. 
We prove the  existence of a compacting  
$C^{0,1}$-cocycle satisfying appropriate regularity properties (Theorem 4.1).  It 
appears 
that the cocycle is in general not 
Fr\'echet differentiable over the space of all $L^2$ functions on the domain (cf. 
[Te], p. 
298). However,  for a subclass of
dissipative non-linearities with a certain dimension requirement,  we further prove 
that 
the cocycle is $C^1$ and possesses 
Oseledec-type  integrability properties (Theorem 4.2). 

 In [F.2], Flandoli  studied the existence of {\it continuous} semi-flows 
for  a class of  spde's with finite dimensional noise and
polynomial nonlinearities of odd degree and with  negative leading coefficients.

Consider the following stochastic reaction diffusion equation in a 
smooth bounded domain $\Cal D \subset \R^d$,
$$
\left. \aligned
du=&\nu \Delta u \, dt+ f(u(t))  \, dt+\sum \limits _{i=1}^{\infty} \sigma _i u\, 
dW_i(t), \quad t > 0 \\
u(0)=&\psi  \\
u(t)|_{\partial \Cal D}=&0, \quad t > 0.
\endaligned \right \} \tag4.1
$$
where $\Delta$ is the Laplacian on $\Cal D$,  
 $\nu>0 $ is a real constant. The initial function $\psi:\Cal D \to \R$ is 
square-integrable with
respect to Lebesgue
measure on $\Cal D$, and a Dirichlet boundary condition is assumed on the boundary 
$\partial \Cal D$. The noise term $\sum\limits_{i=1}^{\infty} \sigma _i u\, dW_i(t)$ 
is 
very similar to the one in (3.2) of section 3, but we assume here that 
$\sigma_i:\Cal 
D\to 
R$, $i\ge 1$, are 
functions in the Sobolev space
$H_0^{s}(\Cal D)$ with $s>2+{d\over 2}$, and 
$$
\sum\limits_{i=1}^{\infty}||\sigma_i||_{H_0^s}^2<\infty.
$$
%
%
The nonlinearity $f: \R \to \R$ satisfies the following classical dissipativity 
conditions: 

\subheading{Conditions (D)}

The function $f$ is $C^2$, and there are positive constants $c_i, i=1,2,3,4,$ and a 
positive integer $p$ such that 
$$
\align
-c_2-c_3 s^{2p} \leq f(s)s &\leq c_2 -c_1 s^{2p} \\
f'(s) &\leq c_4
\endalign
$$
for all $s \in \R$. 

A typical example of a function $f: \R \to \R$ satisfying Conditions (D) is the 
polynomial 
$f(s):= \sum_{k=1}^{2p-1} a_k s^k,
\, s \in \R,$ where $a_{2p-1} < 0$.  (See e.g. [Te], pp. 83-85.)

Solutions of (4.1) are to be understood in a {\it weak} sense as defined below.

Consider the Hilbert space $H:=L^2(\Cal D)$ of all square-integrable functions 
$\psi: \Cal D \to \R$ furnished with the
$L^2$ inner product 
$$ <\psi_1,\psi_2>_H := \int_{\Cal D} \psi_1(\xi)\psi_2 (\xi) \, d\xi,     \quad 
\psi_1,\psi_2 \in H,$$
where $d\xi$ stands for Lebesgue measure on $\Cal D$.  Denote the  induced norm on 
$H$ by 
$$|\psi|_H :=\biggl [\int_{\Cal D} |\psi (\xi)|^2 \, d\xi \biggr ]^{1/2}, \quad \psi 
\in 
H.$$
 Recall $C^{\infty}_0 (\Cal D)$, the set of all smooth 
test functions $\phi : \Cal D \to \R$ which vanish on $\partial \Cal D$. Let 
$L^{\infty} 
(\Cal D)$ stand for all essentially bounded measurable functions $\psi : \Cal D \to 
\R$ 
with the usual norm 
$$\displaystyle \|\psi \|_{\infty}:= \hbox{essup}_{\xi \in \Cal D} |\psi (\xi)|.$$ 
An $(\F_t)_{t \geq 0}$-adapted random field 
$u: \R^+ \times \Cal D \times \O \to \R$ is a {\it weak solution} of (4.1) if 
$u(t,\cdot, \o) \in H$ for a.a. $\o \in \O, t > 0,$ and the 
following identity holds:
$$
\left. \aligned
d<u(t),\phi>_H =& \nu <u(t),\Delta \phi>_H\, dt + <f(u(t)),\phi >_H \, dt + \sum 
\limits 
_{i=1}^{\infty} <\sigma _i u(t),\phi>_H\, 
dW_i(t),\\
    u(0) =& \psi  \in  L^2 (\Cal D ), \\
u(t)|_{\partial \Cal D}=&0, \quad t > 0,
\endaligned \right \} 
$$
for all $\phi \in C^{\infty}_0 (\Cal D)$ a.s..  

Note that, unless $f$ has linear growth ($p=1$ in  Conditions (D)), the Nemytskii 
operator 
$F(u)(\xi):=f(u(\xi)), \,\, \xi \in \Cal D,$ does not 
even map $H=L^2(\Cal D)$ into itself. Thus one cannot view (4.1) as a semilinear see 
on $H$. 
Nevertheless,  we will show that  for each $\psi \in H$, (4.1) admits a unique weak 
solution $u(t) \in H$ 
a.s., for all $t > 0$.  Furthermore,
the ensemble of all weak solutions of (4.1) generates a 
 a globally Lipschitz cocycle ({\it also denoted by the same symbol}) $u: \R^+ 
\times H 
\times \O \to H$ satisfying the assertions of Theorem 4.1 below. In particular, the 
stochastic semiflow $u(t,\cdot, \omega): H\to H$ takes bounded sets into relatively 
 compact sets in $H$, and its global Lipschitz constant has moments of all orders.

As for Fr\'echet differentiability of the cocycle $u: \R^+ \times H \times \O \to H$ 
on the whole of $H$, it appears to be {\it not true} when $f$ is smooth and 
satisfies 
Conditions 
(D) (cf. [Te], p. 298). However, under a stronger dimension requirement on the 
polynomial
growth rate $p$ of $f$, we are able to establish that the cocycle $u$ is $C^1$ on 
$H$ 
(Theorem 4.2). 
Furthermore, it satisfies similar assertions to those of Theorem 2.6. 
 In particular, its Fr\'echet
derivatives $Du(t,\psi,\omega): H\to H$ are compact for all $(t,\psi,\omega)\in 
(0,\i) \times H\times \Omega$.
 
In (4.1), the special case $f(s):= s(1-s), \, s \in \R,$  corresponds to  the 
well-known 
stochastic
KPP equation. It is not covered by the analysis in this section since it only admits  
positive solutions for all time. 
Its random travelling wave and ergodic properties were considered in [E-Z], 
[D-T-Z.1] and 
[O-V-Z].  For the KPP equation
with additive
noise, the reader may refer to [E-H] for the existence of the invariant measure. 


The following lemma reduces (4.1) to a random family of reaction-diffusion 
equations. 

\proclaim{Lemma 4.1}

Recall the process $Q: \R^+ \times \Cal D\times \O \to \R$ defined by 
$Q(t,\xi):= \exp \biggl \{\sum \limits _{i=1}^{\infty}\sigma_i(\xi) W_i(t)- {1\over 
2}\sum \limits _{i=1}^{\infty}\displaystyle
\sigma_i^2(\xi)t \biggr \}, \, t \geq 
0$. Let $u$ be a weak solution of (4.1)
and set $v(t):=Q(t)^{-1} u(t), \, t \geq 0$. Define $\tilde f: \R^+ \times \Cal D 
\times \R 
\times 
\O \to \R$ 
by $\tilde f(t,\xi,s,\o):= Q(t,\xi,\o)^{-1} f(Q(t,\xi,\o)s),  \,  t \in \R^+, \xi 
\in 
\Cal 
D, s \in \R, \o \in \O$. Then $v$ 
is a weak solution of the random reaction-diffusion equation
$$
\left. \aligned
\frac{\partial v}{\partial t} =&\nu Q(t)^{-1}\Delta (Q(t)v) + \tilde f(t,v(t)), 
\quad 
t > 
0 \\
v(0)=&\psi  \in  L^2 (\Cal D) \\
v(t)|_{\partial \Cal D}=&0, \quad t > 0.
\endaligned \right \} \tag4.2
$$
Conversely, every weak solution $v$ of (4.2) corresponds to a weak solution $u$ of 
(4.1) 
given by 
$u(t):=Q(t)v(t), \, t \geq 0$.
\endproclaim


\demo{Proof}

Suppose $u$ is a weak solution of (4.1) with initial function $\psi \in L^2 (\Cal 
D)$. 
Define 
$$v(t,\xi,\o):=Q(t,\xi,\o)^{-1} 
u(t,\xi,\o), \, t \geq 0, \xi \in \Cal D,\o \in \O. \tag{4.3}$$

Assume first that the initial function $\psi : \Cal D \to \R$ is smooth. Then $u$ is 
a 
strong solution of (4.1). Hence by It\^o's formula (as in the proof of Proposition 
3.4), 
it follows that $v$ is a (strong) solution of the random reaction-diffusion equation 
(4.2). 
The case of a general $\psi \in L^2 (\Cal D)$ can be handled by approximating $\psi$ 
in the 
$L^2$-norm  by a sequence of smooth functions $\psi_n : \Cal D \to \R, \, n \geq 1$, 
as in 
the proof of Proposition 3.4.

A similar argument, using It\^o's formula and the relation 
$$dQ(t,\xi)= \sum\limits_{i=1}^{\infty}\sigma _i(\xi)Q(t,\xi)\, dW_i(t), \quad t > 
0, \tag 4.4$$
proves the second assertion  of the lemma. 
\qed
\enddemo

The next lemma shows that the non-linear term $\tilde f$ in (4.2) inherits the 
dissipativity properties of the original non-linear
term $f$ in (4.1).

\proclaim{Lemma 4.2}

Suppose $f$ satisfies Conditions (D).  Let $0 < a < \i$. Then there exist $\Cal 
F$-measurable positive random
variables $\tilde c_i  \in \displaystyle \bigcap_{k=1}^{\infty} L^k (\O, \R),  
i=1,2,3,$  
such that the following is true:
$$
\left. \aligned
-\tilde c_2(\o) -\tilde c_3(\o) s^{2p} \leq \tilde f(t,\xi,s, \o) s &\leq -\tilde 
c_1(\o) 
s^{2p} + \tilde c_2(\o), \\
\frac{\partial \tilde f (t,\xi,s,\o)}{\partial s} &\leq  c_4 
\endaligned \right \} \tag4.5
$$
for  all $t \in [0,a], s \in \R, \o \in \O$.
\endproclaim


\demo{Proof}

Fix $ a \in (0,\i), 0 \leq t \leq a,\xi\in \Cal D, s \in \R, \o \in \O$. Then 
Conditions (D) 
imply that
$$
\aligned
- c_2 Q(t,\xi,\o)^{-2} -c_3 Q(t,\xi,\o)^{(2p-2)} s^{2p} \leq &\tilde
f(t,\xi,s, \o) s\\
\leq & - c_1 
Q(t,\xi,\o)^{(2p-2)} s^{2p} +  c_2 Q(t,\xi,\o)^{-2}, \\
\frac{\partial \tilde f (t,s,\o)}{\partial s} &\leq  c_4. 
\endaligned 
$$
Define 
$$  
\gather
 \tilde c_1 (\o):=c_1 \inf_{0 \leq t \leq a, \xi\in \Cal D} Q(t,\xi,\o)^{(2p-2)},  
\,\,     \tilde c_2 
(\o):=c_2 \sup_{0 \leq t \leq a,\xi\in \Cal D} Q(t,\xi,\o)^{-2}\\
\tilde c_3 (\o):=c_3 \sup_{0 \leq t \leq a,\xi\in \Cal D} Q(t,\xi,\o)^{(2p-2)} 
\endgather
$$
for all $\o \in \O$. By sample continuity of $Q(t,\xi)$ and $Q(t,\xi)^{-1}$, it is 
clear that each 
$\tilde c_i (\o), i=1,2,3,$ is finite for a.a.
$\o \in \O$. The estimates of the lemma follow immediately from the above 
inequalities and 
the definition of $Q(t,\xi)$.  The existence of all moments of $\tilde c_i, 
i=1,2,3$, 
follows from 
Burkholder-Davis-Gundy inequality and 
the fact that $Q(t,\xi)$ and $Q(t,\xi)^{-1}$ satisfy the linear sde's (4.3) and 
(4.4). \qed

\enddemo

In view of Lemmas 4.1 and 4.2, we can now adapt standard methods from deterministic 
pde's 
in order to prove Theorem
4.1 below.  In particular, the existence of the stochastic semiflow for weak 
solutions of 
the spde (4.1) follows from the
regularity properties of solutions to the random reaction diffusion equation (4.2). 
For 
the existence of the semiflow of (4.2), its global Lipschitz
continuity and compactness, we refer the reader to [Te], pp. 80-102, 371-374.  Note 
that 
Lemma 4.2 ensures that the non-linear
time-dependent random term $\tilde f$ in (4.2) satisfies appropriate dissipativity 
estimates which carry sufficient uniformity 
in $t$ to allow for the apriori estimates in [Te] to work.  This renders the  proof 
of 
Theorem 4.1 below an adaptation of the
corresponding arguments in [Te]. Thus, we will only sketch the proof and leave many 
of the 
details to the reader.

\proclaim{Theorem 4.1}  

Assume that $f$ in (4.1) satisfies Conditions (D). Then for each $\psi \in H:=L^2 
(\Cal 
D)$, the spde (4.1)  admits a unique 
$(\F_t)_{t \geq 0}$-adapted weak solution $u(\cdot,\psi, \cdot): \R^+ \times \O \to 
H$ 
such that 
$u(\cdot,\psi,\o) \in L^{2p} ((0,T), L^{2p}(\Cal D))\cap C(\R^+, H)$ for a.a. $\o 
\in \O$. 
The family of all weak solutions 
of (4.1) has  a
$(\B(\R^+)\otimes \B(H) \otimes \F, \B(H))$-measurable
version  $u:\R^+ \times H\times \Omega \to H$
with the following properties:

\item{(i)}For each $\psi \in H$, $u(\cdot,\psi,\cdot): \R^+ \times \Omega \to H$ is  
an 
$(\F_t)_{t\geq 0}$-adapted weak solution
of (4.1). 
 
\item{(ii)}$(u, \theta)$ is a $C^{0,1}$ perfect cocycle on $H$ (in the sense of
Definition 1.2).
\item{(iii)} For each $(t,\omega)\in (0,\i) \times \Omega$, the map $H \ni
\psi \mapsto u(t,\psi,\omega) \in H$ is globally Lipschitz and takes bounded sets in 
$H$ 
into relatively compact sets.

\item{(iv)} For any positive $a, \rho$,
$$
E \log^+ \sup_{0 \leq t \leq a\atop |\psi|_H \leq \rho} 
|u(t,\psi,\cdot)|_{H} 
< \infty.$$

\item{(v)} For each $\o \in \O$, 
%
$$
\limsup_{t \to \infty} \frac{1}{t} \log^+ \sup_{\psi_1 \neq \psi_2,\atop{ \psi_1, 
\psi_2 
\in H}} \frac{|u(t,\psi_1,\o)-
u(t,\psi_2,\o)|_H}{|\psi_1-\psi_2|_H} \leq 
\frac{1}{2}\biggl (c_4 - \nu \lambda_1-\sigma^2 \biggr )
$$
where $\sigma^2 := \displaystyle \inf_{\xi \in \Cal D} \sum_{i=1}^\i \sigma_i^2 
(\xi) 
$.
In particular, if 
$$\displaystyle \sup_{s \in \R} f'(s) - \nu \lambda_1-\sigma^2 < 0,$$
 then the stochastic flow  $u(t, \cdot, \o): H \to H$ is a uniform contraction for 
sufficiently large $t > 0$.
\endproclaim

\demo{Proof} 

The existence and uniqueness of a weak solution of (4.1) follows from the 
corresponding 
result for the random
reaction-diffusion equation (4.2) ([Te], pp. 89-91).  Using the 
dissipativity estimates (4.5) on $\tilde f$, a straightforward modification of the 
Galerkin approximation technique in [Te] (pp.
89-91) gives the existence of a weak solution $u(\cdot, \psi, \o): \R^+  \to H$ of 
the 
random reaction-diffusion equation (4.2)
for each fixed $\o \in \O$ (cf. also [Ro], pp. 221-227). The joint measurability and 
$(\Cal F_t)_{t \geq 0}$-adaptedness of the
solution are also  immediate consequences of the Galerkin approximations. This 
completes 
the proof of assertion (i) of the
theorem. 

To prove assertion (iii), denote by $v(\cdot,\psi,\cdot): \R^+ \times \O \to H, \,\, 
\psi \in H,$ the family of all weak solutions of the random pde (4.2). We will show 
that 
for each $\o \in \O$, the map 
$H \ni \psi \mapsto u(t,\psi,\o) \in H$ is globally Lipschitz uniformly in $t$ over 
bounded sets in $\R^+$. To see this, 
let $\psi_i \in H, i=1,2$. Denote by $v_i(t):= v(t,\psi_i), t \geq 0, i=1,2$, the 
weak 
solutions of the random pde (4.2) starting at
$\psi_i \in H, i=1,2$. Then multiplying both sides of the equation
$$
 \frac{\partial (v_1(t)-v_2(t))}{\partial t} =\nu Q(t)^{-1}\Delta 
(Q(t)(v_1(t)-v_2(t))) + \tilde 
f(t,v_1(t))- \tilde f(t,v_2(t)), \quad t > 0,  
$$
by $v_1(t)-v_2(t)$ and integrating over $\Cal D$, we obtain  
$$
\align
\int_{\Cal D}(v_1(t)-v_2(t))\frac{\partial (v_1(t)-v_2(t))}{\partial t} \, d\xi 
=&\nu 
<Q(t)^{-1}\Delta (Q(t)(v_1(t)-v_2(t))), v_1(t)-v_2(t)>_H \\
&\qquad +\int_{\Cal D}  \tilde f(t,v_1(t))- \tilde  f(t,v_2(t)) (v_1(t)-v_2(t)) \, 
d\xi, 
\endalign 
$$
for all $t > 0$. Using the Mean-Value Theorem and the second estimate in (4.5), it 
follows 
that  
$$
\align
\frac{1}{2} \frac{d}{dt} |v_1(t)-v_2(t)|^2_H =&\nu <Q(t)^{-1}\Delta 
(Q(t)(v_1(t)-v_2(t))), 
v_1(t)-v_2(t)>_H \\
&\qquad +\int_{\Cal D}  \int_0^1  \frac{\partial \tilde f}{\partial s} (t,\lambda 
v_1(t)+(1-\lambda)v_2(t)) \, d \lambda
(v_1(t)-v_2(t))^2 \, d\xi \\
\leq & -\nu \lambda_1|v_1(t)-v_2(t)|^2_H + c_4 |v_1(t)-v_2(t)|^2_H \tag 4.6
\endalign
$$
for all $t > 0$. In the above inequality, $\lambda_1$ is the smallest eigenvalue of 
$-Q(t)^{-1}\Delta (Q(t)\cdot)$. This turns out to be the same as the smallest 
eigenvalue of $-\Delta$.
Applying Gronwall's 
lemma to (4.6), we get 
$$
|v_1(t,\o)-v_2(t,\o)|^2_H \leq |\psi_1-\psi_2|^2_H \exp \{(c_4 - \nu \lambda_1)t\} 
\tag 4.7
$$
for all $t \geq 0$ and all $\o \in \O$. Using the relations $u(t,\psi_i,\o)= 
Q(t,\xi,\o) 
v_i(t,\o), i=1,2,$ in (4.7), we deduce that 
$$
|u(t,\psi_1,\o)- u(t,\psi_2,\o)|_H \leq |\psi_1-\psi_2|_H \exp \biggl \{
\frac{1}{2}\biggl (c_4 - \nu
\lambda_1 \biggr )t \biggr \}\sup\limits_{\xi\in \Cal D}
Q(t,\xi,\omega) 
$$
for all $t \geq 0, \,\o \in \O, \, \psi_1, \psi_2 \in H$. For any $a > 0$, 
define the random variable
$$
c_5(\o):= \sup_{0 \leq t \leq a}\exp \biggl \{\frac{1}{2}\biggl (c_4 - 
\nu\lambda_1 \biggr )t \biggr \}\sup\limits_{\xi\in \Cal D}
Q(t,\xi,\omega),
\quad \o \in \O.
$$
Then it is easy to see that $E \log^+  c_5 < \infty,$ and 
$$
|u(t,\psi_1,\o)- u(t,\psi_2,\o)|_H \leq c_5(\o) |\psi_1-\psi_2|_H   \tag 4.8
$$
for all $t \in [0,a], \,\o \in \O, \, \psi_1, \psi_2 \in H$. This proves the first 
assertion in (iii). (Note that (4.8) implies 
pathwise uniqueness of the weak solution to the spde (4.1): Just put 
$\psi_1=\psi_2=\psi$, 
a given initial function in $H$.)
The local compactness of the semiflow $H \ni \psi \mapsto u(t, \psi,\o) \in H, \, t 
> 0, 
\o,$ follows from the fact that 
$H \ni \psi \mapsto v(t, \psi,\o) \in H, \, t > 0, \o,$ takes bounded sets in $H$ to 
relatively compact sets.

We next prove the perfect cocycle property in (ii).  To this end, fix $\psi \in H, 
\o \in 
\O, t_1, t_2  \geq 0$. Define 
$$ 
Y(t):= v(t+t_1,\psi,\o), \quad Z(t):= Q(t_1,\o)^{-1} v(t, Q(t_1,\o) v(t_1,\psi,\o), 
\theta 
(t_1,\o))
$$
for all $t \geq 0$. Recall the perfect cocycle identity:
$$
Q(t_1+t_2,\xi,\o) = Q(t_2,\xi,\theta (t_1,\o))Q(t_1,\xi,\o) \quad t_1,t_2 \geq 0,  
\xi 
\in 
\Cal D. \tag 4.9
$$
By the definition of $\tilde f$ in Lemma 4.1 and the above cocycle property,  one 
gets
$$
\align
\tilde f(t,\xi,s,\theta (t_1,\o)) =&Q(t,\xi,\theta (t_1,\o))^{-1} f(Q(t,\xi,\theta 
(t_1,\o))s) \\
=&Q(t_1,\xi,\o) Q(t+t_1,\xi,\o) ^{-1}   f(Q(t+t_1,\xi,\o) Q(t_1,\xi,\o)^{-1} s) \\
=&Q(t_1,\xi,\o) \tilde f (t+t_1,\xi, Q(t_1,\xi,\o)^{-1} s, \o) \tag 4.10
\endalign
$$
for all $s \in \R, t \geq 0, \xi \in \Cal D$.  
 
We now claim that the weak solution of the random reaction-diffusion equation (4.2) 
satisfies the following identity
$$
v(t+t_1,\psi, \o)=Q(t_1,\o)^{-1} v(t, Q(t_1,\o) v(t_1,\psi, \o), \theta (t_1,\o)) 
\tag 
4.11
$$
for all $t \geq 1$. This says that $Y(t)=Z(t)$ for all $t \geq 0$. Using (4.11) and 
the 
relation between $u$ and $v$, it is easy
to check that $(u,\theta)$ is a perfect cocycle. So we need only prove (4.11).  By 
the 
definition of $Z$ and (4.10), 
it follows that 
$$
\aligned
\frac{\partial Z}{\partial t}=& \nu  Q(t_1,\o)^{-1} Q(t,\theta (t_1,\omega))^{-1} 
\Delta (Q(t,\theta (t_1,\omega))v(t, Q(t_1,\o) v(t_1,\psi, \o), 
\theta (t_1,\o)))  \\
&\qquad  + Q(t_1,\o)^{-1} \tilde f(t, v(t, Q(t_1,\o) v(t_1,\psi, \o), \theta 
(t_1,\o)), 
\theta (t_1,\o))\\
=&\nu Q(t+t_1,\omega)^{-1}\Delta (Q(t+t_1,\omega) Z(t)) \\
&+  \tilde f(t+t_1, Q(t_1,\o)^{-1} v(t, Q(t_1,\o) v(t_1,\psi, \o), \theta 
(t_1,\o)), \o) \\
=&\nu Q(t+t_1,\omega)^{-1}\Delta (Q(t+t_1,\omega) Z(t)) +  \tilde f(t+t_1, Z(t), 
\o), \quad t > 0;
\endaligned
$$
 and $Z(0)= v(t_1,\psi,\o)$. Now from its definition, $Y$ also satisfies the same 
random 
pde:
$$
\frac{\partial Y}{\partial t}=\nu Q(t+t_1,\omega)^{-1}\Delta (Q(t+t_1,\omega) Y(t)) 
+  \tilde f(t+t_1, Y(t), \o),\quad t > 0
$$
with the same initial condition $Y(0)= v(t_1,\psi,\o)$. Therefore, by uniqueness of 
weak 
solutions to the above pde,
we must have $Y(t)=Z(t)$ for all $t \geq 0$. This proves our claim, and hence $(u, 
\theta)$ is a perfect cocycle on $H$.

Assertion (v) of the theorem follows easily from (4.8). This completes the proof of 
the 
theorem. \qed

\medskip

Our next result establishes Fr\'echet differentiability of the cocycle generated by 
the reaction diffusion equation:
$$
\left. \aligned
du=&\nu \Delta u \, dt+ (1-|u|^{\alpha})u\, dt+\sum \limits 
_{i=1}^{\infty}\sigma_i(\xi) u\, dW_i(t), \quad t > 0 \\
u(0)=&\psi  \in  H:=L^2 (\Cal D)\\
u(t)|_{\partial \Cal D}=&0, \quad t > 0.
\endaligned \right \} \tag4.12
$$
%
%
%
where $\nu > 0$ is a positive constant and $\Delta$ is the Laplacian on a smooth 
bounded 
domain $\Cal D$ with Dirichlet 
boundary conditions. The result is established under the dimension requirement 
$\alpha < {4\over d}$. It is not clear whether this condition is necessary for 
Fr\'echet differentiability of the cocycle.

\proclaim{Theorem 4.2}


In (4.12), assume that $\alpha < {4\over d}$. 
  Then for each $\psi \in H:=L^2 (\Cal D)$, the spde (4.12)  admits a unique 
$(\F_t)_{t \geq 0}$-adapted weak solution $u(\cdot,\psi, \cdot): \R^+ \times \O \to 
H$ such that 
$u(\cdot,\psi,\o) \in L^{2p} ((0,T), L^{2p}(\Cal D))\cap C(\R^+, H)$ for a.a. $\o 
\in \O$. 
The family of all weak solutions 
of (4.12) has  a $(\B(R^+)\otimes \B(H) \otimes \F, \B(H))$-measurable
version  $u:\R^+ \times H\times \Omega \to H$
with the following properties:

\item{(i)}For each $\psi \in H$, $u(\cdot,\psi,\cdot): \R^+ \times \Omega \to H$ is  
an 
$(\F_t)_{t\geq 0}$-adapted weak solution of (4.12).
\item{(ii)}$(u, \theta)$ is a $C^1$ perfect cocycle on $H$ (in the sense of 
Definition 
1.2).
\item{(iii)} For each $(t,\omega)\in (0,\i) \times \Omega$, the map $H \ni
\psi \mapsto u(t,\psi,\omega) \in H$ is globally Lipschitz and takes bounded sets in 
$H$ 
into relatively compact sets.

\item{(iv)} For each $(t,\psi,\omega) \in (0,\i) \times H
\times \Omega$, the Fr\'echet derivative $D u(t,\psi,\omega) \in L(H)$ is compact, 
and the 
map
$$
[0,\infty) \times H\times \Omega \ni (t,\psi,\omega) \mapsto
D u(t,\psi,\omega) \in L(H)
$$
is strongly measurable.
\item{(v)} For any positive $a, \rho$,
$$
E \log^+ \sup_{0 \leq t \leq a\atop |\psi|_H \leq \rho} \biggl
\{|u(t,\psi,\cdot)|_{H}+ \|D u(t,\psi,\cdot)\|_{L(H)}\biggr \} < \infty.
$$
\endproclaim


\demo{Proof} 

Fix any $\psi \in H=L^2 (\Cal D)$. The existence and uniqueness of the solution to 
(4.12)  in $L^2(\Cal D)$ is well-known as the nonlinear term satisfies the 
dissipativity 
condition ([D-Z.1]). This also follows by a similar argument to the proof of Theorem 
4.1. 
So assertion (i) follows easily. The main
purpose is to prove assertions (ii), (iii) and (iv). Recall that  
$Q(t,\xi):= \exp \biggl \{\sum \limits _{i=1}^{\infty}\sigma_i(\xi) W_i(t)- {1\over 
2}\sum \limits _{i=1}^{\infty}\displaystyle
\sigma_i(\xi)^2 t \biggr \}, \, t \geq 
0, \xi \in \Cal D$, and let
$v(t)=u(t)Q^{-1}(t), t \geq 0.$

For simplicity of notation and till further notice, we will suppress the dependence 
of the random fields $u,v$, etc.. on $\o$.  

Observe that $v(t,\psi)$ is a weak solution of the  random reaction diffusion 
equation
$$
\frac{\partial v}{\partial t} =\nu \Delta v+2\nu Q(t)^{-1}\nabla Q(t)\nabla v+v(\nu 
Q(t)^{-1}\Delta
Q(t)+1-Q^{\alpha}(t)|v|^{\alpha}),\quad t > 0. \tag4.13
$$
By the Feynman-Kac formula, we have 
$$
v(t, \psi)(\xi )=\hat E[\chi_{\tau_t=t} \psi(x_t)e^{\int
_0^t(\nu Q(t-s,x_s)^{-1}\Delta 
Q(t-s,x_s)+1-Q^{\alpha}(t-s,x_s)|v|^{\alpha}(t-s,\psi)(x_s))ds}]
\tag4.14
$$
where $x$ is the solution of the following stochastic differential equation
$$
dx (s)=\sqrt {2\nu} dB(s)+2\nu\nabla \log Q(t-s,x_s)ds, \ \ x_0=\xi  \in \Cal D,
$$
and  $B$ is a Brownian motion in $\R^d$ independent of the $W_i, i \geq 1$. In 
(4.14), 
$\tau_t := \min (\tau, t)$, where 
$\tau$ is the first time the diffusion $x$ hits $\partial \Cal D$. 
%
%
Define $\beta :=\nu \sup\limits_{ 0 \leq s\leq t \leq a,\xi\in \Cal D}{ 
\displaystyle \Delta Q (t-s,\xi)\over Q(t-s,\xi)}$ for any $a > 0$.
It follows from Jensen's inequality and (4.14) that
$$
\align
|v(t,\psi)|_H^2\leq& \int _{\Cal D} (\hat E \chi_{\tau_t=t} 
|\psi(x_t)|e^{(\beta+1)t})^2d\xi \\
\leq &e^{2(\beta+1)t} \int _{\Cal D}\biggl (\int _{\Cal D}p_t(\xi ,y)|\psi(y)|dy 
\biggr 
)^2d\xi \\
\leq &e^{2(\beta+1)t} \int _{\Cal D}\int _{\Cal D}p_t(\xi ,y)(\psi(y))^2dyd\xi \\
\leq &e^{2(\beta+1)t} |\psi|_H^2, \qquad  0 \leq  t  \leq a.\tag4.15
\endalign
$$
In the above inequalities, $p_t(\xi ,y)$ denotes the heat kernel associated with
$\nu \Delta +2\nu(\nabla \log Q(t))\nabla $
on $\Cal D$ with Dirichlet boundary condition. Define the induced heat semigroup 
$T_t 
; H 
\to H, t 
\geq 0$, by 
$$(T_t\psi)(\xi):=\int_\Cal D p_t(\xi,y)\psi(y)dy, \quad \psi \in H, \xi \in \Cal D, 
t 
\geq 0.$$
Note that there exists a constant $c>0$ such that
$$
p_t(\xi ,y)\leq {c\over t^{d\over 2}}, \quad  \xi ,y \in \Cal D, t > 0.
\tag4.16
$$
It is easy to see, using Jensen's inequality and (4.16), that
$$\align
|v(t,\psi)(\xi )|^2\leq &(\hat E \chi_{\tau_t=t}  |\psi(x_t)|e^{(\beta +1)t})^2\\
\leq & e^{2(\beta +1)t}\hat E\chi_{\tau_t=t}  (\psi(x_t))^2\\
\leq &  e^{2(\beta +1)t}\int _{\Cal D} p_t(\xi ,y)(\psi(y))^2dy\\
\leq &  e^{2(\beta +1)t}{c\over t^{d\over 2}}\int _{\Cal D}\psi^2(y)\, dy
\endalign
$$
for all $\xi \in \Cal D$ and $t > 0$. Hence
$$
\|v(t,\psi)\|_{\infty} \leq {\sqrt c \, e^{(\beta +1)t}\over t^{d\over 4}} 
 |\psi|_H \tag4.17
$$
for $t > 0$. 

Now let $\psi,g\in L^2(\Cal D)$ with $|g|_H \leq 1$, and h be a small 
real number. Since $v$ is a mild solution of (4.13), it follows that
$v(t,\psi+hg)-v(t,\psi)$ satisfies the following convolution equation in $H$:
$$
\align
v(t,\psi+hg)-v(t,\psi)=& hT_tg +\int_0^t T_{t-s}(v(s,\psi+hg)-v(s,\psi))(1+{\Delta 
Q(s)\over Q(s)})ds \\
&+\int_0^t T_{t-s}Q^{\alpha}(s)[ 
v(s,\psi)|v(s,\psi)|^{\alpha}-v(s,\psi+hg)|v(s,\psi+hg)|^{\alpha}]ds, t > 0.
\tag4.18
\endalign
$$
Define $m(x):=x|x|^{\alpha}$ for each $x \in \R$.
Then $m^{\prime}(x)=(\alpha +1)|x|^{\alpha}, \, x \in \R$.  By the Mean-Value 
Theorem, 
we 
have 
$$
\aligned
 v(s,\psi+hg)&|v(s,\psi+hg)|^{\alpha}-v(s,\psi)|v(s,\psi)|^{\alpha}\\
=&(\alpha 
+1)\int_0^1|rv(s,\psi+hg)+(1-r)v(s,\psi)|^{\alpha}dr(v(s,\psi+hg)-v(s,\psi))
\endaligned \tag4.19
$$
for all $s \in \R^+$. Combining (4.18) and (4.19) we obtain  
$$
\align
|&v(t,\psi+hg)-v(t,\psi)|_H\\
\leq &h |g|_H+c_1(\omega)\int_0^t |v(s,\psi+hg)-v(s,\psi)|_H\, ds\\
&+c_1(\o) (\alpha +1) \int_0^t (\|v(s,\psi+hg)\|_{L^{\infty}} + 
\|v(s,\psi)\|_{L^{\infty}})^{\alpha} |v(s,\psi+hg)-v(s,\psi)|_H\, ds,
\endalign
$$
for all $0 \leq t\leq a$, where  $c_1(\omega)$ is a positive random constant.
By virtue of (4.17),  we get
$$
\aligned
|v(t,\psi+hg)-v(t,\psi)|_H \leq &
h |g|_H +c_1(\omega)\int_0^t |v(s,\psi+hg)-v(s,\psi)|_H \, ds\\
&+c_2(\o) (\alpha +1) \int_0^t \frac{1}{s^{{d\over 4}\alpha}} 
|v(s,\psi+hg)-v(s,\psi)|_H\, ds,
\endaligned
\tag4.20
$$
for all $0 \leq t\leq a$, where $c_2(\o ) > 0$ is a random constant depending on the 
ball  
$\{g \in 
H: |g|_H \leq 1\}$.  Using Gronwall's lemma and the requirement $\alpha d<4$, we 
obtain 
$$
\sup_{0\leq t\leq a}\, \sup\Sb{ g\in H} \\  |g|_H \leq 1 \endSb
|v(t,\psi+hg)-v(t,\psi)|_H \leq C(\o)h.
$$
Consider the $L(H)$-valued integral  equation:
$$
G_t(\psi)=T_t+\int _0^tT_{t-s}\bigg (Q(s)^{-1}\Delta Q(s)+1-(\alpha 
+1)Q^{\alpha}(s)|v(s,\psi)|^{\alpha}\bigg )G_s(\psi)ds, t \geq 0,
\tag4.21
$$
where $1-(\alpha +1)Q^{\alpha}(s)|v(s,\psi)|^{\alpha}$ is regarded as a 
multiplication operator on $L^2(\Cal D)$
whose operator norm satisfies the inequality
$$
||Q(s)^{-1}\Delta Q(s)+1-(\alpha +1)Q^{\alpha}(s)|v(s,\psi)|^{\alpha}||_{\infty}\leq 
C_a(\o) \frac{1}{s^{{d\over 4}\alpha}}
|\psi|_H^{\alpha}, \quad 0 < s \leq a,
\tag4.22
$$
for a positive random constant $C_a(\o)$.

\noindent
{\it Claim:}  
There exists a unique, continuous solution $[0,\i) \ni t \mapsto G_t(\psi) \in L(H)$ 
to 
equation (4.21).  Moreover for $t>0$, $G_t(\psi): H \rightarrow H$ is compact.

\noindent
{\it Proof of claim}: 
Let $G_t^1(\psi)=T_t, t\geq 0$. Define for $n\geq 1$
$$ G_t^{n+1}(\psi)=T_t+\int _0^tT_{t-s}\bigg (Q(s)^{-1}\Delta Q(s)+1-(\alpha 
+1)Q^{\alpha}(s)|v(s,\psi)|^{\alpha}\bigg )G_s^n(\psi)ds.
\tag4.23
$$
Then (4.22)  and (4.23) imply that
$$
||G_t^{n+1}(\psi)||_{L(H)}\leq 1 + C_a(\o )||\psi||^{\alpha}\int_0^t 
\frac{1}{s^{{d\over 4}\alpha}}||G_s^{n}(\psi)||_{L(H)} ds., \quad 0 < t \leq a
\tag 4.24
$$
Since ${\alpha d\over 4} <1$, then  by the standard successive approximation 
technique it follows that the sequence $\{G^n_{\cdot}(\psi) \}_{n=1}^{\i}$ converges 
to 
the unique solution of (4.12).   Next  we prove that $G_t(\psi)$ is compact for each 
$ 
t > 
0$. 
It suffices  to 
show that a Cauchy sequence can be  extracted from the set 
$\{G_t(\psi)(g): |g|_H\leq 1\}$ for each $ t > 0$.   Let $\delta_m, m\geq 1 $ be a 
sequence of positive numbers decreasing to zero. Since $T_t$ is 
compact for every $t>0$, by a diagonal process there exists a 
sequence $g_n\in H$ with
$|g_n|_H \leq 1$ such that $T_{\delta_m}g_n, n\geq 1$ is a Cauchy 
sequence for every $m$. Since $T_t, t\geq 0$,  is a contraction 
semigroup on $H$, it is easy to see  that   $T_tg_n, n\geq 1$ is 
a Cauchy sequence for every $t>0$. Now consider
$$
\aligned
&
G_t(\psi)(g_n)-G_t(\psi)(g_m)=
T_t(g_n-g_m)
\\
&+\int _0^tT_{t-s}\bigg (Q(s)^{-1}\Delta Q(s)+1-(\alpha 
+1)Q^{\alpha}(s)|v(s,\psi)|^{\alpha}\bigg )(G_s(\psi)(g_n)-G_s(\psi)(g_m))ds, t \geq 
0.
\endaligned
\tag4.25
$$
Hence,
$$
\aligned
|G_t(\psi)&(g_n)-G_t(\psi)(g_m)|_H\\
\leq & |T_tg_n-T_tg_m|_H+C_a(\o ) |\psi|_H^{\alpha}\int_0^t 
\frac{1}{s^{{d\over 4}\alpha}} |G_s(\psi)(g_n)-G_s(\psi)(g_m)|_H \,ds.
\endaligned \tag 4.26
$$
for all $t \in [0,a]$.
Set $l(t):=\displaystyle \limsup_{n,m\rightarrow 
\infty} |G_t(\psi)(g_n)-G_t(\psi)(g_m)|_H, \, 0 \leq t \leq a$. Taking 
$\displaystyle \limsup_{ n,m\rightarrow 
\infty}$ on both sides of (4.26) we obtain
$$
l(t)\leq C_a(\o ) |\psi|_H^{\alpha}\int_0^t \frac{1}{s^{{d\over 4}\alpha}} l(s)ds, 
\quad 0 
< t \leq a.
$$
This implies that $l(t)=0$ for all $t>0$, and 
completes the proof of the claim.

Next we  show that  $v$ is Fr\'echet differentiable and $Dv(t,\psi)=G_t(\psi)$ for 
all 
$t 
\geq 0$.
First we note that by using the Feynman-Kac formula  and a similar 
argument as in the proof of (4.15) and (4.17),  one has
$$
\int _{\Cal D}G_t(\psi)(g)(\xi )d\xi \leq e^{2(\beta+1)t}|g|_H^2,\tag4.27
$$
and
$$
G_t(\psi)(g)(\xi )\leq {\sqrt c \, e^{(\beta+1)t} \over t^{d\over 4}}|g|_H^2, \quad 
0 
< t 
\leq a.
\tag4.28
$$
Denote
$$
\mu _t(\psi,g)=\frac{1}{h}(v(t,\psi+hg)-v(t,\psi))-G_t(\psi)(g), \quad t > 0.
\tag4.29
$$
It is easy to see that $\mu $ satisfies the following integral equation:
$$
\align
\mu _t(\psi,g)=&\int_0^tT_{t-s}(1+{\Delta Q(s)\over Q(s)})\mu _s(\psi,g)ds\\
&- \int_0^tT_{t-s} Q^{\alpha}(s)[{1\over 
h}(v(s,\psi+hg)|v(s,\psi+hg)|^{\alpha}-v(s,\psi)|v(s,\psi)|^{\alpha})]ds\\
& + (\alpha +1) \int_0^t 
T_{t-s}Q^{\alpha}(s)[|v(s,\psi)|^{\alpha}G_s(\psi)(g)]ds, \quad t \geq 0 .\tag4.30
\endalign
$$
Using $m(y)-m(x)=\int_0^1 m^{\prime}(ry+(1-r)x)dr(y-x)$  it follows 
from (4.30) that
$$
\aligned
\mu _t (\psi,g)=&\int_0^tT_{t-s}(1+{\Delta Q(s)\over Q(s)}) \mu _s(\psi,g)ds\\
&- (\alpha +1)\int_0^t 
T_{t-s}Q^{\alpha}(s)[\int_0^1|rv(s,\psi+hg)+(1-r)v(s,\psi)|^{\alpha}dr\mu 
_s(\psi,g)]ds\\
& + (\alpha +1) \int_0^t 
T_{t-s}Q^{\alpha}(s)[\int_0^1(|v(s,\psi)|^{\alpha}-|rv(s,\psi+hg)+(1-r)v(s,\psi)|^{
\alpha}
)dr\\
&\hskip4cm G_s(\psi)(g)]ds, \quad t \geq 0.
\endaligned
\tag4.31
$$
Set $D(t):=\displaystyle \sup_{|g|_H \leq 1} |\mu _t(\psi,g)|_H, \,\, t \geq 0$.  
Using the 
$L^{\infty}$ bound on $v(s,\psi+hg)$ this implies  that for $0 < t \leq a$, one 
has
$$
\aligned
D(t) \leq &C(\omega)\int_0^tD(s)ds+ C(\o) \int_0^t\frac{1}{s^{{d\over 4}\alpha}}D(s) 
ds\\
& + C(\o )\sup_{|g|_H \leq 1} \int_0^t 
|\int_0^1(|v(s,\psi)|^{\alpha}\\
&\hskip2cm 
-|rv(s,\psi+hg)+(1-r)v(s,\psi)|^{\alpha})dr
G_s(\psi)(g)|_H \, ds.
\endaligned
\tag4.32
$$
Again by  Gronwall's lemma, it follows that there is a 
random constant $C(\o)$ such that 
$$
\aligned
D&(t) \\ 
\leq & C(\o)\sup_{|g|_H\leq 1} \int_0^t 
|(\int_0^1(|v(s,\psi)|^{\alpha}-|rv(s,\psi+hg)+(1-r)v(s,\psi)|^{\alpha})dr)G_s(\psi)
(g)|_H
\, ds.
\endaligned
\tag 4.33
$$
for all $t \in [0,a]$.
To complete the proof of assertions (ii) and (iv), it suffices to show that
$$
\lim_{h\rightarrow 0}\sup_{|g|_H \leq 1} \int_0^t 
|(\int_0^1(|v(s,\psi)|^{\alpha}-|rv(s,\psi+hg)+(1-r)v(s,\psi)|^{\alpha})dr)G_s(\psi)
(g)|_H
\, ds =0   \tag 4.34
$$
for all $t \in [0,a]$. Let us prove (4.24) for $\alpha \leq 1$ and $\alpha >1$ 
separately. 
Assume first $\alpha \leq 1$. By H\"o{}lder inequality,
$$
\aligned
&
\sup_{|g|_H\leq 1} \int_0^t 
|(\int_0^1(|v(s,\psi)|^{\alpha}-|rv(s,\psi+hg)+(1-r)v(s,\psi)|^{\alpha})dr)G_s(\psi)
(g)|_H 
\, ds\\
\leq & \sup_{|g|_H \leq 1} \int_0^t 
|(|v(s,\psi)-v(s,\psi+hg)|^{\alpha})G_s(\psi)(g)|_H \, 
ds \\
\leq & \sup_{|g|_H\leq 1} \int_0^t 
\|G_s(\psi)(g)\|_{L^{\infty}}^{\alpha} 
(|v(s,\psi)-v(s,\psi+hg)|_H^{\alpha})|G_s(\psi)(g)|^{1-\alpha}||ds
\\
\leq &\sup_{|g|_H \leq 1} \int_0^t \|G_s(\psi)(g)\|_{L^{\infty}}^{\alpha} 
|v(s,\psi)-v(s,\psi+hg)|_{H}^{\alpha}
|G_s(\psi)(g)|_H^{1-\alpha}\, ds.
\endaligned
$$
By virtue of (4.20) and (4.28), we get
$$
\aligned
\sup_{|g|_H \leq 1}& \int_0^t 
|(\int_0^1(|v(s,\psi)|^{\alpha}-|rv(s,\psi+hg)+(1-r)v(s,\psi)|^{\alpha})dr)G_s(\psi)
(g)||d
s\\
\leq &C_{\alpha}(\o)h^{\alpha} \int_0^t \frac{1}{s^{{d\over 4}\alpha}} ds
=C_{\alpha}(\omega){1\over 1-{\alpha d\over 4}}t^{1-{\alpha d\over 
4}}h^{\alpha}, \quad 0 < t \leq a,
\endaligned
$$
where $C_{\alpha}(\o)$ is a random constant depending on the set 
$\{g\in L^2(\Cal D): ||g||_{L^2} \leq 1\}$.
This implies (4.24). 

Assume now that $\alpha >1$.  Then
$$
\aligned
&\int_0^t 
|(\int_0^1(|v(s,\psi)|^{\alpha}-|rv(s,\psi+hg)+(1-r)v(s,\psi)|^{\alpha})dr)
G_s(\psi)(g)|_H \, ds\\
\leq & \int_0^tds\int_0^1\, dr 
|(|v(s,\psi)|^{\alpha}-|rv(s,\psi+hg)+(1-r)v(s,\psi)|^{\alpha})G_s(\psi)(g)|_H 
\\ \leq & \int_0^tds \int_0^1dr\int_0^1dk\alpha 
|(k|v(s,\psi)|+(1-k)|rv(s,\psi+hg)+(1-r)v(s,\psi)|)^{\alpha -1}\\
&\hskip2cm
|v(s,\psi+hg)-v(s,\psi)| G_s(\psi)(g)|_H
\\
\leq & \int_0^tds\int_0^1dr\int_0^1d k\alpha 
\|(k|v(s,\psi)|+(1-k)|rv(s,\psi+hg)+(1-r)v(s,\psi)|)^{\alpha 
-1}\|_{L^{\infty}}\\
&\hskip2cm
||G_s(\psi)(g)||_{L^{\infty}} |v(s,\psi+hg)-v(s,\psi)|_H, \quad 0 < t \leq a .
\endaligned
\tag4.35
$$
By (4.17), (4.20) and (4.28) it follows from (4.35) that
$$
\sup_{|g|_H \leq 1} \int_0^t 
|(\int_0^1(|v(s,\psi)|^{\alpha}-|rv(s,\psi+hg)+(1-r)v(s,\psi)|^{\alpha})dr)
G_s(\psi)(g)|_H\, ds
$$
$$
\leq C_{\alpha}(\o) h \int_0^t \frac{1}{s^{{d\over 
4}(\alpha-1)}}\frac{1}{s^{{d\over 4}}}ds, 0 < t \leq a.
\tag4.36
$$
This implies  (4.34). So assertion (iv) holds.

To establish (iii), use (4.18), (4.19) and a similar argument to the proof of 
(4.20), to 
obtain the following inequality 
$$
\aligned
|v(t,\psi_m)-v(t,\psi_n)|_H \leq &
|T_t\psi_m-T_t\psi_n|_H+C_1(\omega)\int _0^t |v(s,\psi_m)-v(s,\psi_n)|_H\,ds\\
&+C_2(\omega)(\alpha+1)\int _0^t{1\over {s^{{d\over 4}\alpha}}}
|v(s,\psi_m)-v(s,\psi_n)|_H\, ds, \quad 0 < t \leq a,
\endaligned
\tag4.37
$$
for  $\psi_n, \psi_m \in H$ such that $|\psi_m|_H, |\psi_n|_H \leq 1$.
As  in the proof of the compactness of $Dv(t,\psi)$, 
we can select a subsequence denoted also by 
$\{\psi_{n}\}\subset \{\psi: |\psi|_H \leq 1\}$ such that for each $t>0$, 
$|T_t\psi_n-T_t\psi_m|_H \to 0$ as $n,m\to \infty$. One then
can prove from
(4.32) that $\lim\limits _{n,m\to \infty} |v(t,\psi_n)-v(t,\psi_m)|_H=0$.  
Therefore $v(t, \cdot): H \to  H$ is compact or each $t > 0$..
This implies the compactness of each Fr\' echet derivative $Du(t, \psi,\o): H \to H, 
\, t > 
0, \o \in 
\O.$  Hence the first assertion in (iv) holds.

To prove the strong measurability assertion in (iv), we now highlight the dependence 
of $u$ 
on $\o$. 
Note first that the map 
$$
[0,\infty) \times H\times \Omega \ni (t,\psi,\omega) \mapsto  u(t,\psi,\omega) \in H
$$
is jointly measurable. This is a consequence of the (uniform) continuity of 
$$
[0,a] \times H \ni (t,\psi) \mapsto  u(t,\psi,\omega) \in H, \quad  \o \in 
\O,
$$
and the measurability of 
$$
 \Omega \ni \omega \mapsto  u(t,\psi,\omega) \in H, \quad (t,\psi) \in \R^+ \times 
H.
$$
Secondly, the joint strong measurability of 
$$
[0,\infty) \times H\times \Omega \ni (t,\psi,\omega) \mapsto
D u(t,\psi,\omega) \in L(H)
$$
follows from the relation
$$
D u(t,\psi,\omega)(\eta)= \lim_{h \to 0} \frac{1}{h} [u(t, \psi + h \eta, 
\o)-u(t,\psi,\o)], \quad (t,\omega) \in \R^+ \times \O,
\psi, \eta \in H. \tag 4.38
$$
Finally, note that the integrability estimate in (v) follows from the Lipschitz 
property 
of $u(t,\cdot,\o): H \to H, \, (t,\o) \in \R^+ \times \O$. In particular, (4.38) and 
the 
above Lipschitz property give 
$$\|D u(t,\psi,\omega)\|_{L(H)}  \leq c_5(\o)$$
for all $(t,\psi,\o) \in [0,a] \times H \times \O$, with $E \log^+  c_5 < \infty$.
\qed

\enddemo

\remark{Remarks}

\item{(i)} It is easy to see that above proof is also valid for the initial
boundary value problem
with  Neumann
boundary condition. Note the exact formula of the heat kernel was not
needed in the proof.
Only estimates such as (4.6) and (4.7) were actually needed. These
kind of estimate holds
for Laplacian operator on a bounded domain with smooth boundary and
Neumann boundary condition.
The generalized solution of (4.1) can be defined following Freidlin [Fr]:
$$
\aligned
u(t,&\psi)(\xi )\\
=&\hat E \biggl [\psi(x_t^*)e^{\int
_0^t({\Delta (t-s,x_s^*)\over Q(t-s,x_s^*)}+1-|u|^{\alpha}(t-s,\psi)(x_s^*))ds 
 -{1\over 2} \sum \limits_{i=1}^{\infty}\int _0^t \sigma_i^2(x_s^*)ds+\sum\limits 
_{i=1}^{\infty}\sigma _i(x_s^*)dW_i(t-s)} \biggr ]
\endaligned
$$
a.s. Here $x_t^*$ is a diffusion process starting at $\xi \in \Cal D$ with 
reflection 
on the boundary $\partial \Cal D$ generated with the operator $\nu \Delta+2\nu 
\nabla \log Q(t)\nabla $.
One can see that the analysis in the proof 
of Theorem 4.1 carries through for this case as well.

\item{(ii)} The dimension restriction is used only to guarantee the Fr\'echet 
differentiability 
of the semiflow in Theorem 4.2. This condition is not needed for the existence of 
the 
globally Lipschitz 
flow in Theorem 4.1. The  conditions in Theorem 4.2 are stronger than those in 
Theorem 
4.1, 
and 
accordingly the result.

\endremark
 
\bigskip

\noindent
{\it (b) Burgers equation with additive noise}

 The stochastic Burgers equation has been considered intensively by many researchers 
in recent years ([B-C-J], [B-C-F],  [D-T-Z], [D-Z.2], [D-D-T],
[E-V], [H-L-O-U-Z], [Si], [T-Za], [T-Z]). Here we consider the following stochastic 
Burgers
equation on the interval $[0,1]$,
$$
\aligned
du+u{\partial  u\over \partial \xi}\, dt=&\nu \Delta udt+dW(t), \quad t > 0,\\
u(0,\psi)(\xi)=&\psi(\xi),\\
u(t,\psi)(0)=&u(t,\psi)(1)=0
\endaligned
\tag4.39
$$
where the viscosity $\nu$ is a positive constant. 
Here $W$ is an infinite dimensional Brownian motion in $L^2[0,1]$ on 
a probability space $(\Omega, {\Cal F}, P)$:
$$
W(t):=\sum _{k=1}^{\infty} \sqrt {\lambda _k}e_{k}W^{k}(t).\tag4.40
$$
In (4.40), each $e_k$ is an eigenfunction of $-\nu \Delta$ associated with its eigenvalue
$\alpha _k$, $\{W^k: k=1,2,\cdots \}$ are mutually independent one dimensional 
Brownian motions and
$$
\sum _{k=1}^{\infty} 
{\lambda _k\over \alpha _k}<\infty.\tag4.41
$$

Following ([D-Z.2], pp. 260-265), we  will transform the mild solution of the stochastic 
Burgers equation (4.39) to that  of the  random Burgers equation (4.42) below.
  Let $T_t: L^2([0,1]) \to L^2([0,1]), \, t \geq 0,$ be the  heat semi-group on $[0,1]$ with Dirichlet boundary conditions.
Let
$$
W_p(t)(\xi ):=\int _0^tT_{t-s}dW(s)(\xi).
$$
Then $W_p(t)$ has an $C[0,1]$-valued version with H\"older 
continuous paths ([D-Z.2], Theorem 14.2.4).
Set
$$
v(t,\xi ):=u(t,\xi)-W_p(t)(\xi ), \quad t > 0, \xi \in [0,1].
$$
Then $v(t,\xi )$ is a mild solution of the following equation 
$$
{\partial v\over \partial t} =\nu \Delta
v-{1\over 2}{\partial\over \partial \xi }(v+W_p(t)(\xi ))^2,
\tag4.42
$$
 in the sense of  ([D-Z.2], pp. 260-265).
 
Viewing  equation (4.42) as a random Burgers equation, it is not hard to see 
that, for each initial $\psi \in L^2([0,1])$, it has a unique global solution
$v(\cdot,\psi, \o) \in C(\R^+, L^2([0,1]))\cap L^2([0,1],H^1[0,1])$ 
for each $\o \in \O$; and  for any  $a \in \R^+$ 
and  any bounded set $S \subset L^2([0,1])$, the following holds
 $$\sup_{t \in [0,a] \atop \psi \in S }||v(t,\psi,\o)||_{L^2([0,1])} < \infty  
\tag4.43
$$
for all $\o \in \O$ (cf. [Ta], Chapter 15, Proposition 1.3; [D-Z.2], pp. 260-265).

 A continuous semi-flow for a stochastic Burgers equation with skew-symmetric 
 noise was obtained in [B-C-F].  However, this is not 
sufficient for our purposes, since we 
seek  to construct random  families of {\it differentiable} stable/unstable
 manifolds near 
hyperbolic stationary solutions of (4.39).
In the following theorem, we establish the existence of a perfect $C^1$ compacting 
cocycle for (4.39). In 
Part II of this work ([M-Z-Z]), this fact will enable us  
to use multiplicative ergodic theory techniques in order  to prove  a local 
stable/unstable 
manifold theorem 
near stationary solutions of the stochastic Burgers equation (4.39).
 
\proclaim{Theorem 4.3}

Consider the stochastic Burgers equation (4.39) with $L^2[0,1]$-valued Brownian 
(4.40). 
Then equation (4.39)     has a mild solution with  a
$(\B(\R^+)\otimes \B(L^2([0,1])) \otimes \F, \B(L^2([0,1])))$-measurable
version  $u:\R^+ \times L^2([0,1]) \times \Omega \to L^2([0,1])$
having the following properties:

\item{(i)}For each $\psi \in L^2([0,1])$, $u(\cdot,\psi,\cdot): \R^+
\times \Omega \to L^2([0,1])$ is  $(\F_t)_{t\geq 0}$-adapted.
\item{(ii)}$(u, \theta)$ is a $C^1$ perfect cocycle on $L^2([0,1])$
(in the sense of
Definition 1.2).
\item{(iii)} For each $(t,\omega)\in (0,\i) \times \Omega$, the map 
$L^2([0,1])\ni \psi\mapsto
u(t,\psi,\omega)\in L^2([0,1])$ takes bounded sets into relatively compact sets.

\item{(iv)} For each $(t,\psi,\omega) \in (0,\i) \times L^2([0,1])
\times \Omega$, the Fr\'echet
derivative $D u(t,\psi,\omega) \in L(L^2([0,1]))$ is compact. Furthermore,  the map 
$$
[0,\infty) \times L^2([0,1])\times \Omega \ni (t,\psi,\omega) \mapsto
D u(t,\psi,\omega) \in L(L^2([0,1]))
$$
is strongly measurable. 
\item{(v)} For any positive reals $a, \rho$,
$$
E \log^+ \sup_{0 \leq t \leq a\atop \|\psi\|_{L^2([0,1])} \leq \rho} \biggl
\{\|u(t,\psi,\cdot)\|_{L^2([0,1])}+ \|D u(t,\psi,\cdot)\|_{L(L^2([0,1]))}\biggr \} < \infty.
$$
\endproclaim

\demo{Proof} 

For simplicity of notation, we will assume throughout this proof that 
$\nu=\displaystyle  \frac{1}{2}$.

Assertion (i) follows easily from the global existence of solutions to (4.42).

To prove (ii), consider  Burgers equation (4.42).  Denote by $p_t(\xi ,y)$ the heat kernel for the 
Laplacian $\nu \Delta$ on $[0,1]$ with Dirichlet boundary conditions.
 Recall that there are positive constants $c_1, c_2$ such that
$$
\biggl |{\partial p_t(\xi,y)\over \partial y} \biggr | \leq {c_1\over  {t}}
 {\text {\rm e}}^{-{(\xi-y)^2\over 2c_2t}} \tag{4.44}
$$
for all $t > 0, \, \xi, y \in [0,1]$ (c.f. [L-S-U], p. 413). Then pick a positive 
constant 
$c_3$ such that 
$$
\int _0^1{c_3\over {\sqrt { t}}}{\text {\rm e}}^{-{y^2\over 2c_2t}}\,
dy \, \leq \, 1 
$$
for all $ t > 0$.

Using (4.42), variation of parameters, and integration by parts, we get  
$$
\align
v(t,\psi)(\xi )=&T_t\psi(\xi )-{1\over 2}\int _0^tT_{t-s}\nabla
v^2(s,\psi)(\xi )ds\\
&
+\int _0^tT_{t-s}(-\nabla (W_p(s) v(s,\psi))-W_p(s)\nabla W_p(s))(\xi )ds
\\
=&\int _0^1 p_t(\xi ,y)\psi(y)dy-{1\over 2}\int _0^t\int_
0^1p_{t-s}(\xi ,y)\nabla
v^2(s,\psi)(y)dyds\\
& +\int _0^t\int _0^1p_{t-s}(\xi ,y)(-\nabla (W_p(s)v(s,\psi))-W_p(s)\nabla 
W_p(s))(y)dyds\\
=&\int _0^1 p_t(\xi ,y)\psi(y)dy+{1\over 2}\int _0^t\int_0^1\nabla p_{t-s}(\xi ,y)
v^2(s,\psi)(y)dyds\\
& +\int _0^t\int _0^1\nabla p_{t-s}(\xi ,y)(W_p(s)v(s,\psi)+{1\over 
2}W_p(s)^2)(y)dyds
\endalign
$$
for all $ t \geq 0$.
Thus
$$
\align
v(t,\psi+&hg)(\xi )-v(t,\psi)(\xi )\\
=&h\int _0^1 p_t(\xi ,y)g(y)dy+{1\over 2}\int _0^t\int_0^1\nabla p_{t-s}(\xi ,y)
(v^2(s,\psi+hg)(y)-v^2(s,\psi)(y))dyds\\
& +\int _0^t\int _0^1\nabla p_{t-s}(\xi 
,y)(W_p(s)(v(s,\psi+hg)-v(s,\psi))(y)dyds, \quad t > 0.
\endalign 
$$
 %
 Squaring both sides of the above equality and integrating with respect $\xi \in 
[0,1]$, we obtain
 %
$$
\allowdisplaybreaks
\align
||v(t,&\psi+hg)-v(t,\psi)||_{L^2([0,1])}^2\\
\leq &3h^2||g||_{L^2([0,1])}^2+{3\over 4}\int _0^1 \biggl (\int _0^t\int _0^1\nabla 
p_{t-s}(\xi,y)(v^2(s,\psi+hg)(y)-v^2(s,\psi)(y))dyds \biggr )^2d\xi\\
&+3\int _0^1 \biggl (\int _0^t\int _0^1 \nabla p_{t-s}(\xi 
,y)W_p(s)(y)(v(s,\psi+hg)(y)-v(s,\psi)(y))dyds \biggr )^2d\xi\\
\leq &
3h^2||g||_{L^2([0,1])}^2\\
&+{3\over 4}\int _0^1 \biggl (\int _0^t{1\over \sqrt {t-s}}\int _0^1
{c_1\over {\sqrt { t-s}}}{\text {\rm e}}^{-{(\xi-y)^2\over 2c_2(t-s)}}
(v^2(s,\psi+hg)(y)-v^2(s,\psi)(y))dyds \biggr )^2d\xi
\\
&+3\int _0^1 \biggl (\int _0^t{1\over \sqrt {t-s}}
\int _0^1 {c_1\over {\sqrt { t-s}}}{\text {\rm e}}^{-{(\xi-y)^2\over 2c_2(t-s)}}
W_p(s)(y)\times \\
&\hskip4cm  \times (v(s,\psi+hg)(y)-v(s,\psi)(y))\,dyds \biggr )^2d\xi,
\endalign
$$
for all $t > 0$.  Now use Cauchy-Schwartz inequality, the heat kernel estimate 
(4.44)  
and Fubini's theorem to obtain  
$$
\align
||v(t,\psi+hg)-&v(t,\psi)||_{L^2([0,1])}^2
\leq  3h^2||g||_{L^2([0,1])}^2 +{3\over 4}\int _0^1\int _0^t{1\over {(t-s)}^{3\over 
4}}ds 
\int _0^t{1\over {(t-s)}^{1\over 4}}\\
&\biggl (\int _0^1
{c_1\over {\sqrt { t-s}}}{\text {\rm e}}^{-{(\xi-y)^2\over 2c_2(t-s)}}
(v^2(s,\psi+hg)(y)-v^2(s,\psi)(y))dy \biggr )^2dsd\xi
\\
&+3\int _0^1\int _0^t{1\over {(t-s)}^{3\over 4}}ds 
\int _0^t{1\over {(t-s)}^{1\over 4}}\\
& \biggl (\int _0^1 {c_1\over {\sqrt { t-s}}}{\text {\rm e}}^{-{(\xi-y)^2\over 
2c_2(t-s)}}
W_p(s)(v(s,\psi+hg)-v(s,\psi))(y)dy \biggr )^2dsd\xi \\
\leq  & 3h^2||g||_{L^2([0,1])}^2
+Ct^{1\over 4}\int _0^t{1\over {(t-s)}^{3\over 4}} 
||v(s,\psi+hg)-v(s,\psi)||_{L^2([0,1])}^2ds\\
&+Ct^{1\over 4}\int _0^t{1\over {(t-s)}^{1\over 
4}}||v(s,\psi+hg)-v(s,\psi)||_{L^2([0,1])}ds, \quad t > 0
\endalign
$$
for all $t > 0$  and some positive (random) constant $C$.  
Iterating the above computation, we get
%
$$
\allowdisplaybreaks
\align
||v(t,\psi+hg)&-v(t,\psi)||_{L^2([0,1])}^2\\
\leq  &3h^2||g||_{L^2([0,1])}^2+Ct^{1\over 4}
\int _0^t \biggl ({1\over (t-s)^{3\over 4}}+{1\over (t-s)^{1\over 4}} \biggr )
3h^2||g||_{L^2([0,1])}^2ds\\
& +C^2t^{1\over 4}\int _0^t \int _0^s s^{1\over 4} \biggl ({1\over
(t-s)^{3\over 4}}+{1\over (t-s)^{1\over 4}} \biggr )
\biggl ({1\over (s-r)^{3\over 4}}+{1\over (s-r)^{1\over 4}} \biggr ) \times
\\
&\hskip2cm \times ||v(r,\psi+hg)-v (r,\psi)||_{L^2([0,1])}^2drds, \quad t > 0.
\endalign
$$
Consider now the elementary estimate
$$
\int _r^t {s^{\alpha }\over (t-s)^{\beta}(s-r)^{\gamma}}
ds=\int _0^{t-r} {(s+r)^{\alpha }\over (t-r-s)^{\beta}s^{\gamma}}
ds\leq \frac{C_1}{(t-r)^{\beta + \gamma-1}}, \quad t \geq r > 0,
$$
which holds for any $\alpha \geq 0$, $0\leq \beta <1$,
$0\leq \gamma <1$, and where $C_1>0$ is a positive (deterministic) constant. 
Using the above estimate together with Fubini's theorem, gives
$$
\align
||v(t,\psi+&hg)-v(t,\psi)||_{L^2([0,1])}^2  \\
&\leq 3h^2||g||_{L^2([0,1])}^2
+3C(4t^{1\over 2}+{4\over 3} t)h^2||g||_{L^2([0,1])}^2\\
& +C^2t^{1\over 4}\int _0^t
||v(r,\psi+hg)-v(r,\psi)||_{L^2([0,1])}^2 \times \\
&\hskip1cm \times \int _r^t  s^{1\over 4} \biggl ({1\over
(t-s)^{3\over 4}}+{1\over (t-s)^{1\over 4}} \biggr )
\biggl ({1\over (s-r)^{3\over 4}}+{1\over (s-r)^{1\over 4}} \biggr )dsdr\\
\leq & 3h^2||g||_{L^2([0,1])}^2+3C(4t^{1\over 2}+{4\over 3}
t)h^2||g||_{L^2([0,1])}^2\\
& +C_2 (t+t^{1\over 4})\int _0^t {1 \over (t-r)^{1/2} }
||v(r,\psi+hg)-v(r,\psi)||_{L^2([0,1])}^2dr
\tag4.45
\endalign
$$
for all $t \in (0,a], \, a \in \R^+$. Iterating the above process once more and 
applying Gronwall's lemma, we obtain   
$$ \sup_{0 \leq t \leq a, g \in L^2([0,1])\atop \|g\|_{L^2} \leq 
1}||v(t,\psi+hg)-v(t,\psi)||_{L^2([0,1])}^2\leq M h^2,
\tag4.46
$$
for any $a \in \R^+$, where $M$ is a positive random constant depending on $a$.

For fixed $\psi, g \in L^2([0,1])$,  define 
$G := G(t,\psi)(g)(\xi), \, \,t > 0, \xi \in [0,1],$ to be the weak solution of the
``linearized" Burgers equation 
$$
{\partial G\over \partial t} + 
{\partial(v(t,\psi)G) \over \partial\xi }
={1\over 2}\Delta G-{\partial (W_p G) \over \partial \xi },\  \ G(0,\psi)(g)=g \in 
L^2([0,1]).
$$
Set
$$
\mu _t(\psi,g):=v(t,\psi+hg)-v(t,\psi)-hG(t,\psi)(g), \quad |h| < 1, t \geq 0.
$$
Then it is easy to see that
$$
\allowdisplaybreaks
\align
\mu _t(\psi,g)(\xi )=&-\int _0^t\int _0^1p_{t-s}(\xi ,y)({1\over 2}\nabla
(v(s,\psi+hg)(y)-v(s,\psi)(y))^2\\
&\hskip1cm +
\nabla (v(s,\psi)(y)\mu _s(\psi,g)(y))+\nabla (W_p(s)(y)\mu _s(\psi,g)(y)))dyds
\\
=&\int _0^t\int _0^1\nabla p_{t-s}(\xi ,y)({1\over 
2}(v(s,\psi+hg)(y)-v(s,\psi)(y))^2\\
&\hskip1cm +
v(s,\psi)(y)\mu _s(\psi,g)(y)+ W_p(s)(y)\mu _s(\psi,g)(y))dyds.
\endalign
$$
for all $t > 0$. So using the Cauchy-Schwartz inequality and (4.46), we obtain
$$
\allowdisplaybreaks
\align
||\mu _t(&\psi,g)||_{L^2([0,1])}^2 \\
\leq &{3\over 4}\int _0^1(\int _0^t\int _0^1\nabla p_{t-s}(\xi 
,y)(v(s,\psi+hg)(y)-v(s,\psi)(y))^2dyds)^2d\xi\\
&+3\int _0^1(\int _0^t\int _0^1\nabla p_{t-s}(\xi ,y)(v(s,\psi)(y)\mu 
_s(\psi,g)(y)dyds)^2d\xi\\
&+3\int _0^1(\int _0^t\int _0^1\nabla p_{t-s}(\xi ,y)( W_p(s)(y)\mu 
_s(\psi,g)(y))dyds)^2d\xi
\\
\leq &{3\over 4}\int _0^1\int _0^t{1\over (t-s)^{3\over 4}}ds\int _0^t{1\over 
(t-s)^{1\over 4}}\biggl (\int _0^1{c_1\over \sqrt {t-s}}
e^{-{(y-\xi)^2\over c_2(t-s)}}\\
&\hskip5cm (v(s,\psi+hg)(y)-v(s,\psi)(y))^2dy \biggr )^2dsd\xi\\
&+3
\int _0^1\int _0^t{1\over (t-s)^{3\over 4}}ds\int _0^t{1\over (t-s)^{1\over 4}}
\biggl 
(\int 
_0^1{c_1\over \sqrt {t-s}} e^{-{(y-\xi)^2\over
c_2(t-s)}}\\
&
\hskip5cm (v(s,\psi)(y)\mu _s(\psi,g)(y)dy \biggr )^2dsd\xi\\
&+3\int _0^1\int _0^t{1\over (t-s)^{3\over 4}}ds\int _0^t{1\over (t-s)^{1\over 
4}}\biggl (\int _0^1{c_1\over \sqrt {t-s}}
e^{-{(y-\xi)^2\over c_2(t-s)}}\\
&
\hskip5cm ( W_p(s)(y)\mu _s(\psi,g)(y))dy \biggr )^2dsd\xi.
\endalign
$$
for all $t > 0$. Thus
$$
\allowdisplaybreaks
\align
 ||\mu _t(\psi,g)||_{L^2([0,1])}^2 \leq &
{3c_1}t^{1\over 4}\int _0^t{1\over (t-s)^{3\over 4}} \biggl (\int 
_0^1(v(s,\psi+hg)(y)-v(s,\psi)(y))^2dy \biggr )^2ds\\
&+12c_1t^{1\over 4}
\int _0^t{1\over (t-s)^{3\over 4}}\int _0^1v^2(s,\psi)(y)dy\int _0^1\mu^2 
_s(\psi,g)(y)dyds\\
&+12 t^{1\over 4}||b||_{\infty}^2 \int _0^t{1\over (t-s)^{1\over 4}}\int _0^1\mu 
^2_s(\psi,g)(y)dyds\\
\leq & Ch^4+C
\int _0^t({1\over (t-s)^{3\over 4}}+{1\over (t-s)^{1\over 4}})||\mu 
_s(\psi,g||_{L^2([0,1])}^2ds, \quad 0 \leq t \leq a,
\endalign
$$
where  $C=C(\o,a)$ is a positive random constant.  Using the previous iteration  
argument 
followed by Gronwall's lemma,
we obtain the following estimate
$$
\sup_{  {g \in L^2([0,1]), ||g||_{L^2} \leq 1} \atop 0 \leq t \leq a} ||\mu 
_t(\psi,g)||_{L^2([0,1])}^2\leq M_1h^4, \quad |h| < 1
$$
for some positive random constant $M_1=M_1(\o,a)$.  This implies that  
$\displaystyle {\mu _t (\psi, g) \over h}$ converges to $0$ as $h \to 0$  in
$L^2([0,1])$, uniformly in $(t,g) \in [0,a] \times \{g \in L^2 ([0,1]): ||g||_{L^2} 
\leq 1 
\}$. Therefore,
$\displaystyle {v(t,\psi+hg)-v(t,\psi)\over h}\to G_t(\psi,g)$ as $h \to 0$,
uniformly for $g\in \{g;||g||_{L^2([0,1])}\leq 1\} $. Hence,  
$v$ is Fr\' echet differentiable at $\psi \in L^2 ([0,1])$, with Fr\'echet 
derivative  
$Dv(t,\psi): L^2([0,1])\to L^2([0,1])$ satisfying the  
$L(L^2([0,1]))$-valued linear equation
$$
\align
  Dv(t,\psi)
=&T_t-\int _0^t T_{t-s} \biggl ({\partial v(s,\psi) \over \partial  
\xi} Dv(s,\psi)+v(s,\psi){\partial  Dv(s,\psi)\over\partial \xi} \biggr )ds
\\
&-\int _0^t T_{t-s}\biggl (Dv(s,\psi){\partial W_p(s) \over \partial \xi} 
+W_p(s){\partial 
Dv(s,\psi) \over \partial \xi} \biggr )ds
\tag4.47
\endalign
$$
for $(t,\psi) \in \R^+ \times L^2([0,1])$.

In order to complete the proof of assertion (ii) of the theorem, it remains to prove 
that 
$(u,\theta)$ is a perfect cocycle in $L^2([0,1])$. It is
easy to see from (4.39) that
$$
u(t,\psi)(\omega)=T_t\psi-\int _0^tT_{t-s}u(s,\psi,\o)\nabla u(s,\psi,\o)\,ds
+ \biggl [\int_0^tT_{t-s}dW(s) \biggr ](\o)
$$
for $t > 0, \o \in \O, \psi \in L^2([0,1])$. We need to prove that
$$
u(t,u(t_1,\psi,\omega),\theta(t_1,\omega))=u(t+t_1,\psi,\omega). \tag4.48
$$
for $t, t_1 \geq 0, \o \in \O, \psi \in L^2([0,1])$. To see this, fix 
$t_1 \geq 0, \o \in \O, \psi \in L^2([0,1]),$ and denote
$$
Y(t):=u(t,u(t_1,\psi,\omega),\theta(t_1,\omega)), \ \ Z(t):=u(t+t_1,\psi,\omega), 
\quad t > 
0.
$$
Then
$$
\align
Y(t)=&T_tu(t_1,\psi,\omega)-\int
_0^tT_{t-s}u(s,u(t_1,\psi,\omega),\theta(t_1,\omega))\frac{\partial 
u(s,u(t_1,\psi,\omega),\theta(t_1,\omega))}{\partial y}ds\\
&+\biggl [\int _{t_1}^{t+t_1}T_{t+t_1-s}dW(s)\biggr ](\o)\\
=&T_{t+t_1}\psi-\int _0^{t_1}T_{t+t_1-s}u(s,\psi,\omega)\frac{\partial
u(s,\psi,\omega)}{\partial y} ds\\
&-\int _{t_1}^{t+t_1}T_{t+t_1-s}Y(s-t_1)\frac{\partial
Y(s-t_1)}{\partial y} ds
+\biggl [\int _{0}^{t+t_1}T_{t+t_1-s}dW(s)\biggr ](\o), \quad t > 0.
\endalign
$$
Also,
$$
\align
Z(t)=&T_{t+t_1}\psi-\int _0^{t_1}T_{t+t_1-s}u(s,\psi,\omega)\frac{\partial
u(s,\psi,\omega)}{\partial y} ds\\
&-\int _{t_1}^{t+t_1}T_{t+t_1-s}Z(s-t_1)\frac{\partial
Z(s-t_1)}{\partial y} ds
+\biggl [\int _{0}^{t+t_1}T_{t+t_1-s}dW(s)\biggr ](\o), \quad t > 0.
\endalign
$$
Therefore,
$$
\align
Y(t)(\xi)-Z(t)(\xi)
=&-\frac{1}{2}\int _{t_1}^{t+t_1}T_{t+t_1-s} \biggl (\frac{\partial 
Y^2(s-t_1)}{\partial 
y}-\frac{\partial Z^2(s-t_1)}{\partial y} \biggr )(\xi)ds\\
=&\frac{1}{2}\int _{t_1}^{t+t_1}ds\int_0^1\frac{\partial 
p_{t+t_1-s}(\xi,y)}{\partial y}(Y^2(s-t_1)(y)-Z^2(s-t_1)(y))dy\\
=&\frac{1}{2}\int _{0}^{t}ds\int_0^1\frac{\partial p_{t-s}(\xi,y)}{\partial 
y}(Y^2(s)(y)-Z^2(s)(y))dy,
\endalign
$$
for $t > 0, \xi \in \Cal D$.
Using Cauchy-Schwartz inequality and (4.43), we have
$$
\align
&
\vert\vert Y(t)-Z(t)\vert\vert_{L^2([0,1])}^2=\int_0^1(Y(t)(\xi)-Z(t)(\xi))^2d\xi \\
\leq &C \int_0^1 d\xi(\frac{1}{2}\int 
_{0}^{t}ds\int_0^1\frac{1}{\sqrt{t-s}}\frac{c_1}{\sqrt{t-s}}e^{-\frac{(\xi-y)^2}{2c_
2(t-s)}}
\vert Y^2(s)(y)-Z^2(s)(y)\vert dy)^2\\
\leq & C \int_0^1 \int_0^t \frac{1}{(t-s)^{\frac{3}{4}}}ds\int
_{0}^{t}\frac{1}{(t-s)^{\frac{1}{4}}} \biggl 
(\int_0^1\frac{c_1}{\sqrt{t-s}}e^{-\frac{(\xi-y
)^2}{2c_2(t-s)}} \vert Y^2(s)(y)-Z^2(s)(y)\vert
dy \biggr )^2dsd\xi \\
\leq & C \int_0^1 \int 
_{0}^{t}\frac{1}{(t-s)^{\frac{1}{4}}}\int_0^1\frac{c_1}{\sqrt{t-s}}e^{-\frac{(\xi-y)
^2}{2c_2(t-s)}}(
Y(s)(y)-Z(s)(y))^2 dy\\
&\hskip1cm \times \int_0^1\frac{c_1}{\sqrt{t-s}}e^{-\frac{(\xi-y)^2}{2c_2(t-s)}}( 
Y(s)(y)+Z(s)(y))^2 dydsd\xi\\
\leq & C\sup_{0\leq s\leq T}[\vert\vert Y(s)+Z(s)\vert\vert_{L^2([0,1])}^2]  \int 
_{0}^{t}\frac{1}{(t-s)^{\frac{3}{4}}}\vert\vert
Y(s)-Z(s)\vert\vert_{L^2([0,1])}^2ds\\
\leq & M(\omega) \int _{0}^{t}\frac{1}{(t-s)^{\frac{3}{4}}}\vert\vert 
Y(s)-Z(s)\vert\vert_{L^2([0,1])}^2ds,
\tag4.49
\endalign
$$
for all $t \in (0,a]$, where $C$ is a generic constant that may change from line to 
line. Note that in the above computation, we have also used the fact that 
$$\sup_{0\leq s\leq a}[\vert\vert Y(s)+Z(s)\vert\vert_{L^2([0,1])}^2]<\infty, \quad 
a 
\in 
\R^+ .$$
As in the proof of (4.45) it follows from (4.49) and Gronwall's lemma   that 
$Y(t)=Z(t)$ for all $t \geq 0$. This completes the proof of assertion (ii) of the 
theorem.   

To prove assertion (iii), it is easy to see from the proof of (4.45) that
$$
\aligned
||v(t,&\psi_m, \o)-v(t,\psi_n,\o)||_{L^2([0,1])}^2\\
\leq &3||T_t\psi_m-T_t\psi_n||_{L^2([0,1])}^2+Ct\int _0^t \biggl ({1\over 
\sqrt{t-s}}+1 
\biggr )
||v(s,\psi_m,\o)-v(s,\psi_n,\o)||_{L^2([0,1])}^2ds\\
&+Ct^{1\over 4}\int _0^t \biggl ({1\over (t-s)^{3\over 4}}+{1\over (t-s)^{1\over 4}} 
\biggr 
)||v(s,\psi_m,\o)-v(s,\psi_n,\o)||_{L^2([0,1])}^2ds,
\endaligned
\tag4.50
$$
for $t \in [0,a]$ where $C$ is a positive random constant. Now using (4.50) and the 
same 
argument as in 
the proof of compactness of $Dv(t,\psi,\omega)$
in Theorem 4.1, one can show $v(t,\cdot, \o): L^2([0,1])\to L^2([0,1]), \, t > 0,$ 
takes 
bounded sets into relatively compact sets.
The only difference is that  we 
have to iterate (4.50) once before we can use Gronwall's lemma. Details of the 
proof are omitted.
\hfill \qed

\enddemo


%
%


\bigskip



\newpage

\centerline{\bf{ REFERENCES}}

\baselineskip=14truept
\parskip=3truept

\medskip

\medskip
\item {[B-C-J]\quad}
Bertini, L., Cancrin, N. and Jona-Lasinio, G.:
The stochastic
Burgers equation. {\it Commum. Math. Phys.}, 165 (1994), 211-232.
\medskip


\item{[B-C-F]\quad} Brzezniak, Z\., Capinski, M, and Flandoli, F\.,
Stochastic Navier-Stokes equations with multiplicative noise, {\it Stochastic 
Analysis 
and Applications}, 10 (1992), 523-532.

\medskip

\item{[B-F]\quad} Brzezniak, Z\., and Flandoli, F\., Regularity of
solutions and
random evolution operator for stochastic parabolic equations, In:
{\it Stochastic
partial differential equations and applications (Trento, 1990)},
pp54--71, Pitman Res. Notes Math. Ser., 268, Longman Sci. Tech., Harlow, 1992.

  \medskip

    \item{[B-F.1]\quad} Bensoussan, A\. and Flandoli, F\., Stochastic
inertial manifold.
Stochastics Stochastics Rep. 53 (1995), no. 1-2, 13--39.

    \medskip






\item{[D-T-Z]\quad} Davies, I. M., Truman, A., Zhao, H.Z., 
Stochastic generalized KPP equations, {\it Proc. R. Soc. 
Edinb.}, Vol. 126A (1996), No. 5, 
957-84.
\medskip


\item {[D-D-T]\quad}
Da Prato, G., Debusche, A. and Temam, R.: Stochastic Burgers
equation. {\it Nonlinear Differential Eq. Appl}, 1 (1994), 
389.
\medskip

\item{[D-Z.1]\quad} Da Prato, G., and Zabczyk, J., {\it 
Stochastic Equations in Infinite Dimensions}, Cambridge University 
Press (1992).

\medskip

\item{[D-Z.2]\quad} Da Prato, G., and Zabczyk, J., {\it 
Ergodicity for Infinite Dimensional Systems}, Cambridge University 
Press (1996).



\medskip

\item{[E-V]\quad} E, W\., Vanden Eijnden, E., 
Statistical theory for
the stochastic
Burgers equation in the 
inviscid limit, {\it Commun. Pure Appl.
Math.} 53(2000), 
852-901.

\medskip
\item{[E-H]\quad} Eckmann, J.-P., Hairer, M., 
Invariant measures for 
stochastic partial
differential equations in 
unbounded domains, {\it Nonlinearity}, Vol. 
14 (2001), 
133-151.
\medskip

\item{[E-Z]\quad} Elworthy,  K. D. and Zhao, H. 
Z., The propagation 
of travelling
waves for stochastic generalized 
KPP equations,
  {\it Mathl. and Comput. Modelling,} Vol. 20 (1994), 
No.4/5,
131-66.
\medskip

    \item{[F.1]\quad} Flandoli, F\.,  {\it 
Regularity theory and
stochastic flows for parabolic
SPDE's}, 
Stochastics Monographs, 9, Gordon and Breach Science
Publishers, 
Yverdon, 1995.

\medskip

    \item{[F.2]\quad} Flandoli, F\., 
Stochastic flows for nonlinear
second-order parabolic
SPDE, {\it Ann. 
Probab.} 24 (1996), no. 2, 547--558.

\medskip

   \item{[F-S]\quad} 
Flandoli, F\.,  and Schauml\"offel, K\.-U\.,
Stochastic parabolic equations in bounded domains: random evolution operator and 
Lyapunov exponents, {\it  Stochastics and Stochastics Reports} 29 (1990), no. 4, 
461--485.

\medskip
\item{[Fr]\quad} Freidlin, M., 
{\it Functional integration and 
partial differential
equations}, 
Princeton University Press, Princeton, New Jersey 
(1985).
\medskip

\item {[H-L-O-U-Z]\quad}
   Holden, H., Lindstrom, 
T., Oksendal, B., Uboe, J.
and Zhang, T. S.: The Burgers' equation 
with a noise force and the
stochastic heat  equations. {\it Comm. 
PDE.}, 19(1994),  119-141.
\medskip



\item{[I-W] \quad} Ikeda, N\., and 
Watanabe, S\., {\it Stochastic
Differential Equations and Diffusion 
Processes,\/} Second Edition,
North-Holland-Kodansha (1989).

\medskip




\item{[L-S-U] \quad} Ladyzenskaja, O. A., Solonnikov, V. A. and Uralceva, N. N.,
{\it  Linear and quasi-linear equations of parabolic type}, Translations of 
Mathematical 
Monographs, Vol. 23, American Mathematical Society (1968).

\medskip



 

\item{[Mo.1]\quad}  Mohammed, S\.-E\.A\., {\it 
Stochastic Functional
Differential Equations\/}, Research Notes in 
Mathematics, no. 99,
Pitman Advanced Publishing Program, 
Boston-London-Melbourne (1984).

\medskip






\item{[Mo.2]\quad} Mohammed, S\.-E\. A\., Non-Linear Flows 
for Linear
Stochastic
Delay Equations, {\it Stochastics, \/} Vol. l7 
\#3, (1987), 207--212.

\medskip


\item{[M-S.1]\quad}  Mohammed, S\.-E\. A\., and 
Scheutzow, M\. K\. R\., The stable manifold theorem for 
non-linear   stochastic systems with memory. Part I: Existence of the semiflow. Part 
II: 
The local stable manifold theorem" (preprints) (2001).

\medskip

\item{[M-S.2]\quad}  Mohammed, 
S\.-E\. A\., and  Scheutzow, M\. K\. R\., The stable manifold 
theorem for stochastic differential equations, {\it The Annals of 
Probability,\/} Vol. 27, No. 2, (1999), 615-652.
\medskip

\item{[M-Z-Z]\quad}  Mohammed, S\.-E\. A\., Zhang, 
T\. S\. and Zhao, H\. Z\., {\it The stable manifold theorem for semilinear stochastic 
evolution equations and stochastic partial
differential equations II: Existence of stable and unstable manifolds} (preprint) 
(2002).

 \medskip

\item{[O-V-Z]\quad} Oksendal, B., V\"age, G.  and 
Zhao, H. Z.,Two 
properties of stochastic KPP equations: 
ergodicity
and pathwise property, {\it Nonlinearity}, Vol. 14 (2001), 
639-662.
\medskip

\item{[O]\quad} Oseledec, V\. I\.,  A multiplicative 
ergodic theorem. Lyapunov
characteristic numbers for dynamical 
systems, {\it Trudy Moskov. Mat. Ob\v s\v
c.\/} 19 (1968), 179-210. 
English transl. {\it Trans. Moscow Math. Soc.\/} 19
(1968), 
197-221.
\medskip



\item{[Pr]\quad} 
Protter, Ph\. E\.,  {\it Stochastic Integration and
Stochastic 
Differential
Equations: A New Approach\/}, Springer 
(1990).

\medskip

\item{[Ro]\quad} Robinson, J. C., {\it Infinite-Dimensional Dynamical Systems. An 
Introduction to Dissipative 
Parabolic PDE's
and the Theory of Global Attractors}, Cambridge Texts in Applied Mathematics, 
Cambridge 
University Press (2001).

\medskip

\item{[Ru.1]\quad} Ruelle, D., Ergodic theory of 
differentiable dynamical
systems,
{\it Publ. Math.
Inst. Hautes Etud. 
Sci.} (1979), 275-306.

\medskip
\item{[Ru.2]\quad} Ruelle, D\., 
Characteristic exponents and
invariant manifolds
in Hilbert 
space,
{\it Annals of Math. 115\/} (1982), 
243--290.

\medskip

\item{[Si]\quad} Sinai, Ya. G., Burgers system driven by a periodic stochastic 
flow, In: {\it It\^o's stochastic calculus and probability theory},  Springer, Tokyo 
(1996), 347--353.

\medskip

\item{[Sk]\quad} 
Skorohod, A\. V\., {\it Random Linear Operators,\/} 
Riedel
(1984).

\medskip

\item{[Ta]\quad} Taylor, M. E., {\it Partial 
Differential Equations III Nonlinear
Equations}, Springer 
(1996).

\medskip

\item{[Te]\quad} Temam, R., {\it Infinite-Dimensional Dynamical Systems in Mechanics 
and 
Physics}, 
 Springer-Verlag (1988).

\medskip

    \item{[Tw]\quad} Twardowska, K., An extension 
of the Wong-Zakai theorem for
stochastic equations in Hilbert 
spaces, {\it Stochastic Anal. Appl.} 10 (1992), no. 
4,
471--500.
\medskip

\item{[T-Za]\quad} Tribe, R., Zaboronski, O., 
On the large time
asymptotics of decaying Burgers turbulence.,
{\it Comm. Math. Phys.}, 212 (2000), 415--436.
\medskip

\item{[T-Z]\quad} 
Truman, A., Zhao, H. Z., Stochastic Burgers'
equations and 
their
semi-classical expansions, {\it Comm. Math. Phys.}, 
194 (1998), 231-248.
\medskip

\medskip

\bigskip

\noindent

\bigskip

\noindent
Salah-Eldin A. Mohammed \newline
Department of Mathematics,
\newline
Southern Illinois University at Carbondale,
\newline
Carbondale, Illinois 62901.\newline
Email: salah\@sfde.math.siu.edu \newline
Web page: http://sfde.math.siu.edu

\bigskip
\medskip
\noindent
Tusheng Zhang \newline
Department of Mathematics \newline
     University of Manchester,\newline
Oxford Road, Manchester M13 9PL\\ \newline
UK.
     \newline
Email: tzhang\@math.man.ac.uk

\bigskip
\medskip
\noindent
Huaizhong Zhao \newline
Department of Mathematical Sciences \newline
Loughborough University,\newline
LE11 3TU,\newline
UK.
     \newline
Email: H.Zhao\@lboro.ac.uk

\end{document}